\documentclass[hidelinks, 11 pt, a4paper]{article}
\usepackage{amsfonts}
\usepackage{amsmath}
\usepackage{amssymb}
\usepackage{amsthm}
\usepackage{subfigure}
\usepackage[pdftex]{graphicx}
\usepackage{natbib}
\usepackage[hypertexnames=false]{hyperref}
\usepackage{setspace}
\usepackage{pdfsync}
\usepackage{lscape}
\usepackage[hmargin=1.25in,vmargin=1.25in]{geometry}

\hypersetup{colorlinks, linkcolor = {teal}, citecolor = {blue}, urlcolor ={blue!80!black}}

\makeatletter
\renewcommand\paragraph{\@startsection{paragraph}{4}{\z@}%
            {-2.5ex\@plus -1ex \@minus -.25ex}%
            {1.25ex \@plus .25ex}%
            {\normalfont\normalsize\bfseries}}
\makeatother
\usepackage{mathrsfs,dsfont}
\usepackage[dvipsnames,svgnames, x11names]{xcolor}
\usepackage{pgfplots}
\usetikzlibrary{patterns}
\usepackage{caption}

\DeclareGraphicsExtensions{.jpg}

\def\qed{\rule{2mm}{2mm}}

\newtheorem{proposition}{Proposition}[section]
\newtheorem{theorem}{Theorem}[section]
\newtheorem{lemma}{Lemma}[section]
\newtheorem{definition}{Definition}[section]
\newtheorem{example}{Example}[section]
\newtheorem{corollary}{Corollary}[section]
\newtheorem{remark}{Remark}[section]
\newtheorem{assumption}{Assumption}[section]

\newenvironment{packed_enum}{\begin{enumerate} \setlength{\itemsep}{1pt}\setlength{\parskip}{0pt}\setlength{\parsep}{0pt}}{\end{enumerate}}

\begin{document}
\title{Inference on Directionally Differentiable Functions\thanks{We would like to thank Brendan Beare, Xiaohong Chen, Victor Chernozhukov, Bruce Hansen, Han Hong, Hiroaki Kaido, and numerous seminar participants for comments that helped improve this paper.}}
{\author{
Zheng Fang \\ Department of Economics \\ Kansas State University \\ zfang@ksu.edu
\and
Andres Santos\thanks{Research supported by NSF Grant SES-1426882.}\\ Department of Economics\\ U.C. San Diego\\ a2santos@ucsd.edu}}

\date{First Draft: February, 2014\\ This Draft: December, 2015}
\maketitle

\begin{abstract}
This paper studies an asymptotic framework for conducting inference on parameters of the form $\phi(\theta_0)$, where $\phi$ is a known directionally differentiable function and $\theta_0$ is estimated by $\hat \theta_n$. In these settings, the asymptotic distribution of the plug-in estimator $\phi(\hat \theta_n)$ can be readily derived employing existing extensions to the Delta method. We show, however, that (full) differentiability of $\phi$ is a necessary and sufficient condition for bootstrap consistency whenever the limiting distribution of $\hat \theta_n$ is Gaussian. An alternative resampling scheme is proposed which remains consistent when the bootstrap fails, and is shown to provide local size control under restrictions on the directional derivative of $\phi$. We illustrate the utility of our results by developing a test of whether a Hilbert space valued parameter belongs to a convex set -- a setting that includes moment inequality problems, tests of random utility models, and certain tests of shape restrictions as special cases (e.g. tests of monotonicity of the pricing kernel or of parametric conditional quantile model specifications).
\end{abstract}

\begin{center}
\end{center}

\newpage

\section{Introduction}

The Delta method is a cornerstone of asymptotic analysis, allowing researchers to easily derive asymptotic distributions, compute standard errors, and establish bootstrap consistency.\footnote{Interestingly, despite its importance, the origins of the Delta method remain obscure. \cite{hoef2012} recently attributed its invention to the economist Robert Dorfman in his article \cite{dorfman1938}, which was curiously published by the Worcester State Hospital (a public asylum for the insane).} However, an important class of estimation and inference problems in economics fall outside its scope. These problems study parameters of the form $\phi(\theta_0)$, where $\theta_0$ is unknown but estimable and $\phi$ is a known but potentially non-differentiable function. Such a setting arises frequently in economics, with applications including the construction of parameter confidence regions in moment inequality models \citep{Pakes_Porter_etal2006aWP, ciliberto:tamer:2009}, the study of convex partially identified sets \citep{beresteanu:molinari:2008, Bontemps_Magnac_Maurin2012}, and the development of tests of superior predictive ability \citep{white:2000, hansen:2005}, of stochastic dominance \citep{Linton2010}, and of likelihood ratio ordering \citep{BeareandMoon2012}.

The aforementioned examples share a structure common to numerous ``nonstandard" inference problems in economics: the transformation $\phi$ is directionally (but not fully) differentiable in a local neighborhood of $\theta_0$. In this paper, we show this common structure enables us to reduce challenging statistical questions to simple analytical considerations regarding the directional derivative of $\phi$ -- much in the same manner the Delta method and its bootstrap counterpart fundamentally simplify the analysis of applications in which $\phi$ is differentiable. Concretely, we examine a setting in which $\theta_0$ is a possibly infinite dimensional parameter and there exists an estimator $\hat \theta_n$ whose asymptotic distribution we denote by $\mathbb G_0$ -- i.e., for some sequence $r_n \uparrow \infty$, we have that
\begin{equation}\label{intro1}
r_n\{\hat \theta_n - \theta_0\} \stackrel{L}{\rightarrow} \mathbb G_0 ~.
\end{equation}
Within this framework, we study a simple unifying approach for conducting inference on the parameter $\phi(\theta_0)$ by employing $\phi(\hat \theta_n)$ and a suitable estimator of its asymptotic distribution -- a practice common in, for example, the study of moment inequality \citep{andrews:soares:2010}, conditional moment inequality \citep{andrews:shi:2013}, and incomplete linear models \citep{beresteanu:molinari:2008}.

As has been previously noted in the literature, the traditional Delta method readily generalizes to the case where $\phi$ is directionally differentiable \citep{Shapiro1991, Dumbgen1993}. In particular, if $\phi$ is Hadamard directionally differentiable, then
\begin{equation}\label{intro2}
r_n\{\phi(\hat \theta_n) - \phi(\theta_0)\} \stackrel{L}{\rightarrow} \phi_{\theta_0}^\prime(\mathbb G_0) ~,
\end{equation}
where $\phi_{\theta_0}^\prime$ denotes the directional derivative of $\phi$ at $\theta_0$. The utility of the asymptotic distribution of $\phi(\hat \theta_n)$, however, hinges on our ability to consistently estimate it. While it is tempting in these problems to resort to resampling schemes such as the bootstrap of \cite{efron1979}, we know by way of example that they may be inconsistent even if they are valid for the original estimator $\hat \theta_n$ \citep{Bickel_Gotze_Zwet1997, Andrews2000Bootstrap, ham:woutersen:2013}. In our first main result, we establish that these examples reflect a deeper underlying principle. Specifically, we establish that whenever the asymptotic distribution of $\hat \theta_n$ is Gaussian, full differentiability of $\phi$ at $\theta_0$ is in fact a \emph{necessary and sufficient} condition for the consistency of ``standard" bootstrap methods. As a result, we obtain a purely analytical diagnostic for assessing bootstrap consistency in these settings: one need only verify whether $\phi$ is differentiable. An important consequence of our characterization of bootstrap consistency is that, in our setting, ``standard" bootstrap methods in fact fail whenever the asymptotic distribution of $\phi(\hat \theta_n)$ is not Gaussian -- a conclusion that yields an alternative simple way to detect the failure of the bootstrap.

Intuitively, consistently estimating the asymptotic distribution of $\phi(\hat \theta_n)$ requires us to adequately approximate both the law of $\mathbb G_0$ and the directional derivative $\phi_{\theta_0}^\prime$ (see \eqref{intro2}). While a consistent bootstrap procedure for $\hat \theta_n$ enables us to do the former, the bootstrap fails for $\phi(\hat \theta_n)$ due to its inability to properly estimate $\phi_{\theta_0}^\prime$. These heuristics, however, readily suggest a remedy to the problem -- namely to compose a suitable estimator $\hat \phi_n^\prime$ for $\phi_{\theta_0}^\prime$ with the bootstrap approximation to the asymptotic distribution of $\hat \theta_n$. We formalize this intuition, and provide conditions on $\hat \phi_n^\prime$ that ensure the proposed approach yields consistent estimators of the asymptotic distribution of $\phi(\hat \theta_n)$ and its quantiles. Moreover, we further show that popular existing inferential procedures developed in the context of specific applications in fact follow precisely this approach -- these include \cite{andrews:soares:2010} for moment inequalities, \cite{Linton2010} for tests of stochastic dominance, and \cite{Kaido2013dual} for convex partially identified models. Consequently, our proposed resampling scheme may be interpreted as a generalization of procedures that have become prevalent within their respective literatures.

As argued by \cite{imbens:manski:2004}, pointwise asymptotic approximations may be unreliable, in particular when $\phi(\hat \theta_n)$ is not regular. Heuristically, if the asymptotic distribution of $\phi(\hat \theta_n)$ is sensitive to local perturbations of the data generating process, then employing \eqref{intro2} as the basis for inference may yield poor size control in finite samples. We thus examine the ability of our proposed procedure to provide local size control in the context of employing $\phi(\hat \theta_n)$ as a test statistic for the hypothesis
\begin{equation}\label{intro3}
H_0 : \phi(\theta_0) \leq 0 \hspace{1 in } H_1 : \phi(\theta_0) > 0 ~.
\end{equation}
Special cases of \eqref{intro3} include inference in moment inequality models and tests of stochastic dominance -- instances in which our framework encompasses procedures that provide local, in fact uniform, size control \citep{andrews:soares:2010, Linton2010, andrews:shi:2013}. We show that the common structure linking these applications is that $\phi_{\theta_0}^\prime$ and $\hat \theta_n$ are respectively subadditive and regular. Indeed, we more generally establish that these two properties suffice for guaranteeing the ability of our procedure to locally control size along contiguous perturbations. Thus, our analysis again reduces a challenging statistical problem (establishing local size control) to a property of the original estimator (regularity) and an analytical calculation concerning the directional derivative (verifying subadditivity). As part of this local analysis, we in addition characterize local power and conclude that, under mild regularity conditions, the bootstrap is valid for $\phi(\hat \theta_n)$ if and only if $\phi(\hat \theta_n)$ is regular.

We illustrate the utility of our analysis by developing a common testing framework for a diverse set of hypotheses that includes whether the pricing kernel is monotone, whether parametric conditional quantile specifications satisfy shape restrictions, and whether a vector of means belongs to a closed convex set -- see, e.g., \cite{wolak:1988} and \cite{kitamura:stoye:2013} for special cases of the latter hypothesis. Formally, these superficially disparate examples are in fact illustrations of the general problem of testing whether a Hilbert space valued parameter $\theta_0$ belongs to a known closed convex set $\Lambda$. By exploiting the directional differentiability of projections onto convex sets \citep{zarantonello}, we can easily study the properties of a test for this general hypothesis that is based on the distance between the estimator $\hat \theta_n$ and the set $\Lambda$. For instance, our analysis readily implies that the bootstrap is inconsistent, but we are nonetheless able to obtain valid critical values by constructing a suitable estimator $\hat \phi_n^\prime$ which we compose with a bootstrap approximation to the limiting distribution of $\hat \theta_n$. In addition, we establish that the directional derivative $\phi_{\theta_0}^\prime$ is always subadditive, and thus conclude that the proposed test is able to locally control size provided $\hat \theta_n$ is regular. A brief simulation study confirms our theoretical findings by showing that the proposed test indeed exhibits good finite sample size control.

In related work, an extensive literature has established the consistency of the bootstrap and its ability to provide a refinement when $\theta_0$ is a vector of means and $\phi$ is a differentiable function \citep{HallBoot, HorowitzBoot}. Our analysis is most closely related to the pioneering work of \cite{Dumbgen1993}, who first examined the validity of the bootstrap for estimating the asymptotic distribution of $\phi(\hat \theta_n)$ under a potential lack of differentiability. The results in \cite{Dumbgen1993} imply a characterization of bootstrap consistency that, unlike ours, applies when $\mathbb G_0$ is not Gaussian but is harder for practitioners to verify as it concerns properties of both $\mathbb G_0$ and $\phi_{\theta_0}^\prime$.\footnote{We discuss the relationship between these conditions in detail in Section \ref{sec:bootfail}; see Remark \ref{rm:dumbgen}.} In more recent studies, applications where $\phi$ is not fully differentiable have garnered increasing attention due to their preponderance in the analysis of partially identified models \citep{manski:2003}. \cite{hirano:porter:2009}, \cite{Song2013Minimax}, and \cite{fang:2014}, for example, explicitly exploit the directional differentiability of $\phi$ as well, though their focus is on estimation rather than inference. Other work studying these irregular models, though not explicitly relying on the directional differentiability of $\phi$, include \cite{chernozhukov:hong:tamer:2007, chernozhukov:lee:rosen:2013}, \cite{romano:shaikh:2008, Romano_Shaikh2010}, \cite{bugni:2010}, and \cite{canay2010} among many others.

An emerging body of research has validated the usefulness of our results by both employing and expanding on them. For instance, \cite{seo2014tests} and \cite{BrendanPrelim} have used our framework to develop tests of stochastic monotonicity and of density ratio ordering respectively.  Other applications of our results also include \cite{hansen2015regression} who studied the asymptotic properties of regression kink models, \cite{jha2015testing} who estimated transaction costs in energy future markets, and \cite{lee:bhattacharya} who proposed methods for estimation welfare changes in partially identified discrete choice models. Finally, in a highly complementary paper, \cite{hong:li:2014} have built on our results and established that a consistent estimator of $\hat \phi_n^\prime$ for the directional derivative $\phi_{\theta_0}^\prime$ can often be obtained through numerical differentiation of $\phi$ at $\hat \theta_n$ (under appropriate conditions on the step size).

The remainder of the paper is organized as follows. Section 2 formally introduces the model we study and contains a minor extension of the Delta method for directionally differentiable functions. In Section 3 we characterize necessary and sufficient conditions for bootstrap consistency, develop an alternative method for estimating the asymptotic distribution of $\phi(\hat \theta_n)$, and study the local properties of this approach. Section 4 applies these results to develop a test of whether a Hilbert space valued parameter belongs to a closed convex set. All proofs are contained in the Appendix.

\section{Setup and Background}\label{Sec: General Setup}

In this section, we introduce our notation and review the concepts of Hadamard and directional Hadamard differentiability as well as their implications for the Delta method.

\subsection{General Setup}

In order to accommodate applications such as conditional moment inequalities and tests of shape restrictions, we must allow for both the parameter $\theta_0$ and the map $\phi$ to take values in possibly infinite dimensional spaces; see Examples \ref{ex:condineq}-\ref{ex:lrorder} below. We therefore impose the general requirement that $\theta_0 \in \mathbb D_\phi$ and $\phi : \mathbb D_\phi \subseteq \mathbb D \rightarrow \mathbb E$ for $\mathbb D$ and $\mathbb E$ Banach spaces with norms $\|\cdot\|_{\mathbb D}$ and $\|\cdot\|_{\mathbb E}$ respectively, and $\mathbb D_\phi$ the domain of $\phi$.

The estimator $\hat \theta_n$ is assumed to be a function of a sequence of random variables $\{X_i\}_{i=1}^n$ into the domain of $\phi$. The distributional convergence
\begin{equation}\label{setup1}
r_n\{\hat \theta_n - \theta_0\} \stackrel{L}{\rightarrow} \mathbb G_0 ~,
\end{equation}
is then understood to be in $\mathbb D$ and with respect to the joint law of $\{X_i\}_{i=1}^n$. For instance, if $\{X_i\}_{i=1}^n$ is an i.i.d. sample and each $X_i \in \mathbf R^d$ is distributed according to $P$, then probability statements for $\hat \theta_n :\{X_i\}_{i=1}^n \rightarrow \mathbb D_\phi$ are understood to be with respect to the product measure $\bigotimes_{i=1}^n P$. We emphasize, however, that our results are applicable to dependent settings as well. In addition, we also note the convergence in distribution in \eqref{setup1} is meant in the Hoffman-J{\o}rgensen sense \citep{Vaart1996}. Expectations throughout the text should therefore be interpreted as outer expectations, though we obviate the distinction in the notation. The notation is made explicit in the Appendix whenever differentiating between inner and outer expectations is necessary.

Finally, we introduce notation that is recurrent in the context of our examples. For a set $\mathbf A$, we denote the space of bounded functions on $\mathbf A$ by
\begin{equation}\label{setup2}
\ell^\infty(\mathbf A) \equiv \{f:\mathbf A \rightarrow \mathbf R \text{ such that } \|f\|_\infty < \infty\} \hspace{0.5 in} \|f\|_\infty \equiv \sup_{a\in \mathbf A} |f(a)| ~,
\end{equation}
and note $\ell^\infty(\mathbf A)$ is a Banach space under $\|\cdot\|_\infty$. If in addition $\mathbf A$ is a compact Hausdorff topological space, then we let $\mathcal C(\mathbf A)$ denote the set of continuous functions on $\mathbf A$,
\begin{equation}\label{setup3}
\mathcal C(\mathbf A) \equiv \{f : \mathbf A \rightarrow \mathbf R \text{ such that } f \text{ is continuous }\} ~,
\end{equation}
which satisfies $\mathcal C(\mathbf A) \subset \ell^\infty(\mathbf A)$ and is also a Banach space when endowed with $\|\cdot\|_\infty$.

\subsubsection{Examples}

In order to fix ideas, we next introduce a series of examples that illustrate the broad applicability of our setting. We return to these examples throughout the paper, and develop a formal treatment of each of them in the Appendix. For ease of exposition, we base our discussion on simplifications of well known models, though we note that our results apply to the more general problems that motivated them.

Our first example is due to \cite{Bickel_Gotze_Zwet1997}, and provides an early illustration of the potential failure of the nonparametric bootstrap.

\begin{example}[Absolute Value of Mean]\label{ex:absmean} \rm
Let $X \in \mathbf R$ be a scalar valued random variable, and suppose we wish to estimate the parameter
\begin{equation}\label{ex:absmean1}
\phi(\theta_0) = |E[X]| ~.
\end{equation}
Here, $\theta_0 = E[X]$, $\mathbb D = \mathbb E = \mathbf R$, and $\phi : \mathbf R\rightarrow \mathbf R$ satisfies $\phi(\theta) = |\theta|$ for all $\theta \in \mathbf R$. \qed
\end{example}

Our next example is a special case of the intersection bounds model studied in \cite{hirano:porter:2009} and \cite{chernozhukov:lee:rosen:2013} among many others.

\begin{example}[Intersection Bounds]\label{ex:interbounds} \rm
Let $X = (X^{(1)},X^{(2)})^\prime \in \mathbf R^2$ be a bivariate random variable, and consider the problem of estimating the parameter
\begin{equation}\label{ex:interbounds1}
\phi(\theta_0) = \max\{E[X^{(1)}], E[X^{(2)}]\} ~.
\end{equation}
In this context, $\theta_0 = (E[X^{(1)}],E[X^{(2)}])^\prime$, $\mathbb D = \mathbf R^2$, $\mathbb E = \mathbf R$, and $\phi : \mathbf R^2 \rightarrow \mathbf R$ is given by $\phi(\theta) = \max\{\theta^{(1)},\theta^{(2)}\}$ for any $(\theta^{(1)},\theta^{(2)})^\prime = \theta\in \mathbf R^2$. Functionals such as \eqref{ex:interbounds1} are also often employed for inference in moment inequality models; see \cite{chernozhukov:hong:tamer:2007}, \cite{romano:shaikh:2008}, and \cite{andrews:soares:2010}. \qed
\end{example}

A related example arises in conditional moment inequality models, as studied in \cite{andrews:shi:2013}, \cite{armstrong:chan:2012}, and \cite{chetverikov:2012}.

\begin{example}[Conditional Moment Inequalities]\label{ex:condineq} \rm
Let $X = (Y,Z^\prime)^\prime$ with $Y \in \mathbf R$ and $Z \in \mathbf R^{d_z}$. For a suitable set of functions $\mathcal F \subset \ell^\infty(\mathbf R^{d_z})$, \cite{andrews:shi:2013} propose testing whether $E[Y|Z] \leq 0$ almost surely, by estimating the parameter
\begin{equation}\label{ex:condineq1}
\phi(\theta_0) = \sup_{f \in \mathcal F} E[Yf(Z)] ~.
\end{equation}
Here, $\theta_0 \in \ell^\infty(\mathcal F)$ satisfies $\theta_0(f) = E[Yf(Z)]$ for all $f \in \mathcal F$, $\mathbb D = \ell^\infty(\mathcal F)$, $\mathbb E = \mathbf R$, and the map $\phi:\mathbb D \rightarrow \mathbb E$ is given by $\phi (\theta) = \sup_{f\in \mathcal F} \theta(f)$. \qed
\end{example}

The following example is an abstract version of an approach pursued in \cite{beresteanu:molinari:2008} and \cite{Bontemps_Magnac_Maurin2012} for studying partially identified models.

\begin{example}[Convex Identified Sets]\label{ex:suppfunction} \rm
Let $\Lambda \subseteq \mathbf R^d$ denote a convex and compact set, $\mathbb S^{d}$ be the unit sphere on $\mathbf R^d$ and $\mathcal C(\mathbb S^d)$ denote the space of continuous functions on $\mathbb S^d$. For each $p\in \mathbb S^{d}$, the support function $\nu(\cdot,\Lambda)\in \mathcal C(\mathbb S^{d})$ of the set $\Lambda$ is then given by
\begin{equation}\label{ex:suppfunction1}
\nu(p,\Lambda) \equiv \sup_{\lambda \in \Lambda} \langle p,\lambda \rangle \quad \quad p \in \mathbb S^d ~.
\end{equation}
As noted by \cite{beresteanu:molinari:2008} and \cite{Bontemps_Magnac_Maurin2012}, the functional
\begin{equation}\label{ex:suppfunction2}
\phi(\theta_0) = \sup_{p\in \mathbb S^d} \{\langle p, \lambda\rangle - \nu(p,\Lambda)\} ~,
\end{equation}
can form the basis for a test of whether $\lambda$ is an element of $\Lambda$, since $\lambda \in \Lambda$ if and only if $\phi(\theta_0) \leq 0$. In the context of this example, $\theta_0 = \nu(\cdot,\Lambda)$, $\mathbb D = \mathcal C(\mathbb S^d)$, $\mathbb E = \mathbf R$, and $\phi(\theta) = \sup_{p\in \mathbb S^d} \{\langle p,\lambda\rangle - \theta(p)\}$ for any $\theta\in \mathcal C(\mathbb S^d)$. \qed
\end{example}

Our next example is based on the \citet{Linton2010} test for stochastic dominance.

\begin{example}[Stochastic Dominance]\label{ex:stochdom} \rm
Let $X = (X^{(1)},X^{(2)})^\prime \in \mathbf R^2$ be continuously distributed, and define the marginal cdfs $F_j(u) \equiv P(X^{(j)} \leq u)$ for $j \in \{1,2\}$. For a positive integrable weighting function $w:\mathbf R\rightarrow \mathbf R_+$, \citet{Linton2010} estimate
\begin{equation}\label{ex:stochdom1}
\phi(\theta_0) = \int_{\mathbf R} \max\{F_1(u) - F_2(u),0\}w(u)du ~,
\end{equation}
to construct a test of whether $X^{(1)}$ first order stochastically dominates $X^{(2)}$. In this example, we set $\theta_0 = (F_1,F_2)$, $\mathbb D = \ell^\infty(\mathbf R)\times \ell^\infty(\mathbf R)$, $\mathbb E = \mathbf R$, and $\phi((\theta^{(1)},\theta^{(2)})) = \int \max\{\theta^{(1)}(u) - \theta^{(2)}(u),0\}w(u)du$ for any $(\theta^{(1)},\theta^{(2)}) \in \ell^\infty(\mathbf R)\times \ell^\infty(\mathbf R)$. \qed
\end{example}

In addition to tests of stochastic dominance, a more recent literature has aimed to examine whether likelihood ratios are monotonic. Our final example is a simplification of a test proposed in \cite{carolan:tebbs:2005} and \cite{BeareandMoon2012}.

\begin{example}[Likelihood Ratio Ordering]\label{ex:lrorder} \rm
Let $X = (X^{(1)},X^{(2)})^\prime \in \mathbf R^2$ have strictly increasing marginal cdfs $F_j(u) \equiv P(X^{(j)} \leq u)$, and define $G \equiv F_1 \circ F_2^{-1}$. Further let $\mathcal M : \ell^\infty([0,1]) \rightarrow \ell^\infty([0,1])$ be the least concave majorant operator, given by
\begin{equation}\label{ex:lrorder1}
\mathcal M f(u) = \inf \{g(u) : g \in \ell^\infty([0,1]) \text{ is concave and } f(u) \leq g(u) \text{ for all } u\in[0,1]\}
\end{equation}
for every $f \in \ell^\infty([0,1])$. Since the likelihood ratio $dF_1/dF_2$ is nonincreasing if and only if $G$ is concave on $[0,1]$ \citep{carolan:tebbs:2005}, \cite{BeareandMoon2012} note
\begin{equation}\label{ex:lrorder2}
\phi(\theta_0) = \Big\{\int_0^1 (\mathcal M G(u) - G(u))^2du\Big\}^{\frac{1}{2}}
\end{equation}
characterizes whether $dF_1/dF_2$ is nonincreasing because $\phi(\theta_0) = 0$ if and only if $G$ is concave.  In this example, $\theta_0 = G$, $\mathbb D = \ell^\infty([0,1])$, $\mathbb E = \mathbf R$, and $\phi:\mathbb D \rightarrow \mathbb E$ satisfies $\phi(\theta) = \{\int_0^1 (\mathcal M \theta(u) - \theta(u))^2du\}^{\frac{1}{2}}$ for any $\theta \in \ell^\infty([0,1])$. \qed
\end{example}

\subsection{Differentiability Concepts}\label{sec:diffconcept}

In all the previous examples, there exist points $\theta \in \mathbb D$ at which the map $\phi : \mathbb D \rightarrow \mathbb E$ is not differentiable. Nonetheless, at all such $\theta$ at which differentiability is lost, $\phi$ actually remains directionally differentiable. This is most easily seen in Examples \ref{ex:absmean} and \ref{ex:interbounds}, in which the domain of $\phi$ is a finite dimensional space. In order to address Examples \ref{ex:condineq}-\ref{ex:lrorder}, however, a notion of directional differentiability that is suitable for more abstract spaces $\mathbb D$ is necessary. Towards this end, we follow \cite{Shapiro1990} and define

\begin{definition}\label{def:had}\rm
Let $\mathbb D$ and $\mathbb E$ be Banach spaces, and $\phi:\mathbb D_\phi\subseteq \mathbb D\to\mathbb E$.
\begin{itemize}
    \item[(i)] The map $\phi$ is said to be {\it Hadamard differentiable} at $\theta \in\mathbb D_\phi$ {\it tangentially} to a set $\mathbb D_0\subseteq\mathbb D$, if there is a continuous linear map $\phi_\theta':\mathbb D_0\to\mathbb E$ such that
    \begin{equation}\label{def:had1}
    \lim_{n\rightarrow \infty}\| \frac{\phi(\theta +t_nh_n)-\phi(\theta)}{t_n} - \phi_\theta'(h)\|_{\mathbb E} = 0 ~,
    \end{equation}
    for all sequences $\{h_n\}\subset\mathbb D$ and $\{t_n\}\subset\mathbf R$ such that $t_n\to 0$, $h_n\to h\in\mathbb D_0$ as $n\to\infty$ and $\theta+t_nh_n\in\mathbb D_\phi$ for all $n$.
    \item[(ii)] The map $\phi$ is said to be {\it Hadamard directionally differentiable} at $\theta \in\mathbb D_\phi$ {\it tangentially} to a set $\mathbb D_0\subseteq\mathbb D$, if there is a continuous map $\phi_\theta':\mathbb D_0\to\mathbb E$ such that
    \begin{equation}\label{def:had2}
    \lim_{n\rightarrow \infty}\|\frac{\phi(\theta +t_n h_n)-\phi(\theta)}{t_n} -\phi_\theta'(h)\|_{\mathbb E} = 0 ~,
    \end{equation}
    for all sequences $\{h_n\}\subset\mathbb D$ and $\{t_n\}\subset\mathbf R_+$ such that $t_n\downarrow 0$, $h_n\to h\in\mathbb D_0$ as $n\to\infty$ and $\theta+t_nh_n\in\mathbb D_\phi$ for all $n$.
\end{itemize}
\end{definition}

As has been extensively noted in the literature, Hadamard differentiability is particularly suited for generalizing the Delta method to metric spaces \citep{reeds1976definition, gill1989non}. It is therefore natural to employ an analogous approximation requirement when considering an appropriate definition of a directional derivative (compare \eqref{def:had1} and \eqref{def:had2}). However, despite this similarity, two key differences distinguish (full) Hadamard differentiability from Hadamard directional differentiability. First, in \eqref{def:had2} the sequence of scalars $\{t_n\}$ must approach $0$ ``from the right", heuristically giving the derivative a direction. Second, the map $\phi^\prime_\theta: \mathbb D_0 \rightarrow \mathbb E$ is no longer required to be linear, though it is possible to show \eqref{def:had2} implies $\phi_\theta^\prime$ must be homogenous of degree one. It is in fact this latter property that distinguishes the two differentiability concepts.

\begin{proposition}\label{pro:equivalence}
Let $\mathbb D$, $\mathbb E$ be Banach spaces, $\mathbb D_0 \subseteq \mathbb D$ be a subspace, and $\phi: \mathbb D_\phi\subseteq\mathbb D\to\mathbb E$. Then, $\phi$ is Hadamard directionally differentiable at $\theta \in\mathbb D_\phi$ tangentially to $\mathbb D_0$ with linear derivative $\phi_\theta^\prime :\mathbb D_0 \rightarrow \mathbb E$ iff $\phi$ is Hadamard differentiable at $\theta$ tangentially to $\mathbb D_0$.
\end{proposition}

Thus, while (full) Hadamard differentiability implies Hadamard directional differentiability, Proposition \ref{pro:equivalence} shows the converse is true if the directional derivative $\phi_\theta^\prime$ is linear. In what follows, we mildly extend existing results that establish linearity of the directional derivative is in fact not important for the validity of the Delta Method \citep{Shapiro1991, Dumbgen1993}. As we will additionally show, however, linearity will play an instrumental role in determining whether the bootstrap is consistent or not.

\begin{remark}\label{rm:contcomment} \rm
Whenever equation \eqref{def:had2} is satisfied, the continuity of the corresponding map $\phi_{\theta}^\prime : \mathbb D_0 \rightarrow \mathbb E$ is in fact implied by $\mathbb D$ being a Banach space; see Proposition 3.1 in \cite{Shapiro1990}. Therefore, showing requirement \eqref{def:had2} holds for some map $\phi^\prime_\theta : \mathbb D_0\rightarrow \mathbb E$ suffices for establishing the Hadamard directional differentiability of $\phi$ at $\theta$. \qed
\end{remark}

\subsubsection{Examples Revisited}

We next revisit the examples to illustrate the computation of the directional derivative. The first two examples are straightforward, since the domain of $\phi$ is finite dimensional.

\noindent {\bf Example \ref{ex:absmean} (cont.)} The Hadamard directional derivative $\phi_\theta^\prime : \mathbf R \rightarrow \mathbf R$ here equals
\begin{equation}\label{ex:absmean2}
\phi_\theta^\prime (h) = \begin{cases}
h   &\text{ if } \theta > 0\\
|h| &\text{ if } \theta = 0\\
-h  &\text{ if } \theta < 0
\end{cases} ~.
\end{equation}
Note that $\phi$ is actually (fully) Hadamard differentiable everywhere except at $\theta = 0$, but that it is still Hadamard directionally differentiable at that point.  \qed

\noindent {\bf Example \ref{ex:interbounds} (cont.)} For $\theta = (\theta^{(1)},\theta^{(2)})^\prime \in \mathbf R^2$, let $j^* = \arg\max_{j\in\{1,2\}} \theta^{(j)}$. For any $h = (h^{(1)},h^{(2)})^\prime \in \mathbf R^2$, it is then straightforward to verify $\phi_\theta^\prime : \mathbf R^2 \rightarrow \mathbf R$ is given by
\begin{equation}\label{ex:interbounds2}
\phi_\theta^\prime (h) = \begin{cases}
h^{(j^*)}   & \text{ if } \theta^{(1)} \neq \theta^{(2)} \\
\max\{h^{(1)},h^{(2)}\}  &\text{ if } \theta^{(1)} = \theta^{(2)}
\end{cases} ~.
\end{equation}
As in \eqref{ex:absmean2}, $\phi_\theta^\prime$ is nonlinear precisely when Hadamard differentiability is not satisfied. \qed

In the next examples the domain of $\phi$ is infinite dimensional, and we sometimes need to employ Hadamard directional tangential differentiability -- i.e. $\mathbb D_0 \neq \mathbb D$.

\noindent {\bf Example \ref{ex:condineq} (cont.)} Suppose $E[Y^2] < \infty$ and that $\mathcal F$ is compact when endowed with the metric $\|f\|_{L^2(Z)} \equiv \{E[f(Z)^2]\}^{\frac{1}{2}}$. Then, $\theta_0 \in \mathcal C(\mathcal F)$, and Lemma \ref{lm:auxmax} in the Appendix implies $\phi$ is Hadamard directionally differentiable tangentially to $\mathcal C(\mathcal F)$ at any $\theta \in \mathcal C(\mathcal F)$. In particular, for $\Psi_{\mathcal F}(\theta) \equiv \arg\max_{f\in \mathcal F} \theta(f)$, the directional derivative is
\begin{equation}\label{ex:condineq2}
\phi_\theta^\prime(h) = \sup_{f \in \Psi_{\mathcal F}(\theta)} h(f) ~.
\end{equation}
Interestingly $\phi_\theta^\prime$ is linear at any $\theta \in \mathcal C(\mathcal F)$ for which $\Psi_{\mathcal F}(\theta)$ is a singleton, and hence $\phi$ is actually Hadamard differentiable at such $\theta$. We note in this example, $\mathbb D_0 = \mathcal C(\mathcal F)$. \qed

\noindent {\bf Example \ref{ex:suppfunction} (cont.)} For any $\theta \in \mathcal C(\mathbb S^d)$ let $\Psi_{\mathbb S^d}(\theta)\equiv \arg\max_{p \in \mathbb S^d} \{\langle p,\lambda\rangle - \theta(p)\}$. Lemma B.8 in \cite{Kaido2013dual} then shows that $\phi_\theta^\prime : \mathcal C(\mathbb S^d) \rightarrow \mathbf R$ is given by
\begin{equation}\label{ex:suppfunction3}
\phi^\prime_\theta(h) = \sup_{p\in \Psi_{\mathbb S^d}(\theta)} -h(p) ~.
\end{equation}
As in Example \ref{ex:condineq}, $\phi : \mathcal C(\mathbb S^d) \rightarrow \mathbf R$ is Hadamard differentiable at any $\theta \in \mathcal C(\mathbb S^d)$ at which $\Psi_{\mathbb S^d}(\theta)$ is a singleton, but is only Hadamard directionally differentiable otherwise. \qed

\noindent {\bf Example \ref{ex:stochdom} (cont.)} For any $\theta = (\theta^{(1)},\theta^{(2)})\in \ell^\infty(\mathbf R)\times \ell^\infty(\mathbf R)$ define the sets $B_0(\theta) \equiv \{u \in \mathbf R : \theta^{(1)}(u) = \theta^{(2)}(u)\}$ and $B_+(\theta) \equiv \{u\in \mathbf R: \theta^{(1)}(u) > \theta^{(2)}(u)\}$. It then follows that $\phi$ is Hadamard directionally differentiable at any $\theta \in \ell^\infty(\mathbf R)\times \ell^\infty(\mathbf R)$, and that
\begin{equation}\label{ex:stochdom2}
\phi^\prime_\theta(h) = \int_{B_+(\theta)}(h^{(1)}(u) - h^{(2)}(u))w(u)du + \int_{B_0(\theta)}\max\{h^{(1)}(u)-h^{(2)}(u),0\}w(u)du
\end{equation}
for $h = (h^{(1)},h^{(2)})\in \ell^\infty(\mathbf R)\times \ell^\infty(\mathbf R)$ -- see Lemma \ref{lm:auxstochdiff} in the Appendix. In particular, if $B_0(\theta)$ has zero Lebesgue measure, then $\phi$ is Hadamard differentiable at $\theta$. \qed

\noindent {\bf Example \ref{ex:lrorder} (cont.)} Lemma 3.2 in \cite{BeareandMoon2012} establishes the Hadamard directional differentiability of $\mathcal M :\ell^\infty([0,1])\rightarrow \ell^\infty([0,1])$ tangentially to $\mathcal C([0,1])$ at any concave $\theta \in \ell^\infty([0,1])$. Since norms are directionally differentiable at zero, we have
\begin{equation}\label{ex:lrorder3}
\phi^\prime_\theta(h) = \Big\{\int_0^1 (\mathcal M_\theta^\prime(h)(u) - h(u))^2du\Big\}^{\frac{1}{2}}
\end{equation}
where $\mathcal M_\theta^\prime : \mathcal C([0,1])\rightarrow \ell^\infty([0,1])$ is the Hadamard directional derivative of $\mathcal M$ at $\theta$. \qed

\subsection{The Delta Method}

While the Delta method for Hadamard differentiable functions has become a standard tool in econometrics \citep{Vaart1998}, the availability of an analogous result for Hadamard directional differentiable maps does not appear to be as well known. To the best of our knowledge, this powerful generalization was independently established in \cite{Shapiro1991} and \cite{Dumbgen1993}, but only recently employed in econometrics; see \cite{BeareandMoon2012} and \cite{Kaido2013dual} for examples.

We next aim to establish a mild extension of the results in \cite{Shapiro1991} and \cite{Dumbgen1993} by showing the Delta method holds in probability -- a conclusion we require for our subsequent derivations. Towards this end, we impose the following:

\begin{assumption}\label{ass:structure}
(i) $\mathbb D$ and $\mathbb E$ are Banach spaces with norms $\|\cdot\|_{\mathbb D}$ and $\|\cdot\|_{\mathbb E}$ respectively; (ii) $\phi :\mathbb D_\phi \subseteq \mathbb D \rightarrow \mathbb E$ is Hadamard directionally differentiable at $\theta_0$ tangentially to $\mathbb D_0$; (iii) The map $\phi_{\theta_0}^\prime$ can be continuously extended to $\mathbb D$ (rather than $\mathbb D_0 \subseteq \mathbb D)$
\end{assumption}

\begin{assumption}\label{ass:paramest}
(i) $\theta_0 \in \mathbb D_\phi$ and there are $\hat \theta_n :\{X_i\}_{i=1}^n \rightarrow \mathbb D_\phi$ such that, for some $r_n \uparrow \infty$,  $r_n\{\hat \theta_n - \theta_0\} \stackrel{L}{\rightarrow } \mathbb G_0$ in $\mathbb D$; (ii) $\mathbb G_0$ is tight and its support is included in $\mathbb D_0$.
\end{assumption}

Assumptions \ref{ass:structure}(i)-(ii) formalize our previous discussion by requiring that the map $\phi:\mathbb D_\phi \rightarrow \mathbb E$ be Hadamard directionally differentiable at $\theta_0$. In turn, Assumption \ref{ass:structure}(iii) allows us to view the map $\phi_{\theta_0}^\prime$ as well defined and continuous on all of $\mathbb D$ (rather than just $\mathbb D_0$), and is automatically satisfied when $\mathbb D_0$ is closed; see Remark \ref{rm:extension}. We emphasize, however, that Assumption \ref{ass:structure}(iii) does not demand differentiability of $\phi$ tangentially to $\mathbb D$ -- i.e. the extension of $\phi_{\theta_0}^\prime$ need not satisfy \eqref{def:had2} for $h \in \mathbb D \setminus \mathbb D_0$.\footnote{For instance, in Example \ref{ex:condineq} $\phi$ is differentiable tangentially to $\mathbb D_0 = \mathcal C(\mathcal F)$, but the map $\phi_{\theta}^\prime$ in \eqref{ex:condineq2} is naturally well defined and continuous on $\mathbb D = \ell^\infty(\mathcal F)$.} In Assumption \ref{ass:paramest}(i), we additionally impose the existence of an estimator $\hat \theta_n$ for $\theta_0$ that is asymptotically distributed according to $\mathbb G_0$ in the Hoffman-J{\o}rgensen sense. The scaling $r_n$ equals $\sqrt n$ in Examples \ref{ex:absmean}-\ref{ex:lrorder}, but may differ in nonparametric problems. Finally, Assumption \ref{ass:paramest}(ii) requires that the support of the limiting process $\mathbb G_0$ be included on the tangential set $\mathbb D_0$, and imposes the regularity condition that the random variable $\mathbb G_0$ be tight.

\begin{remark}\label{rm:extension} \rm
If $\mathbb D_0$ is closed, then the continuity of $\phi_{\theta_0}^\prime : \mathbb D_0\rightarrow \mathbb E$ and Theorem 4.1 in \cite{dugundji:1951} imply that $\phi_{\theta_0}^\prime$ admits a continuous extension to $\mathbb D$ -- i.e. there exists a continuous map $\bar \phi_{\theta_0}^\prime : \mathbb D \rightarrow \mathbb E$ such that $\bar \phi_{\theta_0}^\prime(h) = \phi_{\theta_0}^\prime(h)$ for all $h \in \mathbb D_0$. Thus, if $\mathbb D_0$ is closed, then Assumption \ref{ass:structure}(iii) is automatically satisfied.  \qed
\end{remark}

Assumptions \ref{ass:structure} and \ref{ass:paramest} suffice for establishing the validity of the Delta method.

\begin{theorem}\label{th:delta}
If Assumptions \ref{ass:structure}(i)-(ii), \ref{ass:paramest} hold, then $r_n\{\phi(\hat \theta_n) - \phi(\theta_0)\} \stackrel{L}{\rightarrow} \phi^\prime_{\theta_0}(\mathbb G_0)$. If in addition Assumption \ref{ass:structure}(iii) is also satisfied, then it follows that
\begin{equation}\label{th:deltadisp}
r_n\{\phi(\hat \theta_n) - \phi(\theta_0)\} = \phi_{\theta_0}^\prime(r_n\{\hat \theta_n - \theta_0\}) + o_p(1) ~.
\end{equation}
\end{theorem}

The intuition behind Theorem \ref{th:delta} is the same that motivates the traditional Delta method. Heuristically, the theorem can be obtained from the approximation
\begin{equation}\label{th:deltaintui}
r_n\{\phi(\hat \theta_n) - \phi(\theta_0)\} \approx \phi_{\theta_0}^\prime(r_n\{\hat \theta_n - \theta_0\}) ~,
\end{equation}
Assumption \ref{ass:paramest}(i), and the continuous mapping theorem applied to $\phi_{\theta_0}^\prime$. Thus, the key requirement is not that $\phi_{\theta_0}^\prime$ be linear, or equivalently that $\phi$ be Hadamard differentiable, but rather that \eqref{th:deltaintui} holds in an appropriate sense -- a condition ensured by Hadamard directional differentiability. Following this insight, Theorem \ref{th:delta} can be established using the same arguments as in the proof of the traditional Delta method \citep{Vaart1996}. It is worth noting that directional differentiability of $\phi$ is only assumed at $\theta_0$. In particular, continuity of $\phi_{\theta_0}^\prime$ in $\theta_0$ is not required since such condition is often violated; see, e.g., Example \ref{ex:absmean}. Strengthening the Delta method to hold in probability further requires Assumption \ref{ass:structure}(iii) to ensure $\phi_{\theta_0}^\prime(r_n\{\hat \theta_n - \theta_0\})$ is well defined.\footnote{Without Assumption \ref{ass:structure}(iii), the domain of $\phi_{\theta_0}^\prime$ must include $\mathbb D_0$, but possibly not $\mathbb D\setminus \mathbb D_0$. Thus, since $r_n\{\hat \theta_n - \theta_0\}$ may not belong to $\mathbb D_0$, $\phi_{\theta_0}^\prime(r_n\{\hat \theta_n - \theta_0\})$ may otherwise not be well defined.}

\begin{remark}\label{rm:exdist} \rm
Theorem \ref{th:delta} immediately delivers the relevant asymptotic distributions in Examples \ref{ex:absmean} and \ref{ex:interbounds}, as well as in Examples \ref{ex:condineq} and \ref{ex:stochdom} provided a functional central limit theorem applies. In turn, the asymptotic distribution in Examples \ref{ex:suppfunction} and \ref{ex:lrorder} can be obtained by employing Theorem \ref{th:delta} together with distributional convergence results for $\sqrt n\{\hat \theta_n - \theta_0\}$ as available, for example, in \cite{Kaido_Santos2013} for support functions and \cite{BeareandMoon2012} for ordinal dominance curves. \qed
\end{remark}

\section{The Bootstrap}

While Theorem \ref{th:delta} enables us to obtain an asymptotic distribution, a suitable method for estimating this limiting law is still required. In this section we assume that the bootstrap ``works" for $\hat \theta_n$ and examine how to leverage this result to estimate the asymptotic distribution of $r_n\{\phi(\hat \theta_n) - \phi(\theta_0)\}$. We show that, whenever $\hat \theta_n$ is asymptotically Gaussian, bootstrap consistency is necessarily lost unless $\phi$ is (fully) Hadamard differentiable. As an alternative to the bootstrap, we thus propose a procedure which generalizes existing approaches and remains valid when $\phi$ is not fully differentiable.

\subsection{Bootstrap Setup}

We begin by introducing the general setup under which we examine bootstrap consistency. Throughout, we let $\hat \theta_n^*$ denote a ``bootstrapped version" of $\hat \theta_n$, and assume the limiting distribution of $r_n\{\hat \theta_n - \theta_0\}$ can be consistently estimated by the law of
\begin{equation}\label{eq:bootset1}
r_n\{\hat \theta_n^* - \hat \theta_n\}
\end{equation}
conditional on the data. In order to formally define $\hat \theta_n^*$, while allowing for diverse resampling schemes, we simply impose that $\hat \theta_n^*$ be a function mapping the data $\{X_i\}_{i=1}^n$ and random weights $\{W_i\}_{i=1}^n$ that are independent of $\{X_i\}_{i=1}^n$ into $\mathbb D_\phi$. This abstract definition suffices for encompassing the nonparametric, Bayesian, block, $m$ out of $n$, score, and weighted bootstrap as special cases. 

Formalizing the notion of bootstrap consistency further requires us to employ a measure of distance between the limiting distribution and its bootstrap estimator. Towards this end, we follow \cite{Vaart1996} and utilize the bounded Lipschitz metric. Specifically, for a metric space $\mathbf A$ with norm $\|\cdot\|_{\mathbf A}$, we denote the set of Lipschitz functionals whose level and Lipschitz constant are bounded by one by
\begin{equation}\label{eq:bootset2}
\text{BL}_1(\mathbf A) \equiv \{f : \mathbf A \rightarrow \mathbf R :  |f(a)| \leq 1 \text{ and } |f(a) - f(a^\prime)|\leq \|a - a^\prime\|_{\mathbf A} \text{ for all } a,a^\prime \in \mathbf A\} ~.
\end{equation}
The bounded Lipschitz distance between two measures $L_1$ and $L_2$ on $\mathbf A$ then equals the largest discrepancy in the expectation they assign to functions in $\text{BL}_1(\mathbf A)$, denoted
\begin{equation}\label{eq:bootset3}
d_{\text{BL}}(L_1,L_2) \equiv \sup_{f \in \text{BL}_1(\mathbf A)} |\int f(a) dL_1(a) - \int f(a) dL_2(a)| ~.
\end{equation}

Given the introduced notation, we can measure the distance between the law of $r_n\{\hat \theta_n^* - \hat \theta_n\}$ conditional on $\{X_i\}_{i=1}^n$, and the limiting distribution of $r_n\{\hat \theta_n - \theta_0\}$ by\footnote{More precisely, $E[f(r_n\{\hat \theta_n^* - \hat \theta_n\})|\{X_i\}_{i=1}^n]$ denotes the outer expectation with respect to the joint law of $\{W_i\}_{i=1}^n$, treating the observed data $\{X_i\}_{i=1}^n$ as constant.}
\begin{equation}\label{eq:bootset4}
\sup_{f\in \text{BL}_1(\mathbb D)} |E[f(r_n\{\hat \theta_n^*-\hat \theta_n\})|\{X_i\}_{i=1}^n] - E[f(\mathbb G_0)]| ~.
\end{equation}
Employing the distribution of $r_n\{\hat \theta_n^* - \hat \theta_n\}$ conditional on the data to approximate the distribution of $\mathbb G_0$ is then asymptotically justified if their distance, equivalently \eqref{eq:bootset4}, converges in probability to zero. This type of consistency can in turn be exploited to validate the use of critical values obtained from the distribution of $r_n\{\hat \theta_n^* - \hat\theta_n\}$ conditional on $\{X_i\}_{i=1}^n$ to conduct inference or construct confidence regions; see Remark \ref{rm:weakconv}.

We formalize the above discussion by imposing the following assumptions on $\hat \theta_n^*$:

\begin{assumption}\label{ass:primboot}
(i) $\hat \theta_n^* : \{X_i,W_i\}_{i=1}^n \rightarrow \mathbb D_\phi$ with $\{W_i\}_{i=1}^n$ independent of $\{X_i\}_{i=1}^n$; (ii) $\hat \theta_n^*$ satisfies $\sup_{f \in \text{BL}_1(\mathbb D)} |E[f(r_n\{\hat \theta_n^* - \hat \theta_n\})|\{X_i\}_{i=1}^n] - E[f(\mathbb G_0)]| = o_p(1)$.
\end{assumption}

\begin{assumption}\label{ass:extraboot}
(i) The sequence $r_n\{\hat \theta_n^* - \hat \theta_n\}$ is asymptotically measurable (jointly in $\{X_i,W_i\}_{i=1}^n$); (ii) $f(r_n\{\hat \theta_n^* - \hat \theta_n\})$ is a measurable function of $\{W_i\}_{i=1}^n$ outer almost surely in $\{X_i\}_{i=1}^n$ for any continuous and bounded $f:\mathbb D \rightarrow \mathbf R$ .
\end{assumption}

Assumption \ref{ass:primboot}(i) defines $\hat \theta_n^*$ in accord with our discussion, while Assumption \ref{ass:primboot}(ii) imposes the consistency of the law of $r_n\{\hat \theta_n^* - \hat \theta_n\}$ conditional on the data for the distribution of $\mathbb G_0$ -- i.e. the bootstrap ``works" for the estimator $\hat \theta_n$. In addition, in Assumption \ref{ass:extraboot} we further demand mild measurability requirements on $r_n\{\hat \theta_n^* - \hat \theta_n\}$. These requirements are automatically satisfied whenever $\hat \theta_n$ and $\hat \theta_n^*$ correspond to the empirical and bootstrapped empirical processes indexed by a suitably Donsker class of functions.

\begin{remark}\label{rm:weakconv} \rm
In the special case where $\mathbb D = \mathbf R^d$, Assumption \ref{ass:primboot}(ii) implies that:
\begin{equation}\label{rm:weakconv1}
\sup_{t\in A} |P(r_n\{\hat \theta_n^* - \hat \theta_n\} \leq t | \{X_i\}_{i=1}^n) - P(\mathbb G_0 \leq t)| = o_p(1)
\end{equation}
for any closed subset $A$ of the continuity points of the cdf of $\mathbb G_0$; see \cite{Kosorok2008}. Thus, consistency in the bounded Lipschitz metric implies consistency of the corresponding cdfs. Result \eqref{rm:weakconv1} then readily yields consistency of the corresponding quantiles at points at which the cdf of $\mathbb G_0$ is continuous and strictly increasing. \qed
\end{remark}

\subsection{Bootstrap Failure}\label{sec:bootfail}

When the transformation $\phi : \mathbb D_\phi \rightarrow \mathbb E$ is Hadamard differentiable at $\theta_0$, the consistency of the bootstrap is inherited by the transformation itself. In other words, if Assumption \ref{ass:primboot}(ii) is satisfied, and $\phi$ is (fully) Hadamard differentiable, then the asymptotic distribution of $r_n\{\phi(\hat \theta_n) - \phi(\theta_0)\}$ can be consistently estimated by the law of
\begin{equation}\label{eq:necsuf1}
r_n\{\phi(\hat \theta_n^*) - \phi(\hat \theta_n)\}
\end{equation}
conditional on the data; see \cite{bickel1981some} and \cite{Vaart1996} respectively for the finite and infinite dimensional cases. For conciseness, we refer to the law of \eqref{eq:necsuf1} conditional on the data as the ``standard" bootstrap.

Unfortunately, while the Delta method generalizes to Hadamard directionally differentiable functionals, we know by way of example that the consistency of the standard bootstrap may not \citep{Bickel_Gotze_Zwet1997, Andrews2000Bootstrap}. These examples serve as a warning that the standard bootstrap may fail when $\phi$ is not (fully) Hadamard differentiable, yet can provide little guidance as to whether the bootstrap is actually valid in particular applications. Our first main result establishes that these examples are in fact special cases of a deeper principle, namely that whenever $\mathbb G_0$ is Gaussian the standard bootstrap is consistent \emph{if and only if} $\phi$ is (fully) differentiable at $\theta_0$.

\begin{theorem}\label{th:gauslin}
Let Assumptions \ref{ass:structure}, \ref{ass:paramest}, \ref{ass:primboot}, and \ref{ass:extraboot} hold, and suppose $\mathbb G_0$ is a Gaussian measure whose support is a vector subspace of $\mathbb D$. Then, it follows that $\phi$ is (fully) Hadamard differentiable at $\theta_0 \in \mathbb D_\phi$ tangentially to the support of $\mathbb G_0$ if and only if
\begin{equation}\label{th:gauslindisp}
\sup_{f \in \text{BL}_1(\mathbb E)} |E[f(r_n\{\phi(\hat \theta_n^*) - \phi(\hat \theta_n)\})|\{X_i\}_{i=1}^n] - E[f(\phi_{\theta_0}^\prime(\mathbb G_0))]| = o_p(1) ~.
\end{equation}
\end{theorem}

A powerful implication of Theorem \ref{th:gauslin} is that in verifying whether the standard bootstrap is valid at a conjectured $\theta_0$, a researcher need only verify whether $\phi$ is (fully) differentiable at $\theta_0$. In effect, Theorem \ref{th:gauslin} thus reduces the potentially challenging statistical problem of verifying bootstrap validity to a simple and purely analytical calculation; see Remark \ref{rm:dumbgen}. The theorem requires both that $\mathbb G_0$ be Gaussian and that its support be a vector subspace of $\mathbb D$. The former requirement may be relaxed at the cost of additional notation, and we thus focus on the Gaussian case due to its ubiquity; see Remark \ref{rm:hpconnect}. In turn, we note that under Gaussianity the condition that the support of $\mathbb G_0$ be a vector subspace is equivalent to zero (in $\mathbb D$) belonging to the support of $\mathbb G_0$.

A further implication of Theorem \ref{th:gauslin} that merits discussion follows from exploiting that Gaussianity of $\mathbb G_0$ and bootstrap consistency together imply $\phi$ is (fully) differentiable and hence that $\phi_{\theta_0}^\prime$ must be linear; recall Proposition \ref{pro:equivalence}. In particular, whenever $\phi_{\theta_0}^\prime$ is linear and $\mathbb G_0$ is Gaussian $\phi_{\theta_0}^\prime(\mathbb G_0)$ must also be Gaussian (in $\mathbb E$), and thus bootstrap consistency implies Gaussianity of $\phi_{\theta_0}^\prime(\mathbb G_0)$ or, equivalently, Gaussianity of the asymptotic distribution of $\phi(\hat \theta_n)$. Conversely, we conclude that \emph{the standard bootstrap fails whenever the limiting distribution is not Gaussian}. Thus, our results imply that the presence of a non-Gaussian limiting distribution may be viewed by practitioners as a simple yet reliable signal of the failure of the standard bootstrap. We formalize this conclusion in the following Corollary to Theorem \ref{th:gauslin}.

\begin{corollary}\label{cor:nonnorm}
Let Assumptions \ref{ass:structure}, \ref{ass:paramest}, \ref{ass:primboot}, \ref{ass:extraboot} hold, and $\mathbb G_0$ be a Gaussian measure whose support is a vector subspace of $\mathbb D$. If the limiting distribution of $r_n\{\phi(\hat \theta_n) - \phi(\theta_0)\}$ is not Gaussian, then it follows that the standard bootstrap is inconsistent.
\end{corollary}

\begin{remark}\label{rm:dumbgen} \rm
The limit of the standard bootstrap measure was first studied in \cite{Dumbgen1993}, whose results imply bootstrap consistency is equivalent to the distribution of
\begin{equation}\label{rm:dumbgen1}
\phi_{\theta_0}^\prime(\mathbb G_0 + h) - \phi_{\theta_0}^\prime(h)
\end{equation}
being constant in $h$ for all $h$ in the support of $\mathbb G_0$.\footnote{See Theorem \ref{th:bootiff} in Appendix A for an analogous result under Assumptions \ref{ass:structure}, \ref{ass:paramest}, \ref{ass:primboot}, and \ref{ass:extraboot}.} This characterization of bootstrap consistency does not rely on Gaussianity and thus applies under more general conditions than those of Theorem \ref{th:gauslin}. Our results complement \cite{Dumbgen1993} by showing that, under the additional requirement that $\mathbb G_0$ be Gaussian, bootstrap consistency is in fact also equivalent to $\phi$ being (fully) differentiable at $\theta_0$. Thus, we obtain a purely analytical and simpler to verify condition than \eqref{rm:dumbgen1} that enables practitioners to easily assess bootstrap validity in the ubiquitous case where $\mathbb G_0$ is Gaussian. \qed
\end{remark}

\begin{remark} \label{rm:hpconnect} \rm
Gaussianity of $\mathbb G_0$ plays an important role in the proof of Theorem \ref{th:gauslin} in enabling us to relate the distribution of \eqref{rm:dumbgen1} to that of $\phi_{\theta_0}^\prime(\mathbb G_0)$ through the Cameron-Martin theorem. A similar insight was employed by \cite{Vaart1991differentibility} and \cite{hirano:porter:2009} who compare characteristic functions in a limit experiment to conclude regular estimability of a functional implies its differentiability. More generally, Theorem \ref{th:gauslin} can be shown to hold provided the support of $\mathbb G_0$ is a vector subspace of $\mathbb D$ and that the Radon-Nikodym derivative of the distribution of $\mathbb G_0 + h$ with respect to that of $\mathbb G_0$ is suitably smooth for all $h$ in a dense subspace of the support of $\mathbb G_0$. \qed
\end{remark}

\subsection{An Alternative Approach}

While Theorem \ref{th:gauslin} establishes that standard bootstrap procedures are inconsistent, alternative resampling schemes such as the $m$ out of $n$ bootstrap \citep{shao1994bootstrap} and subsampling \citep{PWR1999} can nonetheless remain valid.\footnote{Alternatively, though potentially quite conservative, projection methods also remain valid under virtually no assumptions on $\phi$; see, e.g., \cite{dufour2005projection}, \cite{romano:shaikh:2008, Romano_Shaikh2010}, and \cite{ham:woutersen:2013}.} However, in the moment inequalities literature -- an important special case of our setting -- these approaches have been superseded by alternatives based on generalized moment selection \citep{andrews:soares:2010}.  In what follows, we show that generalized moment selection can be interpreted as a procedure that implicitly composes an estimate $\hat \phi_n^\prime$ of the directional derivative $\phi_{\theta_0}^\prime$ with a bootstrap approximation to the asymptotic distribution of $\hat \theta_n$. Based on this reinterpretation of generalized moment selection, we in turn develop a resampling scheme that is applicable to our general setting.

\subsubsection{Consistent Alternative}

Heuristically, the inconsistency of the standard bootstrap arises from its inability to properly estimate the directional derivative $\phi_{\theta_0}^\prime$. However, the underlying bootstrap process $r_n\{\hat \theta_n^* - \hat \theta_n\}$ still provides a consistent estimator for the law of $\mathbb G_0$. Intuitively, it should therefore be possible to obtain a consistent estimator for the limiting distribution found in Theorem \ref{th:delta} by employing the law, conditional on the data, of
\begin{equation}\label{eq:alt1}
\hat \phi_n^\prime(r_n\{\hat \theta_n^* - \hat \theta_n\})
\end{equation}
for $\hat \phi_n^\prime :\mathbb D\rightarrow \mathbb E$ a suitable estimator of the directional derivative $\phi_{\theta_0}^\prime : \mathbb D_0 \rightarrow \mathbb E$. This simple intuition is actually implicit in generalized moment selection and other influential inferential methods designed for specific examples of $\phi : \mathbb D \rightarrow \mathbb E$. In fact, in Section \ref{sec:altbootex} below we show how such existing procedures may be recast as a special case of \eqref{eq:alt1}.

In order for the described intuition to be valid, we require $\hat \phi_n^\prime$ to satisfy the following:

\begin{assumption}\label{ass:derest}
$\hat \phi_n^\prime : \mathbb D \rightarrow \mathbb E$ is a function of $\{X_i\}_{i=1}^n$, satisfying for every compact set $K\subseteq \mathbb D_0$, $K^{\delta} \equiv \{a \in \mathbb D : \inf_{b \in K} \|a - b\|_{\mathbb D} < \delta\}$, and every $\epsilon > 0$, the property:
\begin{equation}\label{ass:derestdisp1}
\lim_{\delta \downarrow 0} \limsup_{n\rightarrow \infty} P\Big( \sup_{h \in K^{\delta}} \| \hat \phi_n^\prime(h) - \phi^\prime_{\theta_0}(h)\|_{\mathbb E} > \epsilon\Big) = 0 ~.
\end{equation}
\end{assumption}

Unfortunately, the requirement in \eqref{ass:derestdisp1} is complicated by the presence of the $\delta$-enlargement of $K$. Without such enlargement, requirement \eqref{ass:derestdisp1} could just be interpreted as demanding that $\hat \phi_n^\prime$ be uniformly consistent for $\phi_{\theta_0}^\prime$ on compact sets $K \subseteq \mathbb D_0$. Heuristically, the need to consider $K^\delta$ arises from $r_n\{\hat \theta_n^* - \hat \theta_n\}$ only being guaranteed to lie in $\mathbb D$ and not necessarily $\mathbb D_0$. However, because $\mathbb G_0$ lies in compact subsets of $\mathbb D_0$ with arbitrarily high probability, it is possible to conclude that $r_n\{\hat \theta_n^* - \hat \theta_n\}$ will eventually be ``close" to such subsets of $\mathbb D_0$. Thus, $\hat \phi_n^\prime$ need only be well behaved in arbitrary small neighborhoods of compact sets in $\mathbb D_0$, which is the requirement imposed in Assumption \ref{ass:derest}. It is worth noting, however, that in many applications stronger, but simpler, conditions than \eqref{ass:derestdisp1} can be easily verified. For instance, under appropriate additional requirements, the $\delta$ factor in \eqref{ass:derestdisp1} may be ignored, and it may even suffice to just verify $\hat \phi_n^\prime(h)$ is consistent for $\phi_{\theta_0}^\prime(h)$ for every $h \in \mathbb D_0$; see Remarks \ref{rm:D0orD} and \ref{rm:neighok}. Moreover, in certain applications the specific choice of estimator $\hat \phi_n$ may possess additional structure that enables us to replace Assumption \ref{ass:derest} with ``lower level" conditions -- see, e.g., \cite{hong:li:2014} who propose estimating $\hat \phi_n^\prime$ through (directional) numerical differentiation.

\begin{remark}\rm \label{rm:D0orD}
In certain applications, it is sufficient to require $\hat \phi_n^\prime:\mathbb D\rightarrow \mathbb E$ to satisfy
\begin{equation}\label{rm:D0orD1}
\sup_{h\in K}\|\hat \phi_n^\prime(h) - \phi_{\theta_0}^\prime(h)\|_{\mathbb E} = o_p(1) ~,
\end{equation}
for any compact set $K \subseteq \mathbb D$. For instance, if $\mathbb D = \mathbf R^d$, then the closure of $K^\delta$ is compact in $\mathbb D$ for any compact $K\subseteq \mathbb D_0$, and hence \eqref{rm:D0orD1} implies \eqref{ass:derestdisp1}. Alternatively, if $\mathbb D$ is separable, $r_n\{\hat \theta_n^* - \hat \theta_n\}$ is Borel measurable as a function of $\{X_i,W_i\}_{i=1}^n$ and tight for each $n$, then $r_n\{\hat \theta_n^* - \hat \theta_n\}$ is uniformly tight and \eqref{rm:D0orD1} may be used in place of \eqref{ass:derestdisp1}.\footnote{Under uniform tightness, for every $\epsilon > 0$ there is a compact set $K$ such that $\limsup_{n\rightarrow \infty} P(r_n\{\hat \theta_n^* - \hat \theta_n\} \notin K) < \epsilon$. In general, however, we only know $r_n\{\hat \theta_n^* - \hat \theta_n\}$ to be asymptotically tight, in which case we are only guaranteed $\limsup_{n\rightarrow \infty} P(r_n\{\hat \theta_n^* - \hat \theta_n\} \notin K^{\delta}) < \epsilon$ for every $\delta >0$.} \qed
\end{remark}

\begin{remark}\rm \label{rm:neighok}
Assumption \ref{ass:derest} greatly simplifies whenever the modulus of continuity of $\hat \phi_n^\prime : \mathbb D \rightarrow \mathbb E$ can be controlled outer almost surely. For instance, if $\|\hat \phi_n^\prime(h_1) - \hat \phi_n^\prime(h_2)\|_{\mathbb E} \leq C\|h_1-h_2\|_{\mathbb D}$ for some $C < \infty$ and all $h_1,h_2 \in \mathbb D$, then showing that for any $h \in \mathbb D_0$
\begin{equation}\label{rm:neighok1}
\|\hat \phi_n^\prime(h) - \phi_{\theta_0}^\prime(h)\|_{\mathbb E} = o_p(1)
\end{equation}
suffices for establishing \eqref{ass:derestdisp1} holds; see Lemma \ref{lm:simpsuff} in the Appendix. This observation is particularly helpful in the analysis of Examples \ref{ex:condineq} and \ref{ex:suppfunction}; see Section \ref{sec:altbootex}. \qed
\end{remark}

Given Assumption \ref{ass:derest} we can establish the validity of the proposed procedure.

\begin{theorem}\label{th:modboot}
Under Assumptions \ref{ass:structure}, \ref{ass:paramest}, \ref{ass:primboot}, \ref{ass:extraboot} and \ref{ass:derest}, it follows that
\begin{equation}\label{th:modbootdisp1}
\sup_{f \in \text{BL}_1(\mathbb E)} |E[f(\hat \phi_n^\prime(r_n\{\hat \theta_n^* - \hat \theta_n\}))|\{X_i\}_{i=1}^n] - E[f(\phi_{\theta_0}^\prime(\mathbb G_0))]| = o_p(1) ~.
\end{equation}
\end{theorem}

Theorem \ref{th:modboot} shows that the law of $\hat \phi_n^\prime(r_n\{\hat \theta_n^* - \hat \theta_n\})$ conditional on the data is indeed consistent for the limiting distribution of $r_n\{\phi(\hat \theta_n) - \phi(\theta_0)\}$ derived in Theorem \ref{th:delta}. In particular, when $\phi(\hat \theta_n)$ is a test statistic, and hence scalar valued, Theorem \ref{th:modboot} enables us to compute critical values for inference by employing the quantiles of the simulated finite sample distribution of $\hat \phi_n^\prime(r_n\{\hat \theta_n^* - \hat \theta_n\})$ conditional on the data $\{X_i\}_{i=1}^n$ (but not $\{W_i\}_{i=1}^n$). The following immediate corollary formally establishes this claim.

\begin{corollary}\label{cor:critval}
Let Assumptions \ref{ass:structure}, \ref{ass:paramest}, \ref{ass:primboot}, \ref{ass:extraboot} and \ref{ass:derest} hold,  $\mathbb E = \mathbf R$, and define
\begin{equation}\label{cor:critvaldisp}
\hat c_{1-\alpha} \equiv \inf\{ c : P(\hat \phi_n^\prime(r_n\{\hat \theta_n^* - \hat \theta_n\}) \leq c |\{X_i\}_{i=1}^n) \geq 1-\alpha\} ~.
\end{equation}
If the cdf of $\phi_{\theta_0}^\prime(\mathbb G_0)$ is continuous and strictly increasing at its $1-\alpha$ quantile denoted $c_{1-\alpha}$, then it follows that $\hat c_{1-\alpha} \stackrel{p}{\rightarrow} c_{1-\alpha}$.
\end{corollary}

It is worth noting that $\phi_{\theta_0}^\prime$ being the directional derivative of $\phi$ at $\theta_0$ is actually never exploited in the proofs of Theorem \ref{th:modboot} or Corollary \ref{cor:critval}. Therefore, these results can more generally be interpreted as providing a method for approximating distributions of random variables that are of the form $\tau(\mathbb G_0)$, where $\mathbb G_0 \in \mathbb D$ is a tight random variable and $\tau:\mathbb D\rightarrow \mathbb E$ is an unknown continuous map; see, e.g., \cite{chen:fang:2015} for an application of this principle to estimate first order degenerate asymptotic distributions. Finally, it is important to emphasize that due to an appropriate lack of continuity of $\phi_{\theta_0}^\prime$ in $\theta_0$, the ``naive" estimator $\hat \phi_n^\prime = \phi^\prime_{\hat \theta_n}$ often fails to satisfy Assumption \ref{ass:derest}. Nonetheless, alternative estimators are still easily obtained as we next discuss.

\subsubsection{Examples Revisited} \label{sec:altbootex}

We revisit Examples \ref{ex:absmean}-\ref{ex:lrorder} to illustrate how our approach is in fact implicit in popular inferential methods designed for special cases of our setting -- in this sense, Theorem \ref{th:modboot} may therefore be viewed as a generalization of such approaches.  For conciseness, we group the analysis of examples that share a similar structure.

\noindent {\bf Examples \ref{ex:absmean} and \ref{ex:interbounds} (cont.)} In the context of Example \ref{ex:interbounds}, let $\{X_i\}_{i=1}^n$ be an i.i.d. sample with $X_i = (X_i^{(1)},X_i^{(2)})^\prime \in \mathbf R^2$, and define $\bar X^{(j)} \equiv \frac{1}{n}\sum_i X_i^{(j)}$ for $j\in \{1,2\}$. Denoting $\hat j^* = \arg\max_{j\in\{1,2\}} \bar X^{(j)}$ and letting $\kappa_n \downarrow 0$ satisfy $\kappa_n\sqrt n\uparrow \infty$, we then define
\begin{equation}\label{ex:interbounds3}
\hat \phi_n^\prime (h) = \begin{cases}
h^{(\hat j^*)}   & \text{ if } |\bar X^{(1)} - \bar X^{(2)}| > \kappa_n  \\
\max\{h^{(1)},h^{(2)}\}  &\text{ if } |\bar X^{(1)} - \bar X^{(2)}| \leq \kappa_n
\end{cases} ~,
\end{equation}
(compare to \eqref{ex:interbounds2}). Under appropriate moment restrictions, it is then straightforward to verify Assumption \ref{ass:derest} holds, since $\hat \phi_n^\prime: \mathbf R^2 \rightarrow \mathbf R$ in fact satisfies
\begin{equation}\label{ex:interbounds4}
\liminf_{n\rightarrow \infty} P\Big( \hat \phi_n^\prime(h) = \phi_{\theta_0}^\prime(h) \text{ for all } h \in \mathbf R^2\Big) = 1 ~.
\end{equation}
If $\{X_i^*\}_{i=1}^n$ is a sample drawn with replacement from $\{X_i\}_{i=1}^n$, and $\bar X^* = \frac{1}{n} \sum_i X_i^*$, then \eqref{eq:alt1} reduces to $\hat \phi_n^\prime(\sqrt n\{\bar X^* - \bar X\})$, which was originally studied in \cite{andrews:soares:2010} and \cite{bugni:2010} for conducting inference in moment inequalities models. Example \ref{ex:absmean} can be studied in a similar manner and we therefore omit its analysis. \qed

\noindent {\bf Examples \ref{ex:condineq} and \ref{ex:suppfunction} (cont.)} In Example \ref{ex:condineq}, recall $\Psi_{\mathcal F}(\theta) \equiv \arg\max_{f\in \mathcal F} \theta(f)$ and suppose $\hat \Psi_{\mathcal F}(\theta_0)$ is a Hausdorff consistent estimate of $\Psi_{\mathcal F}(\theta_0)$ -- i.e. it satisfies\footnote{For subsets $A,B$ of a metric space with norm $\|\cdot\|$, the directed Hausdorff distance is $\vec d_H(A,B) \equiv \sup_{a\in A}\inf_{b\in B}\|a-b\|$, and the Hausdorff distance is $d_H(A,B,\|\cdot\|) \equiv \max\{\vec d_H(A,B,\|\cdot\|), \vec d_H(B,A,\|\cdot\|)\}$.}
\begin{equation}\label{eq:conineq3}
d_H(\Psi_{\mathcal F}(\theta_0), \hat \Psi_{\mathcal F}(\theta_0),\|\cdot\|_{L^2(Z)}) = o_p(1) ~.
\end{equation}
A natural estimator for $\phi^\prime_{\theta_0}$ is then given by $\hat \phi_n^\prime : \ell^\infty(\mathcal F) \rightarrow \mathbf R$ equal to (compare to \eqref{ex:condineq2})
\begin{equation}\label{eq:conineq4}
\hat \phi_n^\prime(h) = \sup_{f \in \hat \Psi_{\mathcal F}(\theta_0)} h(f) ~,
\end{equation}
which can easily be shown to satisfy Assumption \ref{ass:derest}; see Lemma \ref{lm:auxsupest} in the Appendix. If the data is i.i.d., $\{(Y_i^*,Z_i^*)\}_{i=1}^n$ is a sample drawn with replacement from $\{(Y_i,Z_i)\}_{i=1}^n$, and $\sqrt n\{\hat \theta_n^*-\hat \theta_n\}$ is the bootstrapped empirical process, then \eqref{eq:alt1} becomes
\begin{equation}\label{eq:conineq5}
\hat \phi_n^\prime(\sqrt n\{\hat \theta_n^* - \hat \theta_n\}) = \sup_{f\in \hat\Psi_{\mathcal F}(\theta_0)} \frac{1}{\sqrt n} \sum_{i=1}^n \{Y_i^*f(Z_i^*) - \frac{1}{n}\sum_{i=1}^n Y_if(Z_i)\} ~,
\end{equation}
which was originally proposed in \cite{andrews:shi:2013} for conducting inference in conditional moment inequalities models. A similar approach is pursued in \cite{Kaido2013dual} and \cite{Kaido_Santos2013} in the context of Example \ref{ex:suppfunction}. \qed

\noindent {\bf Examples \ref{ex:stochdom} and \ref{ex:lrorder} (cont.)} Recall that in Example \ref{ex:stochdom}, $\theta_0 = (\theta_0^{(1)},\theta_0^{(2)})$ with $\theta_0^{(j)} \in \ell^\infty(\mathbf R)$ for $j\in\{1,2\}$, and that $B_0(\theta_0) = \{u\in \mathbf R : \theta^{(1)}_0(u) = \theta^{(2)}_0(u)\}$ and $B_+(\theta_0) = \{u\in \mathbf R: \theta_0^{(1)}(u) > \theta_0^{(2)}(u)\}$. For $\hat B_0(\theta_0)$ and $\hat B_+(\theta_0)$ estimators of $B_0(\theta_0)$ and $B_+(\theta)$ respectively, it is then natural for any $h\in (h^{(1)},h^{(2)})\in \ell^\infty(\mathbf R)\times \ell^\infty(\mathbf R)$ to define
\begin{equation}\label{ex:stochdom3}
\hat \phi_n^\prime(h) = \int_{\hat B_+(\theta_0)} (h^{(1)}(u) - h^{(2)}(u))w(u)du + \int_{\hat B_0(\theta_0)} \max\{h^{(1)}(u) - h^{(2)}(u),0\}w(u)du
\end{equation}
(compare to \eqref{ex:stochdom2}). For $A\triangle B$ the symmetric set difference between sets $A$ and $B$, it is then straightforward to verify Assumption \ref{ass:derest} is satisfied provided the Lebesgue measure of $B_0(\theta_0)\triangle\hat B_0(\theta_0)$ and $B_+(\theta_0)\triangle\hat B_+(\theta_0)$ converges in probability to zero. When $\sqrt n\{\hat \theta_n^* - \hat \theta\}$ is given by the bootstrap empirical process, $\hat \phi_n^\prime(\sqrt n\{\hat \theta_n^* - \hat \theta_n\})$ reduces to the procedure studied in \cite{Linton2010} for testing stochastic dominance. For a related analysis of Example \ref{ex:lrorder} we refer the reader to \cite{BrendanPrelim}. \qed

\subsection{Local Analysis}

As evidenced in Examples \ref{ex:absmean}-\ref{ex:lrorder}, $\phi(\hat \theta_n)$ is not a regular estimator for $\phi(\theta_0)$ whenever $\phi$ is not (fully) Hadamard differentiable at $\theta_0$. In order to evaluate the usefulness of Theorems \ref{th:delta} and \ref{th:modboot} for conducting inference, it is therefore crucial to complement these results with a study of the asymptotic distribution of $r_n\{\phi(\hat \theta_n)-\phi(\theta_0)\}$ under local perturbations to the underlying distribution of the data. In this section, we first develop such a local analysis and then proceed to examine its implications for inference. As in our previous results, we will show that once again challenging statistical questions can be solved by simple analytical calculations concerning the directional derivative.

\subsubsection{Local Limit}

We begin by introducing the framework that will enable us to conduct the desired local analysis. To this end, we let $\mathbf P^{\infty}$ denote the set of possible distributions for $\{X_i\}_{i=1}^\infty$ and for any $Q^{\infty} \in \mathbf P^{\infty}$ we set $Q^{n}$ to equal the distribution of $\{X_i\}_{i=1}^n$ induced by $Q^{\infty}$. In addition, we make the dependence of the parameter $\theta_0$ on the unknown distribution $P^{\infty}$ of $\{X_i\}_{i=1}^\infty$ explicit by letting $\theta_0$ be the value a known map $\theta : \mathbf P^{\infty} \rightarrow \mathbb D_\phi$ takes at the unknown point $P^{\infty}$ -- i.e. $\theta_0 \equiv \theta(P^{\infty})$.\footnote{For instance, in Examples \ref{ex:absmean} and \ref{ex:interbounds} with $\{X_i\}_{i=1}^\infty$ stationary and distributed according to joint and marginal distributions $P^\infty$ and $P$, the known map $\theta(P^\infty)$ is given by $\theta(P^\infty) \equiv \int xdP(x)$.} Given the introduced notation, the following two assumptions impose the structure we require for our local analysis:

\begin{assumption}\label{ass:localset}
There exist a vector space $\Lambda$ and maps $P_{n,\cdot}^{\infty} : \Lambda \rightarrow \mathbf P^{\infty}$, such that: (i) $P_{n,0}^{\infty} = P^\infty$ for all $n$, where $P^\infty$ denotes the distribution of $\{X_i\}_{i=1}^\infty$; (ii) The random variables $\frac{dP_{n,\lambda}^{n}}{dP_{n,0}^{n}}(\{X_i\}_{i=1}^n)$ are uniformly tight under $P_{n,0}^{n}$ for any $\lambda \in \Lambda$.\footnote{Here, $dP_{n,\lambda}^n/dP^n_{n,0} \equiv p_{n,\lambda}/p_{n,0}$ where for any measure $\mu_n$ dominating $P_{n,0}^n$ and $P_{n,\lambda}^n$ we set $p_{n,0} \equiv dP_{n,0}^n/d\mu_n$, $p_{n,\lambda} \equiv dP_{n,\lambda}^n/d\mu_n$, and define $0/0 = 0$ and $c/0 = \infty$ whenever $c > 0$.}
\end{assumption}

\begin{assumption}\label{ass:paramreg}
There is $\theta :\mathbf P^{\infty} \rightarrow \mathbb D_\phi$, with: (i) $\|r_n\{\theta(P_{n,\lambda}^{\infty}) - \theta(P^\infty_{n,0})\} - \theta^\prime(\lambda)\|_{\mathbb D} = o(1)$ for some linear $\theta^\prime : \Lambda \rightarrow \mathbb D_0$; (ii) $r_n\{\hat \theta_n - \theta(P_{n,\lambda}^{\infty})\} \stackrel{L}{\rightarrow} \mathbb G_0$ under $P_{n,\lambda}^{n}$ for any $\lambda \in \Lambda$.
\end{assumption}

Assumption \ref{ass:localset}(i) introduces the distributions $P_{n,\lambda}^{\infty}$, which intuitively constitute local perturbations to the true distribution $P_{n,0}^{\infty}$ of $\{X_i\}_{i=1}^\infty$. The sense in which $P_{n,\lambda}^\infty$ is local to $P_{n,0}^\infty$ is formalized by Assumption \ref{ass:localset}(ii), which imposes a requirement that is equivalent to the sequence $\{P_{n,\lambda}^{n}\}_{n=1}^\infty$ being contiguous to $\{P_{n,0}^{n}\}_{n=1}^\infty$ for any $\lambda \in \Lambda$ \citep{strasser1985mathematical}; see Remark \ref{rm:iidsimp}. In turn, Assumption \ref{ass:paramreg} contains our requirements on the parameter $\theta : \mathbf P^\infty \rightarrow \mathbb D_\phi$ and the estimator $\hat \theta_n$. Specifically, Assumption \ref{ass:paramreg}(i) demands that $\theta : \mathbf P^\infty \rightarrow \mathbb D_\phi$ be suitably smooth, while Assumption \ref{ass:paramreg}(ii) imposes that the distributional convergence of $\hat \theta_n$ be robust to local perturbations to $P^\infty$. As shown in \cite{Vaart1991differentibility}, these requirements are closely related, whereby Assumption \ref{ass:paramreg}(ii) and mild regularity conditions actually imply Assumption \ref{ass:paramreg}(i).

\begin{remark}\rm \label{rm:iidsimp}
Assumptions \ref{ass:localset} and \ref{ass:paramreg} encompass multiple standard constructions. For instance, suppose $\{X_i\}_{i=1}^n$ is an i.i.d. sequence, with $X_i \sim  P\in \mathbf P$ for some parametric model $\mathbf P$ -- e.g. $\mathbf P \equiv \{P_\beta : \beta \in \mathbf R^{d_\beta}\}$ and $P = P_{\beta_0}$ for some $\beta_0\in \mathbf R^{d_\beta}$. Assumption \ref{ass:localset}(i) is then satisfied by setting $\mathbf P^\infty \equiv \{\bigotimes_{i=1}^\infty P_\beta : \beta \in \mathbf R^{d_\beta}\}$ and employing the local parametrization $P_{n,\lambda}^\infty \equiv \bigotimes_{i=1}^\infty P_{\beta_0 + \lambda/\sqrt n}$, while Assumption \ref{ass:localset}(ii) holds provided the model is differentiable in quadratic mean. An analogous construction can also be employed to accommodate semiparametric or nonparametric models.\footnote{Specifically, in this instance $\Lambda$ corresponds to the tangent set of the model $\mathbf P$ at $P\in \mathbf P$; see an earlier version of this paper at \url{http://arxiv.org/pdf/1404.3763v1.pdf} for details.} It is also worth emphasizing that Assumption \ref{ass:localset} accommodates certain time series applications as well; see, e.g., \cite{bickel1998} and \cite{garel1995local}. \qed
\end{remark}

Given the stated assumptions, we can now conduct the desired local analysis. The following Lemma is a straightforward modification of Theorem \ref{th:delta}, and can be interpreted as extending the local analysis in \cite{Dumbgen1993} to contiguous perturbations.

\begin{lemma}\label{th:localdist}
If Assumptions \ref{ass:structure}, \ref{ass:paramest}, \ref{ass:localset}, and \ref{ass:paramreg} hold, then it follows that
\begin{equation}\label{th:localdistdisp}
r_n\{\phi(\hat \theta_n) - \phi(\theta(P^\infty_{n,\lambda}))\} \stackrel{L_\lambda}{\rightarrow} \phi^\prime_{\theta_0}(\mathbb G_0 + \theta^\prime(\lambda)) - \phi^\prime_{\theta_0}(\theta^\prime(\lambda))
\end{equation}
for any $\lambda \in \Lambda$ and $L_{\lambda}$ denoting convergence in distribution under the laws $P^\infty_{n,\lambda}$.
\end{lemma}

Lemma \ref{th:localdist} characterizes the asymptotic distribution of $\phi(\hat \theta_n)$ under a sequence of local perturbations $P_{n,\lambda}^\infty$ to the distribution $P^\infty$ of $\{X_i\}_{i=1}^\infty$. As expected, the local asymptotic distribution in \eqref{th:localdistdisp} need not equal the pointwise limit derived in Theorem \ref{th:delta}.  Intuitively, the asymptotic approximation in \eqref{th:localdistdisp} reflects the importance of local parameters and for this reason can be expected to provide a better approximation to finite sample distributions -- a point forcefully argued in the study of moment inequality models by \cite{andrews:soares:2010} and \cite{andrews:shi:2013}; see Remark \ref{rm:forunif} below. It is also worth noting that Lemma \ref{th:localdist} enables us to assess whether the estimator $\phi(\hat \theta_n)$ is regular -- i.e. whether the asymptotic distribution of $\phi(\hat \theta_n)$ is robust to the local perturbations $P_{n,\lambda}^\infty$. In fact, provided the index of perturbations $\Lambda$ is sufficiently ``rich", it is possible to show that $\phi(\hat \theta_n)$ is regular if and only if the standard bootstrap is consistent for $\phi(\hat \theta_n)$; see Remark \ref{rm:regular}. Thus, we can conclude that the failure of the standard bootstrap is an innate characteristic of irregular models.

\begin{remark}\label{rm:forunif} \rm
In Example \ref{ex:interbounds}, let $\{X_i\}_{i=1}^n$ be an i.i.d. sequence with $X_i \sim P$ and $\hat \theta_n = \frac{1}{n} \sum_i X_i$. Setting $P^\infty \equiv \bigotimes_{i=1}^\infty P$ and $\theta(P^\infty) \equiv \int x dP(x)$, Theorem \ref{th:delta} yields
\begin{equation}\label{rm:forunif1}
\sqrt n \{\phi(\hat \theta_n) - \phi(\theta(P^\infty))\} \stackrel{L}{\rightarrow} \left\{\begin{array}{ll} \mathbb G_0^{(j^*)} & \text{ if } \theta^{(1)}(P^\infty) \neq \theta^{(2)}(P^\infty) \\ \max\{\mathbb G_0^{(1)},\mathbb G_0^{(2)}\} & \text{ if } \theta^{(1)}(P^\infty) = \theta^{(2)}(P^\infty) \end{array} \right. ~,
\end{equation}
where $\mathbb G_0 = (\mathbb G_0^{(1)}, \mathbb G_0^{(2)})^\prime$ is a normal vector, and $j^{*} = \arg\max_{j\in\{1,2\}} \theta^{(j)}(P^\infty)$ (see \eqref{ex:interbounds2}). As argued in \cite{andrews:soares:2010}, the discontinuity of the pointwise asymptotic distribution in \eqref{rm:forunif1} reflects a poor approximation to the finite sample distribution which depends continuously on $\theta^{(1)}(P^\infty) - \theta^{(2)}(P^\infty)$. An asymptotic analysis local to a $P$ such that $\theta^{(1)}(P^\infty) = \theta^{(2)}(P^\infty)$, however, lets us address this problem. Specifically, for $P_{n,\lambda}^\infty$ satisfying $\theta(P^\infty_{n,\lambda}) = \theta(P^\infty) + \lambda/\sqrt n$ and $\lambda = (\lambda^{(1)}, \lambda^{(2)})^\prime \in \mathbf R^2$, Lemma \ref{th:localdist} implies
\begin{equation}\label{rm:forunif2}
\sqrt n\{\phi(\hat \theta_n) - \phi(\theta(P^\infty_{n,\lambda}))\} \stackrel{L_\lambda}{\rightarrow} \max\{\mathbb G_0^{(1)} + \lambda^{(1)}, \mathbb G_0^{(2)} + \lambda^{(2)}\} - \max\{\lambda^{(1)}, \lambda^{(2)}\} ~.
\end{equation}
Thus, by reflecting the importance of the ``slackness" parameter $\lambda$, \eqref{rm:forunif2} provides a better framework with which to evaluate the performance of our proposed procedure. \qed
\end{remark}

\begin{remark}\label{rm:regular} \rm
The conclusion of Lemma \ref{th:localdist} implies that $\phi(\hat \theta_n)$ is a regular estimator of the parameters $\phi(\theta_0)$ if and only if the distribution of
\begin{equation}\label{rm:regular1}
\phi^\prime_{\theta_0}(\mathbb G_0 + \theta^\prime(\lambda)) - \phi^\prime_{\theta_0}(\theta^\prime(\lambda))
\end{equation}
is constant in $\lambda \in \Lambda$. This ``invariance" requirement is closely related to a necessary and sufficient condition for the consistency of the standard bootstrap that can be derived from results in \cite{Dumbgen1993} (see \eqref{rm:dumbgen1}). In particular, if the closure of $\{\theta^\prime(\lambda) : \lambda \in \Lambda\}$ in $\mathbb D$ equals the support of $\mathbb G_0$, then it can be shown that $\phi(\hat \theta_n)$ is a regular estimator if and only if the standard bootstrap is consistent. Such conclusion complements \cite{beran:1997}, who showed that in finite dimensional likelihood models the parametric bootstrap is consistent if and only if the estimator is regular. \qed
\end{remark}

\subsubsection{Implications for Testing}

As has been emphasized in the moment inequalities literature, the lack of regularity of $\phi(\hat \theta_n)$ can render pointwise asymptotic approximations unreliable \citep{imbens:manski:2004}. However, since in Examples \ref{ex:interbounds}, \ref{ex:condineq}, and \ref{ex:stochdom} our results encompass procedures that are valid uniformly in the underlying distribution, we also know that irregularity of $\phi(\hat \theta_n)$ does not preclude our approach from remaining valid \citep{andrews:soares:2010, Linton2010, andrews:shi:2013}. In what follows, we note that the aforementioned examples are linked by the common structure of $\phi_{\theta_0}^\prime$ being subadditive. More generally, we exploit Lemma \ref{th:localdist} to show that whenever such property holds, the bootstrap procedure of Theorem \ref{th:modboot} locally controls size along a contiguous distributions. Thus, the potentially challenging task of verifying local size control can be reduced to a simple analytical computation of the directional derivative.

We consider hypothesis testing problems in which $\phi$ is scalar valued ($\mathbb E = \mathbf R$), and we are concerned with evaluating whether the distribution $P^\infty$ of $\{X_i\}_{i=1}^\infty$ satisfies
\begin{equation}\label{eq:hyptestdef}
H_0 : \phi(\theta(P^\infty)) \leq 0 \hspace{1 in } H_1 : \phi(\theta(P^\infty)) > 0 ~.
\end{equation}
A natural test statistic for this problem is then $ r_n\phi(\hat \theta_n)$, while Theorem \ref{th:delta} suggests
$$ c_{1-\alpha} \equiv \inf\{c : P(\phi_{\theta_0}^\prime(\mathbb G_0) \leq c ) \geq 1-\alpha \} $$
is an appropriate unfeasible critical value for a $1-\alpha$ level test.\footnote{Note that $c_{1-\alpha}$ is the $1-\alpha$ quantile of the asymptotic distribution of $r_n\phi(\hat \theta_n)$ when $\phi(\theta(P^\infty)) = 0$.} For $\hat c_{1-\alpha}$ the developed bootstrap estimator for $c_{1-\alpha}$ (see \eqref{cor:critvaldisp}), Theorem \ref{th:delta} and Corollary \ref{cor:critval} then establish the (pointwise in $P^\infty$) validity of rejecting $H_0$ whenever $r_n \phi(\hat \theta_n) > \hat c_{1-\alpha}$.

In order to evaluate both the local size and power of such a test, we conduct the analysis local to a $P^\infty$ that is in the ``boundary" of the null and alternative hypotheses -- i.e. we require that $\phi(\theta(P^\infty)) = 0$. This condition, and the additional requirements we impose for our local analysis, are formalized in the following assumption.

\begin{assumption}\label{ass:4test}
(i) $\mathbb E = \mathbf R$ and $\phi(\theta(P^\infty)) = 0$; (ii) The cdf of $\phi_{\theta_0}^\prime(\mathbb G_0)$ is continuous and strictly increasing at $c_{1-\alpha}$; (iii) $\phi_{\theta_0}^\prime(h_1 + h_2) \leq \phi_{\theta_0}^\prime(h_1) + \phi_{\theta_0}^\prime(h_2)$ for all $h_1,h_2 \in \mathbb D_0$.
\end{assumption}

Assumption \ref{ass:4test}(i) formalizes the requirements that $\phi$ be scalar valued and that $\phi(\theta(P^\infty))$ be equal to zero. In turn, in Assumption \ref{ass:4test}(ii) we impose that the cdf of $\phi_{\theta_0}^\prime(\mathbb G_0)$ be strictly increasing and continuous. Strict monotonicity is required to establish the consistency of $\hat c_{1-\alpha}$, while continuity ensures the test controls size at least pointwise in the distribution of $\{X_i\}_{i=1}^\infty$. Assumption \ref{ass:4test}(iii) demands that $\phi_{\theta_0}^\prime$ be subadditive, which represents the key condition that ensures local size control. Since $\phi_{\theta_0}^\prime$ is also positively homogenous of degree one, Assumption \ref{ass:4test}(iii) is in fact equivalent to demanding that $\phi_{\theta_0}^\prime$ be convex, which greatly simplifies verifying Assumption \ref{ass:4test}(ii) when $\mathbb G_0$ is Gaussian; see Remark \ref{rm:gausshelp}. We further note that Assumption \ref{ass:4test} is trivially satisfied when $\phi_{\theta_0}^\prime$ is linear, which by Lemma \ref{th:localdist} also implies $\phi(\hat \theta_n)$ is regular. However, we emphasize that Assumption \ref{ass:4test} can also hold at points $\theta_0$ at which $\phi$ is not Hadamard differentiable, as is easily verified in Examples \ref{ex:absmean}-\ref{ex:lrorder}.

The following Theorem employs the directional derivative $\phi_{\theta_0}^\prime$ to characterize the local asymptotic properties of a test that rejects the null hypothesis in \eqref{eq:hyptestdef} whenever the test statistic $r_n\phi(\hat \theta_n)$ exceeds the critical value $\hat c_{1-\alpha}$ (as in \eqref{cor:critvaldisp}).

\begin{theorem}\label{th:localsize}
If Assumptions \ref{ass:structure}, \ref{ass:paramest}, \ref{ass:primboot}-\ref{ass:paramreg}, and \ref{ass:4test}(i)-(ii) hold, then for any $\lambda \in \Lambda$
\begin{equation}\label{th:localsizedisp1}
\liminf_{n\rightarrow \infty} P_{n,\lambda}^n(r_n\phi(\hat \theta_n) > \hat c_{1-\alpha}) \geq P(\phi_{\theta_0}^\prime(\mathbb G_0 + \theta^\prime(\lambda)) > c_{1-\alpha}) ~.
\end{equation}
If in addition Assumption \ref{ass:4test}(iii) holds and $\phi(\theta(P_{n,\lambda}^\infty)) \leq 0$ for all $n$, then
\begin{equation}\label{th:localsizedisp2}
\limsup_{n\rightarrow \infty} P^n_{n,\lambda}(r_n \phi(\hat \theta_n) > \hat c_{1-\alpha}) \leq \alpha ~.
\end{equation}
\end{theorem}

The first claim of the Theorem derives a lower bound on the power function against local alternatives, with \eqref{th:localsizedisp1} holding with equality whenever $c_{1-\alpha}$ is a continuity point of the cdf of $\phi_{\theta_0}^\prime(\mathbb G_0 + \theta^\prime(\lambda))$. In turn, provided $\phi_{\theta_0}^\prime$ is subadditive, the second claim of Theorem \ref{th:localsize} establishes the ability of the test to locally control size along contiguous sequence. Heuristically, the role of subadditivity can be seen from \eqref{th:localsizedisp1} and noting
$$ P(\phi_{\theta_0}^\prime(\mathbb G_0 + \theta^\prime(\lambda)) > c_{1-\alpha}) \leq P(\phi_{\theta_0}^\prime(\mathbb G_0) + \phi_{\theta_0}^\prime(\theta^\prime(\lambda)) > c_{1-\alpha}) \leq \alpha ~,$$
where the final inequality results from $\phi_{\theta_0}^\prime(\theta^\prime(\lambda)) \leq 0$ provided the contiguous sequence $\{P_{n,\lambda}^\infty\}_{n=1}^\infty$ satisfies the null hypothesis -- i.e. $\phi(\theta(P_{n,\lambda}^\infty)) \leq 0$ for all $n$.\footnote{More precisely, we are exploiting that $\phi_{\theta_0}^\prime(\theta^\prime(\lambda)) = \lim_{n\rightarrow \infty} r_n \{\phi(\theta(P_{n,\lambda}^\infty)) - \phi(\theta(P^\infty))\} \leq 0$.} Thus, $\phi_{\theta_0}^\prime$ being subadditive implies $\lambda = 0$ is the ``least favorable" point in the null, which in turn delivers local size control as in \eqref{th:localsizedisp2}. We note a similar logic can be employed to evaluate confidence regions built using Theorems \ref{th:delta} and \ref{th:modboot}; see Remark \ref{rm:ci}.

Since the results of Theorem \ref{th:localsize} are local to a particular distribution $P^\infty$ of $\{X_i\}_{i=1}^\infty$, their relevance is contingent to them applying to all $P^\infty\in \mathbf P^\infty$ that are deemed possible distributions of the data. We emphasize that the three key requirements in this regard are Assumptions \ref{ass:paramreg}(ii), \ref{ass:4test}(ii), and \ref{ass:4test}(iii) -- i.e. that $\hat \theta_n$ be regular, the cdf of $\phi_{\theta_0}^\prime(\mathbb G_0)$ be continuous and strictly increasing at $c_{1-\alpha}$, and that $\phi_{\theta_0}^\prime$ be subadditive. We view Assumption \ref{ass:4test}(ii) as a technical requirement that can be dispensed with following arguments in \cite{andrews:shi:2013}; see Remark \ref{rm:contcorr}. Regularity of $\hat \theta_n$ and subadditivity of $\phi_{\theta_0}^\prime$, however, are instrumental in establishing the validity of our proposed procedure. In certain applications, such as in Examples \ref{ex:absmean}, \ref{ex:interbounds}, \ref{ex:condineq}, and \ref{ex:stochdom}, both these requirements are seen to be easily satisfied for a large class of possible distributions $P^\infty$. However, in other instances, such as in Example \ref{ex:suppfunction} applied to the estimator in \cite{Kaido_Santos2013}, $\phi_{\theta_0}^\prime$ is always subadditive, but the regularity of $\hat \theta_n$ can fail to hold for an important class of distributions $P^\infty$.

\begin{remark}\label{rm:ci} \rm
We can obtain confidence regions for $\phi(\theta(P^\infty))$ through test inversion of
\begin{equation}\label{eq:rmci1}
H_0 : \phi(\theta(P^\infty)) = c_0 \hspace{1 in} H_1: \phi(\theta(P^\infty)) \neq c_0 ~,
\end{equation}
for different values of $c_0 \in \mathbb E$. Defining $\bar \phi : \mathbb D_\phi \subseteq \mathbb D \rightarrow \mathbf R$ pointwise by $\bar \phi(\theta) \equiv \|\phi(\theta) - c_0\|_{\mathbb E}$, it is then straightforward to see \eqref{eq:rmci1} can be expressed as in \eqref{eq:hyptestdef} with $\bar \phi$ in place of $\phi$. In particular, the chain rule implies $\bar \phi_{\theta_0}^\prime(\cdot) = \|\phi_{\theta_0}^\prime(\cdot)\|_{\mathbb E}$, and hence the subadditivity of $\|\phi_{\theta_0}^\prime(\cdot)\|_{\mathbb E}$ suffices for establishing local size control. \qed
\end{remark}

\begin{remark}\label{rm:gausshelp} \rm
Under Assumptions \ref{ass:structure} and \ref{ass:4test}(iii), it follows that $\phi_{\theta_0}^\prime : \mathbb D_0 \rightarrow \mathbf R$ is a continuous convex functional. Therefore, if $\mathbb G_0$ is in addition Gaussian, then Theorem 11.1 in \cite{davydov:lifshits:smorodina:1998} implies that the cdf of $\phi_{\theta_0}^\prime(\mathbb G_0)$ is continuous and strictly increasing at all points in the interior of its support (relative to $\mathbf R$). \qed
\end{remark}

\begin{remark}\label{rm:contcorr} \rm
In certain applications, such as in Examples \ref{ex:condineq} and \ref{ex:stochdom}, Assumption \ref{ass:4test}(ii) may be violated at distributions $P^\infty$ of interest. To address this problem, \cite{andrews:shi:2013} propose employing the critical value $\hat c_{1-\alpha} + \delta$ for an arbitrarily small $\delta > 0$. It is then possible to show that, even if Assumption \ref{ass:4test}(ii) fails,  we still have
\begin{equation}\label{rm:contcorr1}
\liminf_{n\rightarrow \infty} P(\hat c_{1-\alpha} + \delta \geq c_{1-\alpha})  = 1 ~.
\end{equation}
Therefore, by contiguity it follows that the local size control established in \eqref{th:localsizedisp2} holds without Assumption \ref{ass:4test}(ii) if we employ $\hat c_{1-\alpha} + \delta$ instead of $\hat c_{1-\alpha}$. \qed
\end{remark}

\section{Convex Set Projections}\label{sec:convex}

In this section, we demonstrate the usefulness of our results by constructing a test of whether a Hilbert space valued parameter belongs to a known closed convex set -- a setting that encompasses tests of moment inequalities, shape restrictions, and the validity of random utility models. Despite the generality of the problem, we show our results enable us to develop a valid test by relying on purely analytical calculations.

\subsection{Projection Setup}

In what follows, we let $\mathbb H$ be a Hilbert space with inner product $\langle \cdot,\cdot\rangle_{\mathbb H}$ and norm $\|\cdot\|_{\mathbb H}$. For a known closed convex set $\Lambda \subseteq \mathbb H$, we then consider the hypothesis testing problem
\begin{equation}\label{conv1}
H_0: \theta_0 \in \Lambda \hspace{0.5 in } H_1: \theta_0 \notin \Lambda ~,
\end{equation}
where the parameter $\theta_0 \in \mathbb H$ is unknown, but for which we possess an estimator $\hat \theta_n$. Special cases of this problem have been widely studied in the setting where $\mathbb H = \mathbf R^d$, and to a lesser extent when $\mathbb H$ is infinite dimensional; see Examples \ref{ex:convex:finite}-\ref{ex:convex:constant} below.

We formalize the introduced structure through the following assumption.

\begin{assumption}\label{ass:convexsetup}
(i) $\mathbb D = \mathbb H$ where $\mathbb H$ is Hilbert space with inner product $\langle \cdot,\cdot\rangle_{\mathbb H}$ and corresponding norm $\|\cdot\|_{\mathbb H}$; (ii) $\Lambda\subseteq \mathbb H$ is a known closed and convex set.
\end{assumption}

Since projections onto closed convex sets in Hilbert spaces are attained and unique, we may define the projection operator $\Pi_{\Lambda} : \mathbb H \rightarrow \Lambda$, which for each $\theta \in \mathbb H$ satisfies
\begin{equation}\label{conv2}
\|\theta - \Pi_{\Lambda}\theta\|_{\mathbb H} = \inf_{h \in \Lambda}\|\theta - h\|_{\mathbb H} ~.
\end{equation}
Thus, the hypothesis testing problem in \eqref{conv1} can be rewritten in terms of the distance between $\theta_0$ and $\Lambda$, or equivalently between $\theta_0$ and its projection $\Pi_\Lambda \theta_0$ -- i.e.
\begin{equation}\label{conv3}
H_0: \|\theta_0 - \Pi_{\Lambda}\theta_0\|_{\mathbb H} = 0 \hspace{0.5 in } H_1: \|\theta_0 - \Pi_{\Lambda}\theta_0\|_{\mathbb H} > 0 ~.
\end{equation}
Interpreted in this manner, it is clear that \eqref{conv3} is a special case of \eqref{eq:hyptestdef}, with $\mathbb D = \mathbb H$, $\mathbb E = \mathbf R$, and $\phi : \mathbb H \rightarrow \mathbf R$ given by $\phi(\theta) \equiv \|\theta - \Pi_\Lambda \theta\|_{\mathbb H}$ for any $\theta \in \mathbb H$. The corresponding test statistic $r_n\phi(\hat \theta_n)$ is then simply the scaled distance between the estimator $\hat \theta_n$ and the known convex set $\Lambda$ -- i.e. $r_n\phi(\hat \theta_n) = r_n\|\hat \theta_n - \Pi_{\Lambda}\hat \theta_n\|_{\mathbb H}$.

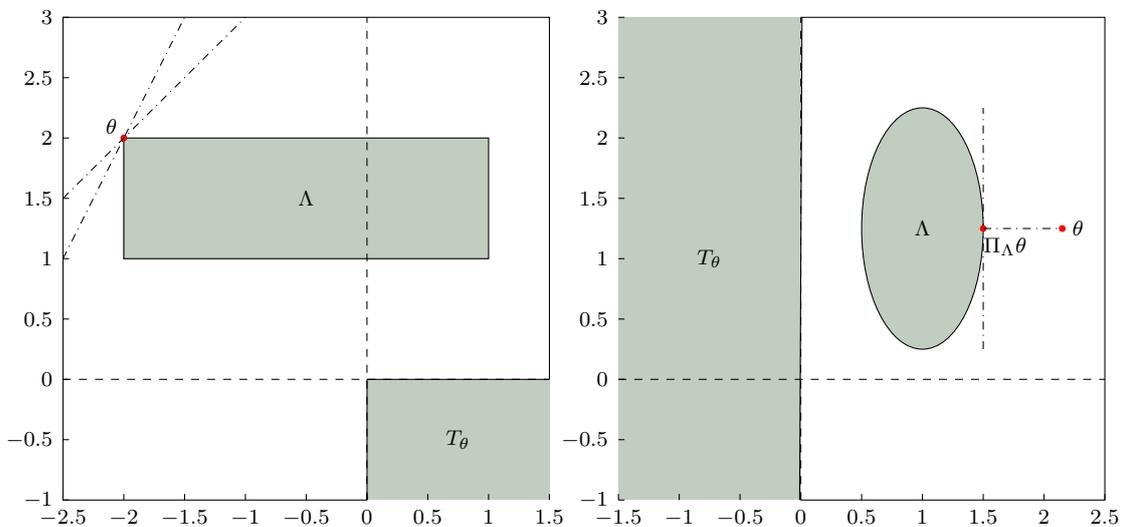
\begin{figure}[t]
\caption{Illustrations of Tangent Cones.}
\flushleft
\begin{minipage}[b]{0.35\linewidth}
\begin{tikzpicture}[scale=1.6]
\draw (0,-1)--(-2.5,-1)--(-2.5,3)--(1.5,3)--(1.5,0); 
\shadedraw[left color=Honeydew3,right color=Honeydew3] (-2,1)rectangle (1,2); 
\draw (-0.5,1.5) node {\scriptsize $\Lambda$};
\fill[red] (-2,2) circle (0.8pt);
\draw (-2.1,2.1) node {\scriptsize $\theta$}; 
\draw[dashdotted] (-2.5,1.5)--(-1,3); 
\draw[dashdotted] (-2.5,1)--(-1.5,3);
\fill[color=Honeydew3] (0,0) rectangle (1.5,-1); 
\draw (0.75,-0.5) node {\scriptsize $T_{\theta}$};
\draw (0,0)--(1.5,0);
\draw (0,0)--(0,-1);
\draw[dashed] (-2.5cm,0)--(1.5cm,0); 
\draw[dashed] (0,-1cm)--(0,3cm);
\foreach \x in {-2.5,-2,...,1.5}
\draw (\x cm,-0.95cm) -- (\x cm,-1cm) node[anchor=north] {\scriptsize $\x$};
\foreach \y in {-1,-0.5,...,3}
\draw (-2.45cm,\y cm) -- (-2.5cm,\y cm) node[anchor=east] {\scriptsize $\y$};
\end{tikzpicture}
\end{minipage}
\qquad \qquad ~ ~
\begin{minipage}[b]{0.35\linewidth}
\begin{tikzpicture}[scale=1.6]
\draw (0,-1)--(2.5,-1)--(2.5,3)--(0,3); 
\shadedraw[left color=Honeydew3,right color=Honeydew3]  (1,1.25) node {\scriptsize $\Lambda$} ellipse [x radius=0.5cm,y radius=1cm]; 
\draw[dashdotted] (1.5,0.25)--(1.5,2.25);
\draw (1.42,1.1) node[anchor=west] {\scriptsize $\Pi_{\Lambda}\theta$}; 
\fill[red] (1.5,1.25) circle (0.8pt);
\draw (2.15,1.25) node[anchor=west] {\scriptsize $\theta$}; 
\fill[red] (2.15,1.25) circle (0.8pt);
\draw[dashdotted] (1.5,1.25)--(2.15,1.25);
\fill[color=Honeydew3] (-1.5,-1) rectangle (0,3); 
\draw (-0.01,-1cm)--(0.01,3cm);
\draw (-0.75,1) node {\scriptsize $T_{\theta}$};
\draw[dashed] (-1.5cm,0)--(2.5cm,0); 
\draw[dashed] (0,-1cm)--(0,3cm);
\foreach \x in {-1.5,-1,...,2.5}
\draw (\x cm,-0.95 cm) -- (\x cm,-1 cm) node[anchor=north] {\scriptsize $\x$};
\foreach \y in {-1,-0.5,...,3}
\draw (-1.45 cm,\y cm) -- (-1.5 cm,\y cm) node[anchor=east] {\scriptsize $\y$};
\end{tikzpicture}
\end{minipage}
\end{figure}

As a final piece of notation, we need to introduce the tangent cone of $\Lambda$ at a $\theta \in \mathbb H$, which plays a fundamental role in our analysis. To this end, for any set $A \subseteq \mathbb H$ let $\overline{A}$ denote its closure under $\|\cdot\|_{\mathbb H}$, and define the tangent cone of $\Lambda$ at $\theta \in \mathbb H$ by
\begin{equation}\label{conv4}
T_{\theta} \equiv \overline{ \bigcup_{\alpha \geq 0} \alpha \{\Lambda - \Pi_{\Lambda} \theta\}} ~,
\end{equation}
which is convex by convexity of $\Lambda$. Heuristically, $T_{\theta}$ represents the directions from which the projection $\Pi_{\Lambda}\theta \in \Lambda$ can be approached from within the set $\Lambda$. As such, $T_{\theta}$ can be seen as a local approximation to the set $\Lambda$ at $\Pi_{\Lambda}\theta$ and employed to study the differentiability properties of the projection operator $\Pi_{\Lambda}$. Figure 4.1 illustrates the tangent cone in two separate cases: one in which $\theta \in \Lambda$, and a second in which $\theta \notin \Lambda$.

\subsubsection{Examples}

In order to aid exposition and illustrate the applicability of \eqref{conv1}, we next provide examples of both well studied and new problems that fit our framework.

\begin{example} \label{ex:convex:finite} \rm
Suppose $X \in \mathbf R^d$ and that we aim to test the moment inequalities
\begin{equation}\label{ex:convex:finite1}
H_0 : E[X] \leq 0 \hspace{0.5 in} H_1 : E[X] \nleq 0 ~,
\end{equation}
where the null is meant to hold at all coordinates, and the alternative indicates at least one coordinate of $E[X]$ is strictly positive. In this instance, $\mathbb H = \mathbf R^d$, $\Lambda$ is the negative orthant in $\mathbf R^d$ $(\Lambda \equiv \{h \in \mathbf R^d : h \leq 0\})$, and the distance of $\theta$ to $\Lambda$ is equal to
\begin{equation}\label{ex:convex:finite2}
\phi(\theta) = \|\Pi_{\Lambda}\theta - \theta\|_{\mathbb H} = \Big\{ \sum_{i=1}^d (E[X^{(i)}])_+^2 \Big\}^{\frac{1}{2}} ~,
\end{equation}
where $(a)_+ = \max\{a,0\}$ and $X^{(i)}$ denotes the $i^{th}$ coordinate of $X$. More generally, the hypothesis in \eqref{conv3} accommodates any regular parameter and any closed convex set in $\mathbf R^d$, such as the test for
moment inequalities on regression coefficients proposed by \cite{wolak:1988} and the test of random utility models developed in \cite{kitamura:stoye:2013}. Analogously, conditional moment inequalities as in Example \ref{ex:condineq} can also be encompassed by employing a weight function on $\mathcal F$ -- this approach leads to the Cramer-von-Mises statistic studied in \cite{andrews:shi:2013}. \qed
\end{example}

The next example is new and concerns quantile models, as employed by \cite{buchinsky:1994} to characterize the U.S. wage structure conditional on levels of education, or by \cite{abadie:angrist:imbens} to estimate the effect of subsidized training on earnings.

\begin{example} \label{ex:convex:monotone} \rm
Let $(Y,D,Z) \in \mathbf R \times \mathbf R \times \mathbf R^{d_z}$ and consider the quantile regression
\begin{equation}\label{ex:convex:monotone1}
(\theta_0(\tau),\beta(\tau)) \equiv \arg\min_{\theta \in \mathbf R, \beta \in \mathbf R^{d_z}} E[\rho_\tau(Y - D\theta - Z^\prime \beta)]
\end{equation}
where $\rho_\tau(u) = (\tau - 1\{u \leq 0\})u$ and $\tau \in (0,1)$. Under appropriate restrictions, the estimator $\hat \theta_n$ for $\theta_0$ converges in distribution in $\ell^\infty([\epsilon,1-\epsilon])$ for any $\epsilon \in(0,1/2)$ \citep{angrist:Chernozhukov:fernanez}.\footnote{This result also holds for the instrumental variables estimator of \cite{chernozhukov:hansen:2005}.} In this setting, it is often of interest to test for shape restrictions on $\theta_0$, which we may accomplish by setting $\mathbb H$ to equal the Hilbert space
\begin{equation}\label{ex:convex:monotone2}
\mathbb H \equiv \{\theta : [\epsilon, 1-\epsilon] \rightarrow \mathbf R : \langle \theta, \theta \rangle_{\mathbb H} < \infty\} \hspace{0.5 in} \langle \theta_1,\theta_2\rangle_{\mathbb H} \equiv \int_\epsilon^{1-\epsilon} \theta_1(\tau)\theta_2(\tau)d\tau ~,
\end{equation}
and considering appropriate convex sets $\Lambda \subseteq \mathbb H$. For instance, in randomized experiments where $D$ is a dummy for treatment, $\theta_0(\tau)$ is the quantile treatment effect and we may test for its constancy or monotonicity; see \cite{karthik2011} for an examination of these features in the evaluation of teacher performance pay. A similar approach may also be employed to test whether the pricing kernel satisfies theoretically predicted restrictions such as a monotonicity \citep{jackwerth:2000}. \qed
\end{example}

Our final example may be interpreted as a generalization of Example \ref{ex:convex:monotone}.

\begin{example} \label{ex:convex:constant} \rm
Let $Z \in \mathbf R^{d_z}$, $\Theta \subseteq \mathbf R^{d_\theta}$, and $\mathcal T \subseteq \mathbf R^{d_\tau}$. Suppose there exists a function $\rho: \mathbf R^{d_z}\times \Theta \times \mathcal T\rightarrow \mathbf R^{d_\rho}$ such that for each $\tau \in \mathcal T$ there is a unique $\theta_0(\tau) \in \Theta$ satisfying
\begin{equation}\label{ex:convex:constant1}
E[\rho(Z,\theta_0(\tau), \tau)] = 0 ~.
\end{equation}
Such a setting arises, for instance, in sensitivity analysis \citep{chen:tamer:torgo:2011}, and in partially identified models where the identified set is a curve \citep{arellano:hansen:sentana:2012} or can be described by a functional lower and upper bound \citep{kline:santos:2013, chandrasekhar:chernozhukov:molinari:schrimpf}. \cite{escanciano:zhu:2013} derives an estimator $\hat \theta_n$ which converges in distribution in $\bigotimes_{i=1}^{d_\theta}\ell^\infty(\mathcal T)$, and hence for an integrable function $w$ also in
\begin{equation}\label{ex:convex:constant2}
\mathbb H \equiv \{\theta : \mathcal T \rightarrow \mathbf R^{d_\theta} : \langle \theta, \theta \rangle_{\mathbb H} < \infty\} \hspace{0.5 in} \langle \theta_1,\theta_2\rangle_{\mathbb H} \equiv \int_{\mathcal T} \theta_1(\tau)^\prime \theta_2(\tau)w(\tau)d\tau ~.
\end{equation}
Appropriate choices of $\Lambda$ then enable us to test, for example, whether the model is identified in \cite{arellano:hansen:sentana:2012}, or whether the identified set in \cite{kline:santos:2013} is consistent with increasing returns to education across quantiles. \qed
\end{example}

\subsection{Theoretical Results}

\subsubsection{Asymptotic Distribution}

Our analysis crucially relies on the seminal work of \cite{zarantonello}, who established the Hadamard directional differentiability of metric projections onto convex sets in Hilbert spaces. Specifically, \cite{zarantonello} showed $\Pi_\Lambda : \mathbb H \rightarrow \Lambda$ is Hadamard directionally differentiable at any $\theta \in \Lambda$, and its directional derivative is equal to the projection operator onto the tangent cone of $\Lambda$ at $\theta$, which we denote by $\Pi_{T_\theta} : \mathbb H \rightarrow T_{\theta}$. Figure 2 illustrates a simple example in which the derivative approximation
\begin{equation}\label{convth1}
\Pi_{\Lambda}\theta_1 - \Pi_{\Lambda} \theta_0 \approx \Pi_{T_{\theta_0}}(\theta_1 - \theta_0)
\end{equation}
actually holds with equality. We note that it is also immediate from Figure 2 that the directional derivative $\Pi_{T_{\theta_0}}$ is not linear, and hence $\Pi_{\Lambda}$ is not fully differentiable.\footnote{For related work where the projection is assumed to be differentiable but $\|\cdot\|_{\mathbb H}$ is not required to be Hilbertian, see \cite{romano1988bootstrap} who studies nonparametric distance tests for empirical distributions.}

Given the result in \cite{zarantonello}, the asymptotic distribution of $r_n\phi(\hat \theta_n)$ can then be obtained as an immediate consequence of Theorem \ref{th:delta}.

\begin{proposition}\label{prop:convexdist}
Let Assumption \ref{ass:paramest} and \ref{ass:convexsetup} hold. If $\theta_0 \in \Lambda$, then it follows that
\begin{equation}\label{prop:convexdistdisp}
r_n\|\hat \theta_n - \Pi_{\Lambda} \hat \theta_n\|_{\mathbb H} \stackrel{L}{\rightarrow} \|\mathbb G_0 - \Pi_{T_{\theta_0}}\mathbb G_0\|_{\mathbb H} ~.
\end{equation}
\end{proposition}

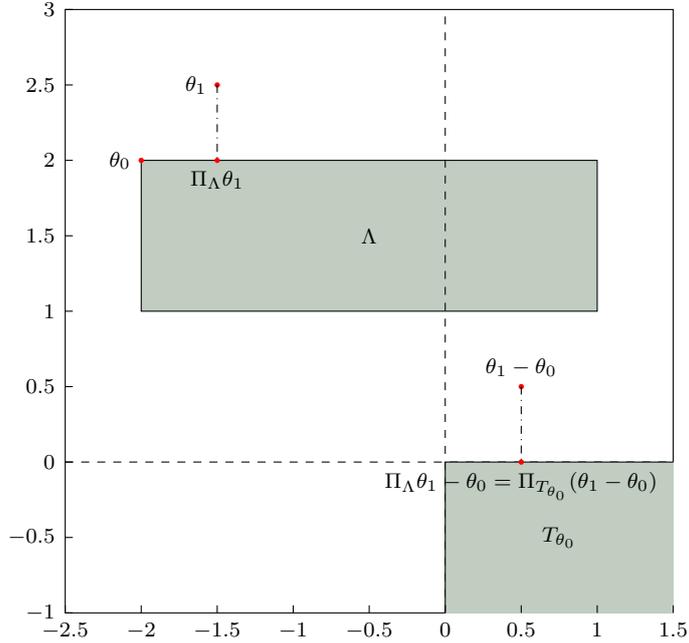
\begin{figure}
\caption{Directional Differentiability}
\centering
\begin{tikzpicture}[scale=2]
\shadedraw[left color=Honeydew3,right color=Honeydew3]  (-2,1) rectangle (1,2); 
\draw (-0.5,1.5) node {\scriptsize $\Lambda$};
\fill[red] (-2,2) circle (0.5pt) node[anchor=east] {\scriptsize\color{black} $\theta_0$}; 
\fill[red] (-1.5,2.5) circle (0.5pt) node[anchor=east] {\scriptsize\color{black} $\theta_1$};
\draw[dashdotted] (-1.5,2.5)--(-1.5,2) node[anchor=north] {\scriptsize $\Pi_{\Lambda}\theta_1$};
\fill[red] (-1.5,2) circle (0.5pt);
\fill[red] (0.5,0.5) circle (0.5pt) node[anchor=south] {\scriptsize\color{black} $\theta_1-\theta_0$};
\draw[dashdotted] (0.5,0.5)--(0.5,0);
\draw[dashed] (-2.5,0)--(1.5,0); 
\draw[dashed] (0,-1)--(0,3);
\fill[color=Honeydew3](0,0)rectangle (1.5,-1); 
\draw (0,-1)--(0,0)--(1.5,0);
\draw (0.75,-0.5) node {\scriptsize $T_{\theta_0}$};
\fill[red] (0.5,0) circle (0.5pt) node[anchor=north] {\scriptsize \color{black} $\Pi_{\Lambda}\theta_1-\theta_0=\Pi_{T_{\theta_0}}(\theta_1-\theta_0)$};
\draw (0,-1)--(-2.5,-1)--(-2.5,3)--(1.5,3)--(1.5,0);
\foreach \x in {-2.5,-2,...,1.5}
\draw (\x cm,-0.95cm) -- (\x cm,-1 cm) node[anchor=north] {\scriptsize $\x$};
\foreach \y in {-1,-0.5,...,3}
\draw (-2.45cm,\y cm) -- (-2.5cm,\y cm) node[anchor=east] {\scriptsize $\y$};
\end{tikzpicture}
\end{figure}


In particular, Proposition \ref{prop:convexdist} follows from norms being directionally differentiable at zero, and hence by the chain rule the directional derivative $\phi_{\theta_0}^\prime : \mathbb H \rightarrow \mathbf R$ satisfies
\begin{equation}\label{convth2}
\phi_{\theta_0}^\prime(h) = \|h - \Pi_{T_{\theta_0}} h\|_{\mathbb H} ~.
\end{equation}
It is interesting to note that $\Lambda \subseteq T_{\theta_0}$ whenever $\Lambda$ is a cone, and hence $\|h - \Pi_{T_{\theta_0}} h\|_{\mathbb H} \leq \|h - \Pi_\Lambda h\|_{\mathbb H}$ for all $h \in \mathbb H$. Therefore, the distribution of $\|\mathbb G_0 - \Pi_{\Lambda}\mathbb G_0\|_{\mathbb H}$ first order stochastically dominates that of $\|\mathbb G_0 - \Pi_{T_{\theta_0}} \mathbb G_0\|_{\mathbb H}$, and by Proposition \ref{prop:convexdist} its quantiles may be employed for potentially conservative inference -- an approach that may be viewed as analogous to assuming all moments are binding in the moment inequalities literature. It is also worth noting that Proposition \ref{prop:convexdist} can be readily extended to study the projection itself rather than its norm, allow for nonconvex sets $\Lambda$, and incorporate weight functions into the test statistic; see Remarks \ref{rm:notconvex} and \ref{rm:weights}.

\begin{remark} \label{rm:notconvex} \rm
\cite{zarantonello} and Theorem \ref{th:delta} can be employed to derive the asymptotic distribution of the projection $r_n\{\Pi_{\Lambda}\hat \theta_n - \Pi_\Lambda \theta_0\}$ itself. However, when studying the projection, it is natural to aim to relax the requirement that $\theta_0 \in \Lambda$. Such an extension, as well as considering non-convex $\Lambda$, is possible under appropriate regularity conditions -- see \cite{shapiro:1994} for the relevant directional differentiability results. \qed
\end{remark}

\begin{remark} \label{rm:weights} \rm
While we do not consider it for simplicity, it is straightforward to incorporate weight functions into the test statistic.\footnote{For instance in \eqref{ex:convex:finite2} we may wish to consider $\{ \sum_{i=1}^d (E[X^{(i)}])_+^2/Var(X^{(i)})\}^{\frac{1}{2}}$ instead.} Formally, a weight function may be seen as a linear operator $A:\mathbb H \rightarrow \mathbb H$, and for any estimator $\hat A_n$ such that $\|\hat A_n - A\|_o = o_p(1)$ for $\|\cdot\|_o$ the operator norm, we obtain by asymptotic tightness of $r_n\{\hat \theta_n - \Pi_{\Lambda} \hat \theta_n \}$ that
\begin{equation}\label{rm:weights1}
r_n\|\hat A_n \{\hat \theta_n - \Pi_{\Lambda} \hat \theta_n\}\|_{\mathbb H} \stackrel{L}{\rightarrow} \| A \{\mathbb G_0 - \Pi_{T_{\theta_0}} \mathbb G_0\}\|_{\mathbb H} ~.
\end{equation}
Thus, estimating weights has no first order effect on the asymptotic distribution. \qed
\end{remark}

\subsubsection{Critical Values}

In order to construct critical values to conduct inference, we next aim to employ Theorem \ref{th:modboot}, which requires the availability of a suitable estimator $\hat \phi_n^\prime$ for the directional derivative $\phi_{\theta_0}^\prime$. To this end, we develop an estimator $\hat \phi_n^\prime$ which, despite being computationally intensive, is guaranteed to satisfy Assumption \ref{ass:derest} under no additional requirements.

Specifically, for an appropriate $\epsilon_n \downarrow 0$, we define $\hat \phi_n^\prime :\mathbb H \rightarrow \mathbf R$ pointwise in $h\in \mathbb H$ by
\begin{equation}\label{convth3}
\hat \phi_n^\prime(h) \equiv \sup_{\theta \in \Lambda : \|\theta - \Pi_{\Lambda} \hat \theta_n\|_{\mathbb H} \leq \epsilon_n } \|h - \Pi_{T_\theta} h\|_{\mathbb H} ~.
\end{equation}
Heuristically, we estimate $\phi_{\theta_0}^\prime(h)= \|h - \Pi_{T_{\theta_0}}h\|_{\mathbb H}$ by the distance between $h$ and the ``least favorable" tangent cone $T_{\theta}$ that can be generated by the $\theta\in \Lambda$ that are in a neighborhood of $\Pi_{\Lambda} \hat \theta_n$. It is evident from this construction that provided $\epsilon_n \downarrow 0$ at an appropriate rate, the shrinking neighborhood of $\Pi_{\Lambda} \hat \theta_n$ will include $\theta_0$ with probability tending to one and as a result $\hat \phi^\prime_n(h)$ will provide a potentially conservative estimate of $\phi_{\theta_0}^\prime(h)$. As the following Proposition shows, however, $\hat \phi_n^\prime(h)$ is in fact not conservative, and $\hat \phi_n^\prime$ provides a suitable estimator for $\phi_{\theta_0}^\prime$ in the sense required by Theorem \ref{th:modboot}.

\begin{proposition}\label{prop:convexinf}
Let Assumptions \ref{ass:paramest}, \ref{ass:convexsetup} hold, and $\phi_{\theta_0}^\prime(h) \equiv \|h- \Pi_{T_{\theta_0}}h\|_{\mathbb H}$. Then,
\begin{enumerate}
    \item[(i)] \vspace{-0.1 in} If $\epsilon_n \downarrow 0$ and $\epsilon_n r_n \uparrow \infty$, then $\hat \phi_n^\prime$ as defined in \eqref{convth3} satisfies Assumption \ref{ass:derest}.
    \item[(ii)] \vspace{-0.05 in} $\phi_{\theta_0}^\prime:\mathbb H \rightarrow \mathbf R$ satisfies $\phi_{\theta_0}^\prime(h_1 + h_2) \leq \phi_{\theta_0}^\prime(h_1) + \phi_{\theta_0}^\prime(h_2)$ for all $h_1,h_2 \in \mathbb H$.
\end{enumerate}
\end{proposition}

The first claim of the Proposition shows that $\hat \phi_n^\prime$ satisfies Assumption \ref{ass:derest}. Therefore, provided the bootstrap is consistent for the asymptotic distribution of $r_n\{\hat \theta_n - \theta_0\}$, Theorem \ref{th:modboot} implies we may employ the $1-\alpha$ quantile (conditional on $\{X_i\}_{i=1}^n$) of
\begin{equation}\label{convth3}
\sup_{\theta \in \Lambda : \|\theta - \Pi_{\Lambda} \hat \theta_n\|_{\mathbb H} \leq \epsilon_n} \|r_n\{\hat \theta_n^* - \hat \theta_n\} - \Pi_{T_{\theta}} \{r_n\{\hat \theta_n^* - \hat \theta_n\}\}\|_{\mathbb H}
\end{equation}
as critical values. In turn, Proposition \ref{prop:convexinf}(ii) establishes that the directional derivative $\phi_{\theta_0}^\prime$ is always subadditive. Thus, one of the key requirement of Theorem \ref{th:localsize} is satisfied, and we can conclude that employing the quantiles of \eqref{convth3} (conditional on $\{X_i\}_{i=1}^n$) as critical values will deliver local size control whenever $\hat \theta_n$ is regular. This latter conclusion continues to hold if an alternative consistent estimator to \eqref{convth3} is employed to construct critical values. Hence, we emphasize that while $\hat \phi_n^\prime$ as defined in \eqref{convth3} is appealing due to its general applicability, its use may not be advisable in instances where simpler estimators of $\phi_{\theta_0}^\prime$ are available; see Remark \ref{rm:tangentineq}.

\begin{remark}\label{rm:tangentineq} \rm
In certain applications, the tangent cone $T_{\theta_0}$ can be easily estimated and as a result so can $\phi_{\theta_0}^\prime$. For instance, in the moment inequalities model of Example \ref{ex:convex:finite},
\begin{equation}\label{eq:tangentineq1}
T_{\theta_0} = \{h \in \mathbf R^d : h^{(i)} \leq 0 \text{ for all } i \text{ such that } E[X^{(i)}] = 0\} ~.
\end{equation}
For $\bar X$ the mean of an i.i.d. sample $\{X_i\}_{i=1}^n$, a natural estimator for $T_{\theta_0}$ is then given by
\begin{equation}\label{eq:tangentineq2}
\hat T_n = \{h \in \mathbf R^d : h^{(i)} \leq 0 \text{ for all } i \text{ such that } \bar X^{(i)} \geq -\epsilon_n\} ~
\end{equation}
for some sequence $\epsilon_n \downarrow 0$ satisfying $\epsilon_n \sqrt n \uparrow \infty$. It is then straightforward to verify that $\hat \phi_n^\prime(h) = \|h - \Pi_{\hat T_n} h\|_{\mathbb H}$ satisfies Assumption \ref{ass:derest} (compare to \eqref{convth2}) and, more interestingly, that the bootstrap procedure of Theorem \ref{th:modboot} then reduces to the generalized moment selection approach of \cite{andrews:soares:2010}. \qed
\end{remark}

\subsection{Simulation Evidence}

In order to examine the finite sample performance of the proposed test and illustrate its implementation, we next conduct a limited Monte Carlo study based on Example \ref{ex:convex:monotone}. Specifically, we consider a quantile treatment effect model in which the treatment dummy $D \in \{0,1\}$ satisfies $P(D = 1) = 1/2$, the covariates $Z = (1,Z^{(1)},Z^{(2)})^\prime \in \mathbf R^3$ satisfy $(Z^{(1)},Z^{(2)})^\prime \sim N(0,I)$ for $I$ the identity matrix, and $Y$ is related by
\begin{equation}\label{mc1}
Y = \frac{\Delta}{\sqrt n} D \times U + Z^\prime \beta + U ~,
\end{equation}
where $\beta = (0,1/\sqrt 2,1/\sqrt 2)^\prime$ and $U$ is unobserved, uniformly distributed on $[0,1]$, and independent of $(D,Z)$. It is then straightforward to verify that $(Y,D,Z)$ satisfy
\begin{equation}\label{mc2}
P(Y \leq D\theta_0(\tau) + Z^\prime \beta(\tau)|D,Z) = \tau ~,
\end{equation}
for $\theta_0(\tau) \equiv \tau \Delta/\sqrt n$ and $\beta(\tau) \equiv (\tau, 1/\sqrt 2, 1/\sqrt 2)^\prime$. Hence, in this context the quantile treatment effect has been set local to zero at all $\tau$, which enables us to evaluate the local power and local size control of the proposed test.

\begin{table}[t]
\begin{footnotesize}
\caption{Empirical Size}
\begin{center}
\begin{tabular}{cc ccc c ccc c ccc}
\hline
\hline
& &  \multicolumn{11}{c}{$n = 200$} \\ \cline{3-13}
\multicolumn{2}{c}{Bandwidth} &  \multicolumn{3}{c}{$\alpha = 0.1$} & & \multicolumn{3}{c}{$\alpha = 0.05$} & & \multicolumn{3}{c}{$\alpha = 0.01$} \\ \cline{3-5} \cline{7-9} \cline{11-13}
$C$ & $\kappa$       & $\Delta = 0$  & $\Delta = 1$ & $\Delta = 2$ & & $\Delta = 0$ & $\Delta = 1$ & $\Delta = 2$ & & $\Delta = 0$ & $\Delta = 1$ & $\Delta = 2$ \\
\hline
1    & 1/4 & 0.042  & 0.017 & 0.006 & & 0.020 & 0.008 & 0.002 & & 0.005 & 0.001 & 0.000 \\
1    & 1/3 & 0.042	& 0.017	& 0.006	& & 0.020 &	0.008 &	0.002 & & 0.005	& 0.001	& 0.000 \\
0.01 & 1/4 & 0.082	& 0.053	& 0.035	& & 0.035 &	0.023 &	0.013 & & 0.007	& 0.002	& 0.001 \\
0.01 & 1/3 & 0.087	& 0.059	& 0.042	& & 0.038 &	0.025 &	0.015 & & 0.007	& 0.002	& 0.001 \\
\multicolumn{2}{c}{Theoretical}
           & 0.100  & 0.042 & 0.015 & & 0.050 & 0.018 & 0.006 & & 0.010 & 0.003 & 0.001 \\ \\

& &  \multicolumn{11}{c}{$n = 500$} \\ \cline{3-13}
\multicolumn{2}{c}{Bandwidth} &  \multicolumn{3}{c}{$\alpha = 0.1$} & & \multicolumn{3}{c}{$\alpha = 0.05$} & & \multicolumn{3}{c}{$\alpha = 0.01$} \\ \cline{3-5} \cline{7-9} \cline{11-13}
$C$ & $\kappa$       & $\Delta = 0$  & $\Delta = 1$ & $\Delta = 2$ & & $\Delta = 0$ & $\Delta = 1$ & $\Delta = 2$ & & $\Delta = 0$ & $\Delta = 1$ & $\Delta = 2$ \\
\hline
1    & 1/4 & 0.051   & 0.020 & 0.007 & & 0.026 & 0.011 & 0.002 & & 0.005 & 0.001 & 0.000 \\
1    & 1/3 & 0.051	& 0.020	& 0.007	& & 0.026 &	0.011 &	0.002 & & 0.005	& 0.001	& 0.000 \\
0.01 & 1/4 & 0.096	& 0.058	& 0.038	& & 0.047 &	0.025 &	0.015 & & 0.009	& 0.005	& 0.001 \\
0.01 & 1/3 & 0.103	& 0.065	& 0.045	& & 0.049 &	0.030 &	0.017 & & 0.009	& 0.005	& 0.001 \\
\multicolumn{2}{c}{Theoretical}
           & 0.100  & 0.042 & 0.015 & & 0.050 & 0.018 & 0.006 & & 0.010 & 0.003 & 0.001 \\ \hline
\end{tabular}
\end{center}
\label{tab:pow}
\end{footnotesize}
\end{table}

We employ the developed framework to study whether the quantile treatment effect $\theta_0(\tau)$ is monotonically increasing in $\tau$, which corresponds to the special case of \eqref{conv1} in which $\Lambda$ equals the set of monotonically increasing functions. For ease of computation, we obtain quantile regression estimates $\hat \theta_n(\tau)$ on a grid $\{0.2, 0.225, \ldots, 0.775, 0.8\}$ and compute the distance of $\hat \theta_n$ to the set of monotone functions on this grid as our test statistic. In turn, critical values for this test statistic are obtained by computing two hundred bootstrapped quantile regression coefficients $\hat \theta_n^*(\tau)$ at all $\tau \in \{0.2, 0.225, \ldots, 0.775, 0.8\}$, and using the $1-\alpha$ quantile across bootstrap replications of the statistic $\hat \phi_n^\prime(\sqrt n\{\hat \theta_n^* - \hat \theta_n\})$, where $\hat \phi_n^\prime$ is computed according to \eqref{convth3} with $\epsilon_n = C n^\kappa$ for different choices of $C$ and $\kappa$. All reported results are based on five thousand Monte Carlo replications.

Table 1 reports the empirical rejection rates for different values of the local parameter $\Delta \in \{0, 1, 2\}$ -- recall that since $\theta_0(\tau) = \tau \Delta/\sqrt n$, the null hypothesis that $\theta_0$ is monotonically increasing is satisfied for all such $\Delta$. The bandwidth parameter $\epsilon_n$ employed in the construction of the estimator $\hat \phi_n^\prime$ is set according $\epsilon_n = Cn^\kappa$ for $C \in \{0.01, 1\}$ and $\kappa \in \{1/4, 1/3\}$. For the explored sample sizes of two and five hundred observations, we observe little sensitivity to the value of $\kappa$ but a more significant effect of the choice of $C$. In addition, the row labeled ``Theoretical" reports the rejection rates we should expect according to the local asymptotic approximation of Theorem \ref{th:localsize}. Throughout the specifications, we see that the test effectively controls size, and Theorem \ref{th:localsize} provides an adequate approximation often in between the rejection probabilities obtained from employing $C = 1$ and those corresponding to the more aggressive selection of $C = 0.01$.

In Table 2, we examine the local power of a 5$\%$ level test by considering values of $\Delta \in \{-1,\ldots, -8\}$. For such choices of the local parameter, the null hypothesis is violated since $\theta_0(\tau) = \tau \Delta /\sqrt n$ is in fact monotonically decreasing in $\tau$ (rather than increasing). In this context, we see that the theoretical local power is slightly above the empirical rejection rates, in particular for small values of $\Delta$. These distortions are most severe for $n$ equal to two hundred, though we note a quick improvement in the approximation error when $n$ is set to equal five hundred. Overall, we find the results of the Monte Carlo study encouraging, though certainly limited in their scope.

\begin{table}[t]
\begin{footnotesize}
\caption{Local Power of 0.05 Level Test}
\begin{center}
\begin{tabular}{cc cccccccc}
\hline
\hline
\multicolumn{2}{c}{Bandwidth}  &   \multicolumn{8}{c}{$n = 200$} \\ \cline{3-10}
$C$ & $\kappa$          & $\Delta = -1$ & $\Delta = -2$ & $\Delta = -3$ & $\Delta = -4$ & $\Delta = -5$ & $\Delta = -6$ & $\Delta = -7$ & $\Delta = -8$ \\ \hline
1       & 1/4           & 0.061        & 0.155        & 0.321        & 0.555        & 0.782        & 0.934        & 0.989        & 1.000        \\
1       & 1/3           & 0.061        & 0.155        & 0.321        & 0.555        & 0.782        & 0.934        & 0.989        & 1.000        \\
0.01    & 1/4           & 0.078        & 0.172        & 0.330        & 0.558        & 0.783        & 0.934        & 0.989        & 1.000        \\
0.01    & 1/3           & 0.081        & 0.174        & 0.331        & 0.559        & 0.783        & 0.934        & 0.989        & 1.000        \\
\multicolumn{2}{c}{Theoretical}
                        & 0.120        & 0.245        & 0.423        & 0.623        & 0.796        & 0.911        & 0.970        & 0.992 \\ \\

\multicolumn{2}{c}{Bandwidth}  &   \multicolumn{8}{c}{$n = 500$} \\ \cline{3-10}
$C$ & $\kappa$          & $\Delta = -1$ & $\Delta = -2$ & $\Delta = -3$ & $\Delta = -4$ & $\Delta = -5$ & $\Delta = -6$ & $\Delta = -7$ & $\Delta = -8$ \\ \hline
1       & 1/4           & 0.071        & 0.181        & 0.355        & 0.576        & 0.789        & 0.925        & 0.981        & 0.997        \\
1       & 1/3           & 0.071        & 0.181        & 0.355        & 0.576        & 0.789        & 0.925        & 0.981        & 0.997        \\
0.01    & 1/4           & 0.094        & 0.201        & 0.370        & 0.583        & 0.791        & 0.925        & 0.981        & 0.997        \\
0.01    & 1/3           & 0.098        & 0.204        & 0.371        & 0.585        & 0.791        & 0.925        & 0.981        & 0.997        \\
\multicolumn{2}{c}{Theoretical}
                        & 0.120        & 0.245        & 0.423        & 0.623        & 0.796        & 0.911        & 0.970        & 0.992 \\ \hline
\end{tabular}
\end{center}
\label{tab:size}
\end{footnotesize}
\end{table}

\section{Conclusion}\label{Sec: Conclusions}

In this paper, we have developed a general asymptotic framework for conducting inference in an important class of irregular models. In analogy with the Delta method, we have shown crucial features of these problems can be understood simply in terms of the asymptotic distribution $\mathbb G_0$ and the directional derivative $\phi_{\theta_0}^\prime$. The utility of these insights were demonstrated by both unifying diverse existing results and easily studying the otherwise challenging problem of testing for convex set membership. Further validation of the potential usefulness of our results has also been provided by a number of manuscripts that have exploited our results; see \cite{hansen2015regression}, \cite{jha2015testing}, and \cite{lee:bhattacharya} among others. We hope these are just the first applications of this framework, which should be of use to theorists and empirical researchers alike in determining statistical properties such as asymptotic distributions, bootstrap validity, and ability of tests to locally control size.

\newpage

\appendix


\begin{center}
{\Large {\sc \hypertarget{mainresults}{Appendix A}} - Proof of Main Results}
\end{center}

\renewcommand{\theequation}{A.\arabic{equation}}
\renewcommand{\thelemma}{A.\arabic{lemma}}
\renewcommand{\thecorollary}{A.\arabic{corollary}}
\renewcommand{\thetheorem}{A.\arabic{theorem}}
\setcounter{lemma}{0}
\setcounter{theorem}{0}
\setcounter{corollary}{0}
\setcounter{equation}{0}
\setcounter{remark}{0}

The following list includes notation and definitions that will be used in the appendix.

\begin{table}[h]
\begin{center}
\begin{tabular}{cl}
$ a \lesssim b$                     & $a\leq Mb$ for some constant $M$ that is universal in the proof.\\
$\|\cdot\|_{L^q(W)}$                & For a random variable $W$ and function $f$, $\|f\|_{L^q(W)} \equiv \{E[|f(W)|^q]\}^{\frac{1}{q}}$.\\
$\mathcal C(\mathbf A)$             & For a set $\mathbf A$, $\mathcal C(\mathbf A) \equiv \{f:\mathbf A \rightarrow \mathbf R : \sup_{a\in A} |f(a)| < \infty \text{ and } f \text{ is continuous}\}$. \\
$\vec{d}_H(\cdot,\cdot,\|\cdot\|)$  & For sets $A,B$, $\vec{d}_H(A,B,\|\cdot\|) \equiv \sup_{a \in A}\inf_{b\in B} \|a-b\|$.\\
$d_H(\cdot,\cdot,\|\cdot\|)$        & For sets $A,B$, $d_H(A,B,\|\cdot\|) \equiv \max\{\vec d_H(A,B,\|\cdot\|), \vec d_H(B,A,\|\cdot\|)\}$.\\
$\ell^\infty(\mathbf A)$            & For a set $\mathbf A$, $\ell^\infty(\mathbf A) \equiv \{f:\mathbf A \rightarrow \mathbf R : \sup_{a\in A} |f(a)| < \infty\}$.
\end{tabular}
\end{center}
\end{table}

\noindent {\sc Proof of Proposition \ref{pro:equivalence}:}  One direction is clear since, by definition, $\phi$ being Hadamard differentiable implies that its Hadamard directional derivative exists, equals the Hadamard derivative of $\phi$, and hence must be linear.

Conversely suppose the Hadamard directional derivative $\phi_\theta':\mathbb D_0 \rightarrow \mathbb E$ exists and is linear. Let $\{h_n\}$ and $\{t_n\}$ be sequences such that $h_n\to h\in\mathbb D_0$, $t_n\to 0$ and $\theta +t_nh_n\in\mathbb D_\phi$ for all $n$. Then note that from any subsequence $\{t_{n_k}\}$ we can extract a further subsequence $\{t_{n_{k_j}}\}$, such that either: (i) $t_{n_{k_j}} > 0$ for all $j$ or (ii) $t_{n_{k_j}} < 0$ for all $j$. When (i) holds, $\phi$ being Hadamard directional differentiable, then immediately yields that:
\begin{equation}\label{pro:equivalence1}
\lim_{j\rightarrow \infty} \frac{\phi(\theta+t_{n_{k_j}} h_{n_{k_j}})-\phi(\theta)}{t_{n_{k_j}}} = \phi_\theta'(h) ~.
\end{equation}
On the other hand, if (ii) holds, then $h \in \mathbb D_0$ and $\mathbb D_0$ being a subspace implies $- h \in \mathbb D_0$. Therefore, by Hadamard directional differentiability of $\phi$ and $-t_{n_{k_j}} > 0$ for all $j$:
\begin{multline}\label{pro:equivalence2}
\lim_{j\rightarrow \infty} \frac{\phi(\theta+t_{n_{k_j}} h_{n_{k_j}})-\phi(\theta)}{t_{n_{k_j}}} \\ = - \lim_{j\rightarrow \infty} \frac{\phi(\theta + (- t_{n_{k_j}})(- h_{n_{k_j}}))-\phi(\theta)}{-t_{n_{k_j}}} = -\phi_\theta'(- h)  = \phi_\theta^\prime(h) ~,
\end{multline}
where the final equality holds by the assumed linearity of $\phi_\theta^\prime$. Thus, results \eqref{pro:equivalence1} and \eqref{pro:equivalence2} imply that every subsequence $\{t_{n_k},h_{n_k}\}$ has a further subsequence along which
\begin{equation}\label{pro:equivalence3}
\lim_{j\rightarrow \infty} \frac{\phi(\theta+t_{n_{k_j}} h_{n_{k_j}})-\phi(\theta)}{t_{n_{k_j}}} = \phi_\theta'(h) ~.
\end{equation}
Since the subsequence $\{t_{n_k},h_{n_k}\}$ is arbitrary, it follows that \eqref{pro:equivalence3} must hold along the original sequence $\{t_n,h_n\}$ and hence $\phi$ is Hadamard differentiable tangentially to $\mathbb D_0$. \qed

\noindent {\sc Proof of Theorem \ref{th:delta}:} The proof closely follows the proof of Theorem 3.9.4 in \cite{Vaart1996}, and we include it here only for completeness. First, let $\mathbb D_n \equiv \{h \in \mathbb D : \theta_0 + h/r_n \in \mathbb D_\phi\}$ and define $g_n : \mathbb D_n \rightarrow \mathbb E$ to be given by
\begin{equation}\label{th:delta1}
g_n(h_n) \equiv r_n\{\phi(\theta_0 + \frac{h_n}{r_n}) - \phi(\theta_0)\}
\end{equation}
for any $h_n \in \mathbb D_n$. Then note that for every sequence $\{h_n\}$ with $h_n \in \mathbb D_n$ satisfying $\|h_n - h\|_{\mathbb D} = o(1)$ with $h\in \mathbb D_0$, it follows from Assumption \ref{ass:structure}(ii) that $\|g_n(h_n) - \phi_{\theta_0}^\prime(h)\|_{\mathbb E} = o(1)$. Therefore, the first claim follows by Theorem 1.11.1 in \cite{Vaart1996} and $\mathbb G_0$ being tight implying that it is also separable by Lemma 1.3.2 in \cite{Vaart1996}.

For the second claim of the Theorem, we define $f_n : \mathbb D_n \times \mathbb D \rightarrow \mathbb E \times \mathbb E$ by:
\begin{equation}\label{th:delta2}
f_n(h_n,h) = (g_n(h_n),\phi^\prime_{\theta_0}(h)) ~,
\end{equation}
for any $(h_n,h) \in \mathbb D_n \times \mathbb D$. It then follows by applying Theorem 1.11.1 in \cite{Vaart1996} again, that as processes in $\mathbb E\times \mathbb E$ we have:
\begin{equation}\label{th:delta3}
\left[\begin{array}{c} r_n\{\phi(\hat \theta_n) - \phi(\theta_0)\} \\ \phi_{\theta_0}^\prime(r_n\{\hat \theta_n - \theta_0\})\end{array}\right] \stackrel{L}{\rightarrow}
\left[\begin{array}{c} \phi_{\theta_0}^\prime(\mathbb G_0) \\ \phi_{\theta_0}^\prime(\mathbb G_0)\end{array}\right] ~.
\end{equation}
In particular, result \eqref{th:delta3} and the continuous mapping theorem allow us to conclude:
\begin{equation}\label{th:delta4}
r_n\{\phi(\hat \theta_n) - \phi(\theta_0)\} - \phi_{\theta_0}^\prime(r_n(\hat \theta_n - \theta_0)) \stackrel{L}{\rightarrow} 0 ~.
\end{equation}
The second claim then follows from \eqref{th:delta4} and Lemma 1.10.2(iii) in \cite{Vaart1996}. \qed

\noindent {\sc Proof of Theorem \ref{th:gauslin}:} Let $P$, $\mathbb D_L \subseteq \mathbb D$, and $\mu_0$ respectively denote the distribution, support, and mean of $\mathbb G_0$. Since $\mathbb D_L$ is a vector space, it follows that $0 \in \mathbb D_L$ and $\mathbb D_L = \mathbb D_L + \mathbb D_L$, and hence Theorem \ref{th:bootiff} implies that result \eqref{th:gauslindisp} holds if and only if
\begin{equation}\label{th:gauslin1}
E[f(\phi_{\theta_0}^\prime(\bar{\mathbb G}_0 + \mu_0))] = E[f(\phi_{\theta_0}^\prime(\bar{\mathbb G}_0 + \mu_0 + a_0) - \phi_{\theta_0}^\prime(a_0))]
\end{equation}
for $\bar{\mathbb G}_0 = \mathbb G_0 - \mu_0$, and all $a_0 \in \mathbb D_L$ and $f\in \text{BL}_1(\mathbb E)$. On the other hand, $\mathbb D_L$ is a subspace of $\mathbb D$, and therefore Proposition \ref{pro:equivalence} implies that $\phi$ is Hadamard differentiable at $\theta_0 \in \mathbb D_\phi$ tangentially to $\mathbb D_L$ if and only if $\phi_{\theta_0}^\prime : \mathbb D_L \rightarrow \mathbb E$ is linear. Thus, the claim of the Theorem will follow from establishing that \eqref{th:gauslin1} holds if and only if $\phi_{\theta_0}^\prime : \mathbb D_L \rightarrow \mathbb E$ is linear. To this end, we note that one direction is trivial, since linearity of $\phi_{\theta_0}^\prime$ implies
\begin{equation}\label{th:gauslin2}
\phi_{\theta_0}^\prime(\bar{\mathbb G}_0 + \mu_0 + a_0) - \phi_{\theta_0}^\prime(a_0) = \phi_{\theta_0}^\prime(\bar{\mathbb G}_0 +\mu_0)
\end{equation}
$P$ almost surely for all $a_0 \in \mathbb D_L$, and thus \eqref{th:gauslin1} must hold for any $f\in \text{BL}_1(\mathbb E)$.

The opposite direction is more challenging and requires us to introduce additional notation which closely follows Chapter 7 in \cite{davydov:lifshits:smorodina:1998}. First, we note that by Lemma \ref{aux:supp} $\mathbb D_L$ is a separable Banach space under $\|\cdot\|_{\mathbb D}$. Next, let $\mathbb D^*_L$ denote the dual space of $\mathbb D_L$, and $\langle d,d^*\rangle_{\mathbb D} = d^*(d)$ for any $d\in \mathbb D_L$ and $d^*\in \mathbb D^*_L$. Similarly denote the dual space of $\mathbb E$ by $\mathbb E^*$ and corresponding bilinear form by $\langle \cdot,\cdot\rangle_{\mathbb E}$. Further let
\begin{equation}\label{th:gauslin3}
\mathbb D^\prime_P \equiv \Big\{d^\prime : \mathbb D_L \rightarrow \mathbf R : d^\prime \text{ is linear, Borel-measurable, and } \int_{\mathbb D} (d^\prime(d))^2 dP(d) < \infty \Big\} ~,
\end{equation}
and with some abuse of notation also write $d^\prime(d) = \langle d^\prime, d\rangle_{\mathbb D}$ for any $d^\prime \in \mathbb D^\prime_P$ and $d\in \mathbb D_L$. In addition, we let $\bar P$ denote the distribution of $\bar {\mathbb G}_0$, and note that $\bar P$ is a centered Gaussian measure whose support also equals $\mathbb D_L$ since $\mu_0 \in \mathbb D_L$ by Lemma \ref{aux:supp}. Finally, for each $h\in \mathbb D_L$ we let $\bar P^h$ denote the law of $\bar{\mathbb G}_0 + h$, write $\bar P^h \ll \bar P$ whenever $\bar P^h$ is absolutely continuous with respect to $\bar P$, and define the set:
\begin{equation}\label{th:gauslin4}
\mathbb H_P \equiv \{h \in \mathbb D_L : \bar P^{rh} \ll \bar P \text{ for all } r\in \mathbf R\} ~.
\end{equation}
To proceed, note that since $\mathbb D_L$ is separable, the Borel $\sigma$-algebra, the $\sigma$-algebra generated by the weak topology, and the cylindrical $\sigma$-algebra all coincide \citep[p. 38]{ledoux:talagrand}. Furthermore, by Theorem 7.1.7 in \cite{bogachev2:2007}, $\bar P$ is Radon and thus by Theorem 7.1 in \cite{davydov:lifshits:smorodina:1998}, it follows that there exists a linear map $I:\mathbb H_P \rightarrow \mathbb D^\prime_P$ such that for every $h\in \mathbb H_P$ we have:
\begin{equation}\label{th:gauslin5}
\frac{d \bar P^h}{d\bar P}(d) = \exp\Big\{ \langle d,Ih\rangle_{\mathbb D} - \frac{1}{2} \sigma^2(h)\Big\} \hspace{0.5 in} \sigma^2(h) \equiv \int_{\mathbb D} \langle d,Ih \rangle_{\mathbb D}^2d \bar P(d) ~.
\end{equation}

Next, fix an arbitrary $e^* \in \mathbb E^*$ and $h \in \mathbb H_P$. Then note that if \eqref{th:gauslin1} holds, then Lemma 1.3.12 in \cite{Vaart1996} implies $\langle e^*, \phi_{\theta_0}^\prime(\bar{\mathbb G}_0 + \mu_0 + rh) - \phi_{\theta_0}^\prime(rh)\rangle_{\mathbb E}$ and $\langle e^*, \phi_{\theta_0}^\prime(\bar{\mathbb G}_0 + \mu_0)\rangle_{\mathbb E}$ must be equal in distribution for all $r\in \mathbf R$. Thus, their characteristic functions must be equal, and hence for all $r\geq 0$ and $t\in \mathbf R$:
\begin{multline}\label{th:gauslin6}
E[\exp\{it\langle e^*,\phi_{\theta_0}^\prime(\bar{\mathbb G}_0 + \mu_0)\rangle_{\mathbb E}\}] = E[\exp\{it\{\langle e^*, \phi_{\theta_0}^\prime(\bar{\mathbb G}_0 + \mu_0 + rh) - \phi_{\theta_0}^\prime(rh)\rangle_{\mathbb E}\}\}] \\
 = \exp\{-itr\langle e^*,\phi_{\theta_0}^\prime(h)\rangle_{\mathbb E}\} E[\exp\{it\langle e^*, \phi_{\theta_0}^\prime(\bar{\mathbb G}_0 + \mu_0 + rh)\rangle_{\mathbb E}\}]~,
\end{multline}
where in the second equality we have exploited that $\phi_{\theta_0}^\prime(rh) = r\phi_{\theta_0}^\prime(h)$ due to $\phi_{\theta_0}^\prime$ being positively homogenous of degree one. Setting $C(t) \equiv  E[\exp\{it\langle e^*,\phi_{\theta_0}^\prime(\bar{\mathbb G}_0+\mu_0)\rangle_{\mathbb E}\}]$ and exploiting result \eqref{th:gauslin6} we can then obtain by direct calculation that for all $t\in \mathbf R$
\begin{multline}\label{th:gauslin7}
itC(t)\times \langle e^*,\phi_{\theta_0}^\prime(h)\rangle_{\mathbb E} = \lim_{r\downarrow 0} \frac{1}{r}\{E[\exp\{it\langle e^*, \phi_{\theta_0}^\prime(\bar{\mathbb G}_0 + \mu_0+ rh)\rangle_{\mathbb E}\}] - C(t)\} \\
= \lim_{r\downarrow 0} \frac{1}{r}\int_{\mathbb D} \Big\{\exp\Big\{it\langle e^*,\phi_{\theta_0}^\prime(d+\mu_0)\rangle_{\mathbb E} + r\langle d,Ih\rangle_{\mathbb D} - \frac{r^2}{2}\sigma^2(h)\Big\} -C(t)\Big\}d\bar P(d)
\end{multline}
where in the second equality we exploited result \eqref{th:gauslin5}, linearity of $I:\mathbb H_P\rightarrow \mathbb D^\prime_P$ and that $h\in \mathbb H_P$ implies $rh \in \mathbb H_P$ for all $r\in \mathbf R$. Furthermore, by the mean value theorem
\begin{align}\label{th:gauslin8}
\sup_{r\in (0,1]} & \frac{1}{r} \Big|\exp\Big\{it\langle e^*,\phi_{\theta_0}^\prime(d + \mu_0)\rangle_{\mathbb E} + r\langle d,Ih\rangle_{\mathbb D} - \frac{r^2}{2}\sigma^2(h)\Big\} -\exp\{it\langle e^*,\phi_{\theta_0}^\prime(d+\mu_0)\rangle_{\mathbb E}\}\Big| \nonumber \\
& \leq \sup_{r\in (0,1]} \Big|\exp\Big\{it\langle e^*,\phi_{\theta_0}^\prime(d+\mu_0)\rangle_{\mathbb E}+ r\langle d,Ih\rangle_{\mathbb D} - \frac{r^2}{2}\sigma^2(h)\Big\}\times \{\langle d,Ih\rangle_{\mathbb D} - r\sigma^2(h)\}\Big| \nonumber \\ & \leq \exp\{|\langle d, Ih\rangle_{\mathbb D}|\}\times \{|\langle d,Ih\rangle_{\mathbb D}| + \sigma^2(h)\} ~,
\end{align}
where the final inequality follows from $\sigma^2(h) > 0$ and $|\exp\{it\langle e^*,\phi_{\theta_0}^\prime(d+\mu_0)\rangle_{\mathbb E}\}|\leq 1$. Moreover, by Proposition 2.10.3 in \cite{bogachevgauss} and $Ih \in \mathbb D^\prime_P$, it follows that $\langle \bar{\mathbb G}_0, Ih\rangle_{\mathbb D} \sim N(0,\sigma^2(h))$. Thus, we can obtain by direct calculation:
\begin{multline}\label{th:gauslin9}
\int_{\mathbb D} \exp\{|\langle d, Ih\rangle_{\mathbb D}|\}\times \{|\langle d,Ih\rangle_{\mathbb D}| + \sigma^2(h)\} d\bar P(d) \\ = \int_{\mathbf R}\frac{\{|u| + \sigma^2(h)\}}{\sigma(h)\sqrt{2\pi}}\times \exp\Big\{|u| -\frac{u^2}{2\sigma^2(h)}\Big\}du < \infty ~.
\end{multline}
Hence, results \eqref{th:gauslin8} and \eqref{th:gauslin9} justify the use of the dominated convergence theorem in \eqref{th:gauslin7}. Also note that $t\mapsto C(t)$ is the characteristic function of $\langle e^*,\phi_{\theta_0}^\prime(\bar{\mathbb G}_0+\mu_0)\rangle_{\mathbb E}$ and hence it is continuous. Thus, since $C(0) = 1$ there exists a $t_0 > 0$ such that $C(t_0)t_0 \neq 0$. For such $t_0$ we then finally obtain from the above results that
\begin{equation}\label{th:gauslin10}
\langle e^*,\phi_{\theta_0}^\prime(h)\rangle_{\mathbb E} = -\frac{iE[\exp\{it_0\langle e^*,\phi_{\theta_0}^\prime(\bar{\mathbb G}_0+\mu_0)\rangle_{\mathbb E}\}\langle \bar{\mathbb G}_0,Ih\rangle_{\mathbb D}]}{t_0C(t_0)} ~.
\end{equation}

To conclude note that $\mathbb H_P$ being a vector space \citep[p. 38]{davydov:lifshits:smorodina:1998} and $I:\mathbb H_P \rightarrow \mathbb D^\prime_P$ being linear imply together with result \eqref{th:gauslin9} that $h \mapsto \langle e^*,\phi_{\theta_0}^\prime(h)\rangle_{\mathbb E}$ is linear on $\mathbb H_P$. Moreover, note that $h \mapsto \langle e^*,\phi_{\theta_0}^\prime(h)\rangle_{\mathbb E}$ is also continuous on $\mathbb D_L$ due to continuity of $\phi_{\theta_0}^\prime$ and having $e^* \in \mathbb E^*$. Hence, since $\mathbb H_P$ is dense in $\mathbb D_L$ by Proposition 7.4(ii) in \cite{davydov:lifshits:smorodina:1998} we can conclude that $\langle e^*,\phi_{\theta_0}^\prime(\cdot)\rangle_{\mathbb E} : \mathbb D_L \rightarrow \mathbf R$ is linear and continuous. Since this result holds for all $e^* \in \mathbb E^*$, Lemma A.2 in \cite{Vaart1991differentibility} implies $\phi_{\theta_0}^\prime : \mathbb D_L \rightarrow \mathbb E$ must be linear and continuous, which establishes the Theorem. \qed

\noindent {\sc Proof of Corollary \ref{cor:nonnorm}:} By Theorem \ref{th:gauslin} and Proposition \ref{pro:equivalence} the bootstrap is consistent if and only if $\phi_{\theta_0}^\prime$ is linear. However, since $\mathbb G_0$ is Gaussian and $\phi_{\theta_0}^\prime : \mathbb D_0 \rightarrow \mathbb E$ is continuous, Lemma 2.2.2 in \cite{bogachevgauss} implies $\phi_{\theta_0}^\prime(\mathbb G_0)$ must be Gaussian (on $\mathbb E$) whenever $\phi_{\theta_0}^\prime$ is linear, and hence the claim of the Corollary follows. \qed

\noindent {\sc Proof of Theorem \ref{th:modboot}:} Fix arbitrary $\epsilon > 0$, $\eta > 0$ and for notational convenience let $\mathbb G_n^* \equiv r_n\{\hat \theta_n^* - \hat \theta_n\}$. By Assumption \ref{ass:paramest}(ii) there is a compact set $K_0 \subseteq \mathbb D_0$ such that
\begin{equation}\label{th:modboot0p5}
P(\mathbb G_0 \notin K_0) < \frac{\epsilon \eta}{2} ~.
\end{equation}
Thus, by Lemma \ref{lm:auxGnboot} and the Portmanteau Theorem, we conclude that for any $\delta > 0$
\begin{equation}\label{th:modboot1}
\limsup_{n\rightarrow \infty} P(\mathbb G_n^* \notin K_0^\delta) \leq P(\mathbb G_0 \notin K_0^\delta) \leq  P(\mathbb G_0 \notin K_0) < \frac{\epsilon \eta}{2} ~.
\end{equation}
On the other hand, since $K_0$ is compact, Assumption \ref{ass:derest} yields that for some $\delta_0 > 0$:
\begin{equation}\label{th:modboot2}
\limsup_{n\rightarrow \infty} P(\sup_{h\in K_0^{\delta_0}} \|\hat \phi_n^\prime(h) - \phi_{\theta_0}^\prime(h)\|_{\mathbb E} > \epsilon) < \eta ~.
\end{equation}
Next, note that Lemma 1.2.2(iii) in \cite{Vaart1996}, $h \in \text{BL}_1(\mathbb E)$ being bounded by one and satisfying $|h(e_1) - h(e_2)| \leq \|e_1 - e_2\|_{\mathbb E}$ for all $e_1,e_2 \in \mathbb E$, imply:
\begin{align}\label{th:modboot3}
\sup_{f \in \text{BL}_1(\mathbb E)} |E[f(\hat \phi_n^\prime(\mathbb G_n^*))& |\{X_i\}]  -  E[f(\phi_{\theta_0}^\prime(\mathbb G_n^*))|\{X_i\}]| \nonumber \\
& \leq \sup_{f \in \text{BL}_1(\mathbb E)} E[|f(\hat \phi_n^\prime(\mathbb G_n^*)) - f(\phi_{\theta_0}^\prime(\mathbb G_n^*))||\{X_i\}] \nonumber \\
& \leq E[2\times 1\{\mathbb G_n^* \notin K_0^{\delta_0}\} + \sup_{f \in K_0^{\delta_0}} \|\hat \phi^\prime_n(f) - \phi^\prime_{\theta_0}(f)\|_{\mathbb E}|\{X_i\}] \nonumber \\
& \leq 2P(\mathbb G_n^* \notin K_0^{\delta_0}|\{X_i\}_{i=1}^n) + \sup_{f\in K_0^{\delta_0}} \|\hat \phi_n^\prime(f) - \phi_{\theta_0}^\prime(f)\|_{\mathbb E} ~,
\end{align}
where in the final inequality we exploited Lemma 1.2.2(i) in \cite{Vaart1996} and $\hat \phi_n^\prime: \mathbb D\rightarrow \mathbb E$ depending only on $\{X_i\}_{i=1}^n$. Furthermore, Markov's inequality, Lemma 1.2.7 in \cite{Vaart1996}, and result \eqref{th:modboot1} yield:
\begin{equation}\label{th:modboot4}
\limsup_{n\rightarrow \infty} P( P(\mathbb G_n^* \notin K_0^{\delta_0}|\{X_i\}_{i=1}^n) > \epsilon) \leq \limsup_{n\rightarrow \infty} \frac{1}{\epsilon} P( \mathbb G_n^* \notin K_0^{\delta_0}) < \eta ~.
\end{equation}
Next, also note that Assumption \ref{ass:primboot}(i) and Theorem 10.8 in \cite{Kosorok2008} imply that:
\begin{equation}\label{th:modboot5}
\sup_{f\in \text{BL}_1(\mathbb E)}|E[f(\phi_{\theta_0}^\prime(\mathbb G_n^*))|\{X_i\}_{i=1}^n] - E[f(\phi_{\theta_0}^\prime(\mathbb G_0))]| = o_p(1) ~.
\end{equation}
Thus, by combining results \eqref{th:modboot2}, \eqref{th:modboot3}, \eqref{th:modboot4} and \eqref{th:modboot5} we can finally conclude:
\begin{equation}\label{th:modboot6}
\limsup_{n\rightarrow \infty} P( \sup_{f\in \text{BL}_1(\mathbb E)} |E[f(\hat \phi_n^\prime(\mathbb G_n^*))|\{X_i\}_{i=1}^n] - E[f(\phi_{\theta_0}^\prime(\mathbb G_0))]| > 3\epsilon) < 3\eta ~.
\end{equation}
Since $\epsilon$ and $\eta$ were arbitrary, the claim of the Theorem then follows from \eqref{th:modboot6}. \qed

\noindent {\sc Proof of Corollary \ref{cor:critval}:} Let $F$ denote the cdf of $\phi^\prime_{\theta_0}(\mathbb G_0)$, and similarly define:
\begin{equation}\label{cor:critval1}
\hat F_n(c) \equiv P(\hat \phi_n^\prime(r_n\{\hat \theta_n^* - \hat \theta_n\}) \leq c |\{X_i\}_{i=1}^n) ~.
\end{equation}
Next, observe that Theorem \ref{th:modboot} and Lemma 10.11 in \cite{Kosorok2008} imply that:
\begin{equation}\label{cor:critval2}
\hat F_n(c) = F(c) + o_p(1) ~,
\end{equation}
for all $c\in \mathbf R$ that are continuity points of $F$. Fix $\epsilon > 0$, and note that since $F$ is strictly increasing at $c_{1-\alpha}$ and the set of continuity of points of $F$ is dense in $\mathbf R$, it follows that there exist points $c_1,c_2\in \mathbf R$ such that: (i) $c_1 < c_{1-\alpha} < c_2$, (ii) $|c_1 - c_{1-\alpha}| < \epsilon$ and $|c_2 - c_{1-\alpha}| < \epsilon$, (iii) $c_1$ and $c_2$ are continuity points of $F$, and (iv) $F(c_1) + \delta < 1-\alpha < F(c_2) - \delta$ for some $\delta > 0$. We can then conclude that:
\begin{multline}\label{cor:critval3}
\limsup_{n\rightarrow \infty} P( |\hat c_{1-\alpha} - c_{1-\alpha}| > \epsilon) \\ \leq \limsup_{n\rightarrow \infty} \{P( |\hat F_n(c_1) - F(c_1)| > \delta) + P( |\hat F_n(c_2) - F(c_2)| > \delta)\} = 0 ~,
\end{multline}
due to \eqref{cor:critval2}. Since $\epsilon > 0$ was arbitrary, the Corollary then follows. \qed

\noindent {\sc Proof of Lemma \ref{th:localdist}:} First note that by Assumption \ref{ass:paramreg}(i) we can conclude:
\begin{equation}\label{th:localdist1}
\lim_{n\rightarrow \infty} \| r_n\{\theta(P_{n,\lambda}^\infty) - \theta_0\} - \theta^\prime(\lambda)\|_{\mathbb D} = 0 ~.
\end{equation}
Hence, letting $t_n \equiv r_n^{-1}$, $h_n \equiv r_n\{\theta(P_{n,\lambda}^\infty) - \theta_0\}$ we note $\theta_0 + t_nh_n = \theta(P_{n,\lambda}^\infty) \in \mathbb D_\phi$, and by \eqref{th:localdist1} that $\|h_n - h\|_{\mathbb D} = o(1)$ for $h \equiv \theta^\prime(\lambda)$. Therefore, Assumption \ref{ass:structure}(ii) yields
\begin{multline}\label{th:localdist2}
\lim_{n\rightarrow \infty} \|r_n\{\phi(\theta(P_{n,\lambda}^\infty)) - \phi(\theta_0)\} - \phi^\prime_{\theta_0}(\theta^\prime(\lambda))\|_{\mathbb E} \\ = \lim_{n\rightarrow \infty} \| \frac{\phi(\theta_0 + t_n h_n) - \phi(\theta_0)}{t_n} - \phi^\prime_{\theta_0}(h)\|_{\mathbb E} = 0 ~.
\end{multline}

Next, note that Theorem \ref{th:delta} and Assumption \ref{ass:localset}(i) implies that we have
\begin{equation}\label{th:localdist3}
r_n \{\phi(\hat \theta_n) - \phi(\theta_0)\} = \phi_{\theta_0}^\prime(r_n\{\hat \theta_n - \theta_0\}) + o_p(1)
\end{equation}
under $P_{n,0}^n$. However, by Assumption \ref{ass:localset}(ii) and Theorem 18.9 in \cite{strasser1985mathematical}, the sequence $P_{n,\lambda}^n$ is contiguous to $P_{n,0}^n$. Hence, from \eqref{th:localdist2} and \eqref{th:localdist3} we obtain
\begin{align}\label{th:localdist4}
r_n \{\phi(\hat \theta_n) - \phi(\theta(P_{n,\lambda}^n)\} & = r_n \{\phi(\hat \theta_n) - \phi(\theta_0)\} - r_n \{\phi(\theta(P_{n,\lambda}^\infty)) - \phi(\theta_0)\} \nonumber \\
& = \phi_{\theta_0}^\prime(r_n\{\hat \theta_n - \theta_0\}) - \phi_{\theta_0}^\prime(\theta^\prime(\lambda)) +  o_p(1)
\end{align}
under $P_{n,\lambda}^n$. Furthermore, by Assumption \ref{ass:paramreg}(ii) and result \eqref{th:localdist1} we also have that:
\begin{equation}\label{th:localdist5}
r_n\{\hat \theta_n - \theta_0\} = r_n\{\hat \theta_n - \theta(P_{n,\lambda}^\infty)\} + r_n\{\theta(P_{n,\lambda}^\infty) - \theta_0\} \stackrel{L_\lambda}{\rightarrow} \mathbb G_0 + \theta^\prime(\lambda) ~.
\end{equation}
Finally, note that since \eqref{th:localdist5} and $P_{n,\lambda}^n$ being contiguous to $P_{n,0}^n$ imply that the support of $\mathbb G_0 + \theta^\prime(\lambda)$ is a (weak) subset of the support of $\mathbb G_0$, Assumption \ref{ass:paramest}(ii) yields that
\begin{equation}\label{th:localdist6}
P(\mathbb G_0 + \theta^\prime(\lambda) \in \mathbb D_0) = 1 ~.
\end{equation}
Thus, the Lemma follows from \eqref{th:localdist4}-\eqref{th:localdist6} and the continuous mapping theorem. \qed

\noindent {\sc Proof of Theorem \ref{th:localsize}:} First note that by Assumption \ref{ass:localset}(ii) and Theorem 18.9 in \cite{strasser1985mathematical}, the sequence $P^n_{n,\lambda}$ is contiguous to $P_{n,0}^n$. Therefore, since Assumption \ref{ass:localset}(i) and Corollary \ref{cor:critval} imply $\hat c_{1-\alpha} \stackrel{p}{\rightarrow} c_{1-\alpha}$ under $P_{n,0}^n$, it follows that
\begin{equation}\label{th:localsize1}
\hat c_{1-\alpha} = c_{1-\alpha} + o_p(1) ~ \text{ under } P_{n,\lambda}^n ~.
\end{equation}
Moreover, since $\phi(\theta_0) = 0$ by Assumption \ref{ass:4test}(i), result \eqref{th:localdist3} implies that under $P_{n,\lambda}^n$:
\begin{equation}\label{th:localsize2}
r_n\phi(\hat \theta_n) = r_n \{\phi(\hat \theta_n) - \phi(\theta_0)\} = \phi_{\theta_0}^\prime(r_n \{\hat \theta_n - \theta_0\}) + o_p(1) \stackrel{L_\lambda}{\rightarrow} \phi_{\theta_0}^\prime(\mathbb G_0 + \theta^\prime(\lambda)) ~,
\end{equation}
where the final result holds for $L_\lambda$ denoting law under $P_{n,\lambda}^n$ by \eqref{th:localdist5}, \eqref{th:localdist6} and the continuous mapping theorem. Thus, the Portmanteau Theorem and \eqref{th:localsizedisp1} imply \eqref{th:localsize2}.

In order to establish \eqref{th:localsizedisp2} holds whenever $\phi(\theta(P_{n,\lambda}^\infty)) \leq 0$ for all $n$, note \eqref{th:localdist2} implies
\begin{equation}\label{th:localsize3}
0 \geq \lim_{n\rightarrow \infty} r_n\{\phi(\theta(P_{n,\lambda}^\infty)) - \phi(\theta_0)\} = \phi_{\theta_0}^\prime(\theta^\prime(\lambda)) ~,
\end{equation}
where we have exploited that $\phi(\theta_0) = 0$ and $\phi(\theta(P_{n,\lambda}^\infty)) \leq 0$ for all $n$. Therefore, result \eqref{th:localsize1} together with the second equality in \eqref{th:localsize2} allow us to conclude
\begin{align}
\limsup_{n\rightarrow \infty}  P_{n,\lambda}^n & (r_n\phi(\hat \theta_n) > \hat c_{1-\alpha}) \nonumber \\ & \leq \limsup_{n\rightarrow \infty} P_{n,\lambda}^n(\phi_{\theta_0}^\prime(r_n\{\hat \theta_n - \theta_0\})  \geq c_{1-\alpha}) \nonumber \\ & \leq \limsup_{n\rightarrow \infty} P_{n,\lambda}^n(\phi_{\theta_0}^\prime(r_n\{\hat \theta_n - \theta(P_{n,\lambda}^\infty)\}) + \phi_{\theta_0}^\prime(r_n\{\theta(P_{n,\lambda}^\infty) - \theta_0\}) \geq c_{1-\alpha}) \nonumber \\
& \leq \limsup_{n\rightarrow \infty} P_{n,\lambda}^n(\phi_{\theta_0}^\prime(r_n\{\hat \theta_n - \theta(P_{n,\lambda}^\infty)\}) \geq c_{1-\alpha}) \nonumber \\
& = P(\phi_{\theta_0}^\prime(\mathbb G_0) \geq c_{1-\alpha}) ~,
\end{align}
where the second inequality follows from subadditivity of $\phi_{\theta_0}^\prime$, the third inequality is implied by \eqref{th:localsize3}, and the final result follows from $r_n\{\hat \theta_n - \theta(P_{n,\lambda}^\infty)\}\stackrel{L_\lambda}{\rightarrow} \mathbb G_0$ by Assumption \ref{ass:paramreg}(ii), the continuous mapping theorem, and $c_{1-\alpha}$ being a continuity point of the cdf of $\phi_{\theta_0}^\prime(\mathbb G_0)$. Since $P(\phi_{\theta_0}^\prime(\mathbb G_0) \geq c_{1-\alpha}) = \alpha$ by construction, result \eqref{th:localsizedisp2} follows. \qed


\begin{theorem}\label{th:bootiff}
Let Assumptions \ref{ass:structure}, \ref{ass:paramest}, \ref{ass:primboot}, and \ref{ass:extraboot} hold, $\mathbb D_L$ denote the support of $\mathbb G_0$, $0\in \mathbb D_L$, and $\mathbb D_0 = \mathbb D_0+\mathbb D_0$. Then, the following statements are equivalent:

\vspace{-0.2 in}
\begin{packed_enum}
    \item[(i)] $E[f(\phi_{\theta_0}^\prime(\mathbb G_0))] = E[f(\phi_{\theta_0}^\prime(\mathbb G_0 + a_0) - \phi_{\theta_0}^\prime(a_0))] \text{ for all } a_0\in \mathbb D_L \text{ and } f\in \text{BL}_1(\mathbb E)$.
    \item[(ii)] $\sup_{f \in \text{BL}_1(\mathbb E)} |E[f(r_n\{\phi(\hat \theta_n^*) - \phi(\hat \theta_n)\})|\{X_i\}_{i=1}^n] - E[f(\phi_{\theta_0}^\prime(\mathbb G_0))]| = o_p(1)$.
\end{packed_enum}
\end{theorem}

\noindent {\sc Proof:} In these arguments we need to distinguish between outer and inner expectations, and we therefore employ the notation $E^*$ and $E_*$ respectively. In addition, for notational convenience we let $\mathbb G_n \equiv r_n\{\hat \theta_n - \theta_0\}$ and $\mathbb G_n^* \equiv r_n\{\hat \theta_n^* - \hat \theta_n\}$. To begin, note that Lemma \ref{lm:auxGnindep} and the continuous mapping theorem imply that:
\begin{multline}\label{th:bootiff1}
(r_n\{\hat \theta_n^* - \theta_0\}, r_n\{\hat \theta_n - \theta_0\}) \\ = (r_n\{\hat \theta_n^* - \hat \theta_n\} + r_n\{\hat \theta_n - \theta_0\},r_n\{\hat \theta_n - \theta_0\}) \stackrel{L}{\rightarrow} (\mathbb G_1 + \mathbb G_2, \mathbb G_2)
\end{multline}
on $\mathbb D\times \mathbb D$, where $\mathbb G_1$ and $\mathbb G_2$ are independent copies of $\mathbb G_0$. Further let $\Phi : \mathbb D_\phi \times \mathbb D_\phi \rightarrow \mathbb E$ be given by $\Phi(\theta_1,\theta_2) = \phi(\theta_1) - \phi(\theta_2)$ for any $\theta_1,\theta_2 \in \mathbb D_\phi\times \mathbb D_\phi$. Then observe that Assumption \ref{ass:structure}(ii) implies $\Phi$ is Hadamard directionally differentiable at $(\theta_0,\theta_0)$ tangentially to $\mathbb D_0\times \mathbb D_0$ with derivative $\Phi^\prime_{\theta_0} : \mathbb D_0\times \mathbb D_0 \rightarrow \mathbb E$ given by
\begin{equation}\label{th:bootdiff2}
\Phi^\prime_{\theta_0}(h_1,h_2) =  \phi_{\theta_0}^\prime(h_1)- \phi^\prime_{\theta_0}(h_2)
\end{equation}
for any $(h_1,h_2) \in \mathbb D_0 \times \mathbb D_0$. Thus, by Assumptions \ref{ass:paramest}(ii) and $\mathbb D_0 = \mathbb D_0 + \mathbb D_0$, Theorem \ref{th:delta}, result \eqref{th:bootiff1}, and $r_n\{\hat \theta_n^* - \theta_0\} = \mathbb G_n^* + \mathbb G_n$ we can conclude that
\begin{multline}\label{th:bootiff3}
r_n\{\phi(\hat \theta_n^*) - \phi(\hat \theta_n)\} = r_n\{\Phi(\hat \theta_n^*,\hat \theta_n) - \Phi(\theta_0,\theta_0)\} \\ = \Phi^\prime_{\theta_0}(\mathbb G_n^* + \mathbb G_n,\mathbb G_n) + o_p(1)  = \phi_{\theta_0}^\prime(\mathbb G_n^* + \mathbb G_n) - \phi_{\theta_0}^\prime(\mathbb G_n) + o_p(1)~.
\end{multline}
Further observe that for any $\epsilon > 0$, it follows from the definition of $\text{BL}_1(\mathbb E)$ that:
\begin{multline}\label{th:bootdiff4}
\sup_{h \in \text{BL}_1(\mathbb E)} | E^*[h(r_n\{\phi(\hat \theta_n^*) - \phi(\hat \theta_n)\}) - h(\phi_{\theta_0}^\prime(\mathbb G_n^* + \mathbb G_n) - \phi_{\theta_0}^\prime(\mathbb G_n))|\{X_i\}_{i=1}^n] |\\ \leq \epsilon + 2P^*(\|r_n\{\phi(\hat \theta_n^*) - \phi(\hat \theta_n)\} - \{\phi_{\theta_0}^\prime(\mathbb G_n^* + \mathbb G_n) - \phi_{\theta_0}^\prime(\mathbb G_n)\}\|_{\mathbb E} > \epsilon|\{X_i\}_{i=1}^n)
\end{multline}
Moreover, Lemma 1.2.6 in \cite{Vaart1996} and result \eqref{th:bootiff3} also yield:
\begin{multline}\label{th:bootdiff5}
E^*[P^*(\|r_n\{\phi(\hat \theta_n^*) - \phi(\hat \theta_n)\} - \{\phi_{\theta_0}^\prime(\mathbb G_n^* + \mathbb G_n) - \phi_{\theta_0}^\prime(\mathbb G_n)\}\|_{\mathbb E} > \epsilon|\{X_i\}_{i=1}^n)] \\ \leq P^*(\|r_n\{\phi(\hat \theta_n^*) - \phi(\hat \theta_n)\} - \{\phi_{\theta_0}^\prime(\mathbb G_n^* + \mathbb G_n) - \phi_{\theta_0}^\prime(\mathbb G_n)\}\|_{\mathbb E} > \epsilon) = o(1) ~.
\end{multline}
Therefore, since $\epsilon > 0$ was arbitrary, we obtain from results \eqref{th:bootdiff4} and \eqref{th:bootdiff5} that:
\begin{multline}\label{th:bootdiff6}
\sup_{h\in \text{BL}_1(\mathbb E)} |E^*[h(r_n\{\phi(\hat \theta_n^*) - \phi(\hat \theta_n)\})|\{X_i\}_{i=1}^n] - E[h(\phi_{\theta_0}^\prime(\mathbb G_0))]| \\ = \sup_{h\in \text{BL}_1(\mathbb E)} |E^*[h(\phi^\prime_{\theta_0}(\mathbb G_n^* + \mathbb G_n) - \phi^\prime_{\theta_0}(\mathbb G_n))|\{X_i\}_{i=1}^n] - E[h(\phi^\prime_{\theta_0}(\mathbb G_0))]| + o_p(1)
\end{multline}
Thus, in establishing the Theorem, it suffices to study the right hand side of \eqref{th:bootdiff6}.

\noindent \underline{First Claim:} We aim to show that (ii) implies (i). To this end, note by Lemma \ref{lm:auxGnindep}
\begin{equation}\label{th:bootdiff7}
(\phi^\prime_{\theta_0}(\mathbb G_n^* + \mathbb G_n) - \phi^\prime_{\theta_0}(\mathbb G_n),\mathbb G_n)\stackrel{L}{\rightarrow} (\phi^\prime_{\theta_0}(\mathbb G_1 + \mathbb G_2) - \phi^\prime_{\theta_0}(\mathbb G_2),\mathbb G_2)
\end{equation}
on $\mathbb E\times \mathbb D$ by the continuous mapping theorem. Let $f \in \text{BL}_1(\mathbb E)$ and $g\in \text{BL}_1(\mathbb D)$ satisfy $f(h_1) \geq 0$ and $g(h_2) \geq 0$ for any $h_1 \in \mathbb E$ and $h_2 \in \mathbb D$. By \eqref{th:bootdiff7} we then have:
\begin{equation}\label{th:bootdiff8}
\lim_{n\rightarrow \infty} E^*[f(\phi^\prime_{\theta_0}(\mathbb G_n^* + \mathbb G_n) - \phi^\prime_{\theta_0}(\mathbb G_n))g(\mathbb G_n)] = E[f(\phi^\prime_{\theta_0}(\mathbb G_1 + \mathbb G_2) - \phi^\prime_{\theta_0}(\mathbb G_2))g(\mathbb G_2)]
\end{equation}
On the other hand, also note that if the bootstrap is consistent, then result \eqref{th:bootdiff6} yields
\begin{equation}\label{th:bootdiff9}
\sup_{h\in \text{BL}_1(\mathbb E)} |E^*[h(\phi^\prime_{\theta_0}(\mathbb G_n^* + \mathbb G_n) - \phi^\prime_{\theta_0}(\mathbb G_n))|\{X_i\}_{i=1}^n] - E[h(\phi^\prime_{\theta_0}(\mathbb G_0))]| = o_p(1) ~.
\end{equation}
Moreover, since $\|g\|_{\infty} \leq 1$ and $\|f\|_\infty \leq 1$, it also follows that for any $\epsilon > 0$ we have:
\begin{align}\label{th:bootdiff10}
\lim_{n\rightarrow \infty} & E^*[|E^*[f(\phi^\prime_{\theta_0}(\mathbb G_n^* + \mathbb G_n) - \phi^\prime_{\theta_0}(\mathbb G_n))|\{X_i\}_{i=1}^n] - E[f(\phi^\prime_{\theta_0}(\mathbb G_0))]| g(\mathbb G_n)] \nonumber \\
& \leq \lim_{n \rightarrow \infty} E^*[|E^*[f(\phi^\prime_{\theta_0}(\mathbb G_n^* + \mathbb G_n) - \phi^\prime_{\theta_0}(\mathbb G_n))|\{X_i\}_{i=1}^n] - E[f(\phi_{\theta_0}^\prime(\mathbb G_0))]|] \nonumber \\
& \leq \lim_{n\rightarrow \infty} 2P^*(|E^*[f(\phi^\prime_{\theta_0}(\mathbb G_n^* + \mathbb G_n) - \phi^\prime_{\theta_0}(\mathbb G_n))|\{X_i\}_{i=1}^n] - E[f(\phi_{\theta_0}^\prime(\mathbb G_0))]| > \epsilon) + \epsilon ~.
\end{align}
Thus, result \eqref{th:bootdiff9}, $\epsilon$ being arbitrary in \eqref{th:bootdiff10}, Lemma \ref{lm:auxmeas}(v), $g(h) \geq 0$ for all $h \in \mathbb D$, and $\mathbb G_n \stackrel{L}{\rightarrow} \mathbb G_2$ by result \eqref{th:bootdiff7} allow us to conclude that:
\begin{multline}\label{th:bootdiff11}
\lim_{n \rightarrow \infty} E^*[E^*[f(\phi^\prime_{\theta_0}(\mathbb G_n^* + \mathbb G_n) - \phi^\prime_{\theta_0}(\mathbb G_n))|\{X_i\}_{i=1}^n]g(\mathbb G_n)] \\ = \lim_{n \rightarrow \infty} E^*[E[f(\phi_{\theta_0}^\prime(\mathbb G_0))]g(\mathbb G_n)] = E[f(\phi_{\theta_0}^\prime(\mathbb G_0))] E[g(\mathbb G_2)] ~.
\end{multline}
In addition, we also note that by Lemma 1.2.6 in \cite{Vaart1996}:
\begin{align}\label{th:bootdiff12}
\lim_{n\rightarrow \infty} E_*[f(\phi^\prime_{\theta_0}& (\mathbb G_n^* + \mathbb G_n) - \phi^\prime_{\theta_0}(\mathbb G_n))g(\mathbb G_n)]   \nonumber \\
& \leq \lim_{n\rightarrow \infty} E^*[E^*[f(\phi^\prime_{\theta_0}(\mathbb G_n^* + \mathbb G_n) - \phi^\prime_{\theta_0}(\mathbb G_n))|\{X_i\}_{i=1}^n] g(\mathbb G_n)] \nonumber \\
& \leq \lim_{n\rightarrow \infty} E^*[f(\phi^\prime_{\theta_0}(\mathbb G_n^* + \mathbb G_n) - \phi^\prime_{\theta_0}(\mathbb G_n))g(\mathbb G_n)]
\end{align}
since $\mathbb G_n$ is a function of $\{X_i\}_{i=1}^n$ only and $g(\mathbb G_n) \geq 0$. However, by \eqref{th:bootdiff7} and Lemma 1.3.8 in \cite{Vaart1996}, $(\phi^\prime_{\theta_0}(\mathbb G_n^* + \mathbb G_n) - \phi^\prime_{\theta_0}(\mathbb G_n),\mathbb G_n)$ is asymptotically measurable, and thus combining results \eqref{th:bootdiff11} and \eqref{th:bootdiff12} we can conclude:
\begin{equation}\label{th:bootdiff13}
\lim_{n\rightarrow \infty} E^*[f(\phi^\prime_{\theta_0}(\mathbb G_n^* + \mathbb G_n) - \phi^\prime_{\theta_0}(\mathbb G_n))g(\mathbb G_n)]
= E[f(\phi_{\theta_0}^\prime(\mathbb G_0))] E[g(\mathbb G_2)] ~.
\end{equation}
Hence, comparing \eqref{th:bootdiff8} and \eqref{th:bootdiff13} with $g \in \text{BL}_1(\mathbb D)$ given by $g(a) = 1$ for all $a \in \mathbb D$,
\begin{align}\label{th:bootdiff14}
E[f(\phi_{\theta_0}^\prime(\mathbb G_0))] E[g(\mathbb G_2)]  & = E[f(\phi_{\theta_0}^\prime(\mathbb G_1 + \mathbb G_2) - \phi^\prime_{\theta_0}(\mathbb G_2))] E[g(\mathbb G_2)] \nonumber \\ & =  E[f(\phi^\prime_{\theta_0}(\mathbb G_1 + \mathbb G_2) - \phi^\prime_{\theta_0}(\mathbb G_2))g(\mathbb G_2)] ~,
\end{align}
where the second equality follows again by \eqref{th:bootdiff8} and \eqref{th:bootdiff13}. Since \eqref{th:bootdiff14} must hold for any $f\in \text{BL}_1(\mathbb E)$ and $g\in \text{BL}_1(\mathbb D)$ with $f(h_1) \geq 0$ and $g(h_2) \geq 0$ for any $h_1 \in \mathbb E$ and $h_2 \in \mathbb D$, Lemma 1.4.2 in \cite{Vaart1996} implies $\phi_{\theta_0}^\prime(\mathbb G_1 + \mathbb G_2) - \phi^\prime_{\theta_0}(\mathbb G_2)$ must be independent of $\mathbb G_2$, and hence (i) must hold by Lemma \ref{lm:auxinvariance}.

\noindent \underline{Second Claim:} To conclude, we show (i) implies (ii). Fix $\epsilon > 0$ and note that by Assumption \ref{ass:paramest}, Lemma \ref{lm:auxGnboot}, and Lemma 1.3.8 in \cite{Vaart1996}, $\mathbb G_n$ and $\mathbb G_n^*$ are asymptotically tight. Hence, there is a compact set $K \subset \mathbb D$ such that:
\begin{equation}\label{th:bootdiff15}
\liminf_{n\rightarrow \infty} P_*(\mathbb G_n^* \in K^\delta) \geq 1-\epsilon \hspace{0.5 in} \liminf_{n\rightarrow \infty} P_*(\mathbb G_n \in K^\delta) \geq 1-\epsilon ~,
\end{equation}
for any $\delta > 0$ and $K^\delta \equiv \{a \in \mathbb D : \inf_{b \in K} \|a - b\|_{\mathbb D} < \delta\}$. Furthermore, by the Portmanteau Theorem we may assume without loss of generality that $K$ is a subset of the support of $\mathbb G_0$ and that $0 \in K$. Next, let $K + K \equiv \{a \in \mathbb D : a = b+c \text{ for some } b,c \in K\}$ and note that the compactness of $K$ implies $K+K$ is also compact. Thus, by Lemma \ref{lm:auxcont} and continuity of $\phi_{\theta_0}^\prime : \mathbb D \rightarrow \mathbb E$, there exist scalars $\delta_0 > 0$ and $\eta_0 > 0$ such that:
\begin{equation}\label{th:bootdiff16}
\sup_{a,b \in (K+K)^{\delta_0} : \|a - b\|_{\mathbb D} < \eta_0} \|\phi_{\theta_0}^\prime(a) - \phi_{\theta_0}^\prime(b)\|_{\mathbb E} < \epsilon ~.
\end{equation}
Next, for each $a \in K$, let $B_{\eta_0/2}(a) \equiv \{b \in \mathbb D : \|a - b\|_{\mathbb D} < \eta_0/2\}$. Since $\{B_{\eta_0/2}(a)\}_{a\in K}$ is an open cover of $K$, there exists a finite collection $\{B_{\eta_0/2}(a_j)\}_{j=1}^J$ also covering $K$. Therefore, since for any $b \in K^{\frac{\eta_0}{2}}$ there is a $\Pi b \in K$ such that $\|b - \Pi b\|_{\mathbb D} < \eta_0/2$, it follows that for every $b \in K^{\frac{\eta_0}{2}}$ there is a $1 \leq j \leq J$ such that $\|b - a_{j}\|_{\mathbb D} < \eta_0$. Setting $\delta_1 \equiv \min\{\delta_0,\eta_0\}/2$, we obtain that if $a \in K^{\delta_1}$ and $b\in K^{\delta_1}$, then: (i) $a+b \in (K+K)^{\delta_0}$ since $K^{\frac{\delta_0}{2}} + K^{\frac{\delta_0}{2}} \subseteq (K+K)^{\delta_0}$, (ii) there is a $1\leq j \leq J$ such that $\|b - a_j\|_{\mathbb D} < \eta_0$, and (iii) $(a+a_j) \in (K+K)^{\delta_0}$ since $a_j \in K$ and $a\in K^{\frac{\delta_0}{2}}$. Therefore, since $0 \in K$, we can conclude from \eqref{th:bootdiff16} that for every $b \in K^{\delta_1}$ there exists a $1 \leq j(b) \leq J$ such that
\begin{multline}\label{th:bootdiff17}
\sup_{a \in K^{\delta_1}} \|\{\phi_{\theta_0}^\prime(a + b) - \phi_{\theta_0}^\prime(b)\} - \{\phi_{\theta_0}^\prime(a + a_{j(b)}) - \phi_{\theta_0}^\prime(a_{j(b)})\}\|_{\mathbb E} \\ \leq \sup_{a,b\in (K + K)^{\delta_0}:\|a-b\|_{\mathbb D} < \eta_0}  2\|\phi_{\theta_0}^\prime(a) - \phi_{\theta_0}^\prime(b)\|_{\mathbb E} < 2\epsilon ~.
\end{multline}
In particular, if we define the set $\Delta_n \equiv \{\mathbb G_n^* \in K^{\delta_1}, \mathbb G_n \in K^{\delta_1}\}$, then \eqref{th:bootdiff17} implies that for every realization of $\mathbb G_n$ there is an $a_j$ independent of $\mathbb G_n^*$ such that:
\begin{multline}\label{th:bootdiff18}
\sup_{f\in \text{BL}_1(\mathbb E)} |(f(\phi^\prime_{\theta_0}(\mathbb G_n^* + \mathbb G_n) - \phi^\prime_{\theta_0}(\mathbb G_n)) - f(\phi^\prime_{\theta_0}(\mathbb G_n^* + a_j) - \phi^\prime_{\theta_0}(a_j)))1\{\Delta_n\}| < 2\epsilon ~.
\end{multline}
Letting $\Delta_n^c$ denote the complement of $\Delta_n$, result \eqref{th:bootdiff18} then allows us to conclude
\begin{multline}\label{th:bootdiff19}
\sup_{f\in \text{BL}_1(\mathbb E)} |E^*[f(\phi^\prime_{\theta_0}(\mathbb G_n^* + \mathbb G_n) - \phi^\prime_{\theta_0}(\mathbb G_n))|\{X_i\}_{i=1}^n] - E[f(\phi^\prime_{\theta_0}(\mathbb G_0))]| \leq 2P^*( \Delta_n^c|\{X_i\}_{i=1}^n)  \\
 + \max_{1\leq j \leq J} \sup_{f\in \text{BL}_1(\mathbb E)} |E^*[f(\phi^\prime_{\theta_0}(\mathbb G_n^* + a_j) - \phi^\prime_{\theta_0}(a_j))|\{X_i\}_{i=1}^n] - E[f(\phi^\prime_{\theta_0}(\mathbb G_0))]| + 2\epsilon
\end{multline}
since $\|f\|_\infty \leq 1$ for all $f\in \text{BL}_1(\mathbb E)$. However, by Assumptions \ref{ass:primboot}(i)-(ii) and \ref{ass:extraboot}(ii), and Theorem 10.8 in \cite{Kosorok2008} it follows that for any $1\leq j \leq J$:
\begin{equation}\label{th:bootdiff20}
\sup_{f\in \text{BL}_1(\mathbb E)} |E^*[f(\phi^\prime_{\theta_0}(\mathbb G_n^* + a_j) - \phi^\prime_{\theta_0}(a_j))|\{X_i\}_{i=1}^n] - E[f(\phi^\prime_{\theta_0}(\mathbb G_0 + a_j) - \phi^\prime_{\theta_0}(a_j))]| = o_p(1) ~.
\end{equation}
Thus, since  $K$ is a subset of the support of $\mathbb G_0$ and property (i) holds by hypothesis, result \eqref{th:bootdiff20}, the continuous mapping theorem, and $J < \infty$ allow us to conclude that:
\begin{equation}\label{th:bootdiff21}
\max_{1\leq j \leq J} \sup_{f\in \text{BL}_1(\mathbb E)} |E^*[f(\phi^\prime_{\theta_0}(\mathbb G_n^* + a_j) - \phi^\prime_{\theta_0}(a_j))|\{X_i\}_{i=1}^n] - E[f(\phi^\prime_{\theta_0}(\mathbb G_0))]|  = o_p(1) ~.
\end{equation}
Moreover, for any $\epsilon \in (0,1)$ we also have by Markov's inequality, Lemma 1.2.6 in \cite{Vaart1996}, $1\{\Delta_n^c\} \leq 1\{\mathbb G_n^* \notin K^{\delta_1}\} + 1\{\mathbb G_n \notin K^{\delta_1}\}$, and \eqref{th:bootdiff15} that:
\begin{multline}\label{th:bootdiff22}
\limsup_{n\rightarrow \infty} P^*(2P^*(\Delta_n^c|\{X_i\}_{i=1}^n) + 2\epsilon > 6 \sqrt \epsilon)  \leq \limsup_{n\rightarrow \infty} P^*(P^*(\Delta_n^c|\{X_i\}_{i=1}^n) > 2 \sqrt \epsilon) \\
\leq \frac{1}{2\sqrt \epsilon} \times \limsup_{n\rightarrow \infty} \{P^*(\mathbb G_n \notin K^{\delta_1}) + P^*(\mathbb G_n^* \notin K^{\delta_1}) \} \leq \sqrt \epsilon ~.
\end{multline}
Since $\epsilon > 0$ was arbitrary, combining \eqref{th:bootdiff6}, \eqref{th:bootdiff19}, \eqref{th:bootdiff21}, and \eqref{th:bootdiff22} imply (ii) holds, thus establishing the claim of the Theorem. \qed

\begin{lemma}\label{lm:auxGnboot}
If Assumptions \ref{ass:structure}(i), \ref{ass:paramest}(ii), \ref{ass:primboot}, \ref{ass:extraboot}(i) hold, then $r_n\{\hat \theta_n^* - \hat \theta_n\} \stackrel{L}{\rightarrow} \mathbb G_0$.
\end{lemma}

\noindent {\sc Proof:} In these arguments we need to distinguish between outer and inner expectations, and we therefore employ the notation $E^*$ and $E_*$ respectively. For notational simplicity also let $\mathbb G_n^* \equiv r_n\{\hat \theta_n^* - \hat \theta_n\}$. First, let $f \in \text{BL}_1(\mathbb D)$, and then note that by Lemma \ref{lm:auxmeas}(i) and Lemma 1.2.6 in \cite{Vaart1996} we have that:
\begin{align}\label{lm:auxGnboot1}
E^*[f(\mathbb G_n^*)] - E[f(\mathbb G_0)] & \geq E^*[E^*[f(\mathbb G_n^*)|\{X_i\}_{i=1}^n]] - E[f(\mathbb G_0)] \nonumber \\
& \geq -E^*[|E^*[f(\mathbb G_n^*)|\{X_i\}_{i=1}^n] - E[f(\mathbb G_0)]|] \nonumber \\  & \geq - E^*[\sup_{f\in \text{BL}_1(\mathbb D)}|E^*[f(\mathbb G_n^*)|\{X_i\}_{i=1}^n] - E[f(\mathbb G_0)]|] ~.
\end{align}
Similarly, applying Lemma 1.2.6 in \cite{Vaart1996} once again together with Lemma \ref{lm:auxmeas}(ii), and exploiting that $f\in \text{BL}_1(\mathbb D)$ we can conclude that:
\begin{align}\label{lm:auxGnboot2}
E_*[f(\mathbb G_n^*)] - E[f(\mathbb G_0)] & \leq E_*[E^*[f(\mathbb G_n^*)|\{X_i\}_{i=1}^n]] - E[f(\mathbb G_0)] \nonumber \\
& \leq E^*[|E^*[f(\mathbb G_n^*)|\{X_i\}_{i=1}^n] - E[f(\mathbb G_0)]|] \nonumber \\
& \leq E^*[\sup_{f\in \text{BL}_1(\mathbb D)}|E^*[f(\mathbb G_n^*)|\{X_i\}_{i=1}^n] - E[f(\mathbb G_0)]|] ~.
\end{align}
However, since $\|f\|_\infty \leq 1$ for all $f\in \text{BL}_1(\mathbb D)$, it also follows that for any $\eta > 0$ we have:
\begin{multline}\label{lm:auxGnboot3}
E^*[\sup_{f\in \text{BL}_1(\mathbb D)}| E^*[f(\mathbb G_n^*)|\{X_i\}_{i=1}^n] - E[f(\mathbb G_0)]|] \\ \leq 2P^*(\sup_{f\in \text{BL}_1(\mathbb D)}| E^*[f(\mathbb G_n^*)|\{X_i\}_{i=1}^n] - E[f(\mathbb G_0)]| > \eta) + \eta ~.
\end{multline}
Moreover, by Assumption \ref{ass:extraboot}(i), $E^*[f(\mathbb G_n^*)] = E_*[f(\mathbb G_n^*)] + o(1)$. Thus, Assumption \ref{ass:primboot}(ii), $\eta$ being arbitrary, and results \eqref{lm:auxGnboot1} and \eqref{lm:auxGnboot2} together imply that:
\begin{equation}\label{lm:auxGnboot4}
\lim_{n\rightarrow \infty} E^*[f(\mathbb G_n^*)] = E[f(\mathbb G_0)]
\end{equation}
for any $f\in \text{BL}_1(\mathbb D)$. Further note that since $\mathbb G_0$ is tight by Assumption \ref{ass:paramest}(ii) and $\mathbb D$ is a Banach space by Assumption \ref{ass:structure}(i), Lemma 1.3.2 in \cite{Vaart1996} implies $\mathbb G_0$ is separable. Therefore, the claim of the Lemma follows from \eqref{lm:auxGnboot4}, Theorem 1.12.2 and Addendum 1.12.3 in \cite{Vaart1996}. \qed

\begin{lemma}\label{lm:auxGnindep}
Let Assumptions \ref{ass:structure}(i), \ref{ass:paramest}, \ref{ass:primboot}, \ref{ass:extraboot}(i) hold, and $\mathbb G_1, \mathbb G_2 \in \mathbb D$ be independent random variables with the same law as $\mathbb G_0$. Then, it follows that on $\mathbb D\times \mathbb D$:
\begin{equation}
(r_n\{\hat \theta_n - \theta_0\}, r_n\{\hat \theta_n^* - \hat \theta_n\}) \stackrel{L}{\rightarrow} (\mathbb G_1, \mathbb G_2) ~.
\end{equation}
\end{lemma}

\noindent {\sc Proof:} In these arguments we need to distinguish between outer and inner expectations, and we therefore employ the notation $E^*$ and $E_*$ respectively. For notational convenience we also let $\mathbb G_n \equiv r_n\{\hat \theta_n - \theta_0\}$ and $\mathbb G_n^* \equiv r_n\{\hat \theta_n^* - \hat \theta_n\}$. Then, note that Assumptions \ref{ass:paramest}(i)-(ii), Lemma \ref{lm:auxGnboot}, and Lemma 1.3.8 in \cite{Vaart1996} imply that both $\mathbb G_n$ and $\mathbb G_n^*$ are asymptotically measurable, and asymptotically tight in $\mathbb D$. Therefore, by Lemma 1.4.3 in \cite{Vaart1996} $(\mathbb G_n,\mathbb G_n^*)$ is asymptotically tight in $\mathbb D\times\mathbb D$ and asymptotically measurable as well. Thus, by Prohorov's theorem (Theorem 1.3.9 in \cite{Vaart1996}), each subsequence $\{(\mathbb G_{n_k},\mathbb G_{n_k}^*)\}$ has an additional subsequence $\{(\mathbb G_{n_{k_j}},\mathbb G_{n_{k_j}}^*)\}$ such that:
\begin{equation}\label{lm:auxGnindep1}
(\mathbb G_{n_{k_j}}, \mathbb G_{n_{k_j}}^*) \stackrel{L}{\rightarrow} (\mathbb Z_1,\mathbb Z_2)
\end{equation}
for a tight Borel random variable $\mathbb Z\equiv(\mathbb Z_1,\mathbb Z_2)\in\mathbb D\times \mathbb D$. Since the sequence $\{(\mathbb G_{n_k},\mathbb G_{n_k}^*)\}$ was arbitrary, the Lemma follows if we show the law of $\mathbb Z$ equals that of $(\mathbb G_1,\mathbb G_2)$.

Towards this end, let $f_1, f_2 \in \text{BL}_1(\mathbb D)$ satisfy $f_1(h) \geq 0$ and $f_2(h) \geq 0$ for all $h\in \mathbb D$. Then note that by result \eqref{lm:auxGnindep1} it follows that:
\begin{equation}\label{lm:auxGnindep2}
\lim_{j\rightarrow \infty} E^*[f_1(\mathbb G_{n_{k_j}})f_2(\mathbb G_{n_{k_j}}^*)] = E[f_1(\mathbb Z_1)f_2(\mathbb Z_2)] ~.
\end{equation}
However, $f_1,f_2 \in \text{BL}_1(\mathbb D)$ satisfying $f_1(h) \geq 0$ and $f_2(h) \geq 0$ for all $h \in \mathbb D$, Lemma 1.2.6 in \cite{Vaart1996}, and Lemma \ref{lm:auxmeas}(iii) imply that:
\begin{align}\label{lm:auxGnindep3}
\lim_{j \rightarrow \infty} E^* [f_1(\mathbb G_{n_{k_j}})& f_2(\mathbb G_{n_{k_j}}^*)]  - E^*[f_1(\mathbb G_{n_{k_j}})E[f_2(\mathbb G_0)]] \nonumber \\
& \geq \lim_{j \rightarrow \infty} E^* [f_1(\mathbb G_{n_{k_j}}) E^*[f_2(\mathbb G_{n_{k_j}}^*)|\{X_i\}_{i=1}^n]]  - E^*[f_1(\mathbb G_{n_{k_j}})E[f_2(\mathbb G_0)]] \nonumber  \\
& \geq - \lim_{j \rightarrow \infty} E^*[f_1(\mathbb G_{n_{k_j}})|E^*[f_2(\mathbb G_{n_{k_j}}^*)|\{X_i\}_{i=1}^n] - E[f_2(\mathbb G_0)]|] \nonumber \\ & \geq - \lim_{j \rightarrow \infty} E^*[\sup_{f\in \text{BL}_1(\mathbb D)}|E^*[f(\mathbb G_{n_{k_j}}^*)|\{X_i\}_{i=1}^n] - E[f(\mathbb G_0)]|] ~,
\end{align}
where in the final inequality we exploited that $f_1 \in \text{BL}_1(\mathbb D)$. Similarly, Lemma 1.2.6 in \cite{Vaart1996}, Lemma \ref{lm:auxmeas}(iv), and $f_1,f_2 \in \text{BL}_1(\mathbb D)$ also imply that:
\begin{align}\label{lm:auxGnindep4}
\lim_{j \rightarrow \infty} E_* [f_1(\mathbb G_{n_{k_j}})& f_2(\mathbb G_{n_{k_j}}^*)]  - E_*[f_1(\mathbb G_{n_{k_j}})E[f_2(\mathbb G_0)]] \nonumber \\
& \leq \lim_{j \rightarrow \infty} E_* [f_1(\mathbb G_{n_{k_j}})E^*[f_2(\mathbb G_{n_{k_j}}^*)|\{X_i\}_{i=1}^n]]  - E_*[f_1(\mathbb G_{n_{k_j}})E[f_2(\mathbb G_0)]] \nonumber\\
& \leq \lim_{j \rightarrow \infty} E^*[f_1(\mathbb G_{n_{k_j}})|E^*[f_2(\mathbb G_{n_{k_j}}^*)|\{X_i\}_{i=1}^n] - E[f_2(\mathbb G_0)]|] \nonumber \\ & \leq \lim_{j \rightarrow \infty} E^*[\sup_{f\in \text{BL}_1(\mathbb D)}|E^*[f(\mathbb G_{n_{k_j}}^*)|\{X_i\}_{i=1}^n] - E[f(\mathbb G_0)]|] ~.
\end{align}
Thus, combining result \eqref{lm:auxGnboot3} together with \eqref{lm:auxGnindep3} and \eqref{lm:auxGnindep4}, and the fact that $(\mathbb G_n,\mathbb G_n^*)$ and $\mathbb G_n$ are asymptotically measurable, we can conclude that:
\begin{align}\label{lm:auxGnindep5}
\lim_{j \rightarrow \infty} E^* [f_1(\mathbb G_{n_{k_j}}) f_2(\mathbb G_{n_{k_j}}^*)] & = \lim_{j \rightarrow \infty} E^* [f_1(\mathbb G_{n_{k_j}})E[f_2(\mathbb G_0)]] \nonumber \\ & = E[f_1(\mathbb G_0)]E[f_2(\mathbb G_0)]  ~,
\end{align}
where the final result follows from $\mathbb G_n \stackrel{L}{\rightarrow} \mathbb G_0$ in $\mathbb D$. Hence, \eqref{lm:auxGnindep2} and \eqref{lm:auxGnindep5} imply
\begin{equation}\label{lm:auxGnindep6}
E[f_1(\mathbb Z_1)f_2(\mathbb Z_2)] = E[f_1(\mathbb G_0)]E[f_2(\mathbb G_0)]
\end{equation}
for all $f_1,f_2 \in \text{BL}_1(\mathbb D)$ satisfying $f_1(h)\geq 0$ and $f_2(h)\geq 0$ for all $h\in \mathbb D$. Since $\mathbb Z$ is tight on $\mathbb D\times \mathbb D$ it is also separable by Lemma 1.3.2 in \cite{Vaart1996} and Assumption \ref{ass:structure}(i), and hence result \eqref{lm:auxGnindep6} and Lemma 1.4.2 in \cite{Vaart1996} imply the law of $\mathbb Z$ equals that of $(\mathbb G_1,\mathbb G_2)$. In view of \eqref{lm:auxGnindep1}, the claim of the Lemma then follows. \qed

\begin{lemma}\label{lm:auxinvariance}
Let Assumptions \ref{ass:structure}(i)-(ii), \ref{ass:paramest}(ii) hold, $\mathbb D_L$ denote the support of $\mathbb G_0$, $0\in \mathbb D_L$, $\mathbb D_0 = \mathbb D_0 + \mathbb D_0$, and $\mathbb G_1$ be an independent copy of $\mathbb G_0$. If $\phi_{\theta_0}^\prime(\mathbb G_0 + \mathbb G_1) - \phi_{\theta_0}^\prime(\mathbb G_1)$ is independent of $\mathbb G_1$, then for any $a_0 \in \mathbb D_L$ and bounded continuous $f:\mathbb E\rightarrow \mathbf R$:
\begin{equation}\label{lm:auxinvariancedisp}
E[f(\phi_{\theta_0}^\prime(\mathbb G_0))] = E[f(\phi_{\theta_0}^\prime(\mathbb G_0 + a_0) - \phi_{\theta_0}^\prime(a_0))] ~.
\end{equation}
\end{lemma}

\noindent {\sc Proof:} We first note that since $\mathbb D_L \subseteq \mathbb D_0$ by Assumption \ref{ass:paramest}(ii) and $\mathbb G_1$ is independent of $\mathbb G_0$, it follows that the support of $\mathbb G_0 + \mathbb G_1$ is included in $\mathbb D_0 +\mathbb D_0 = \mathbb D_0$, and hence $\phi_{\theta_0}^\prime(\mathbb G_0 + \mathbb G_1)$ is well defined. Next, for any $a_0 \in \mathbb D$ and sequence $\{a_n\} \in \mathbb D$ with $\|a_0 - a_n\|_{\mathbb D} = o(1)$, we observe that continuity of $\phi_{\theta_0}^\prime$ and $f$, $f$ being bounded, and the dominated convergence theorem allow us to conclude that:
\begin{equation}\label{lm:auxinvariance1}
\lim_{n\rightarrow \infty} E[f(\phi_{\theta_0}^\prime(\mathbb G_0 + a_n) - \phi_{\theta_0}^\prime(a_n))] = E[f(\phi_{\theta_0}^\prime(\mathbb G_0 + a_0) - \phi_{\theta_0}^\prime(a_0))] ~.
\end{equation}
Hence, letting $B_\epsilon(a_0) \equiv \{a \in \mathbb D : \|a_0 - a\|_{\mathbb D} < \epsilon\}$, we note that result \eqref{lm:auxinvariance1} implies:
\begin{multline}\label{lm:auxinvariance2}
E[f(\phi_{\theta_0}^\prime(\mathbb G_0 + a_0) - \phi_{\theta_0}^\prime(a_0))] = \liminf_{\epsilon \downarrow 0} \inf_{a \in B_{\epsilon}(a_0)} E[f(\phi_{\theta_0}^\prime(\mathbb G_0 + a) - \phi_{\theta_0}^\prime(a))]  \\
\leq  \limsup_{\epsilon \downarrow 0} \sup_{a \in B_{\epsilon}(a_0)} E[f(\phi_{\theta_0}^\prime(\mathbb G_0 + a) - \phi_{\theta_0}^\prime(a))] = E[f(\phi_{\theta_0}^\prime(\mathbb G_0 + a_0) - \phi_{\theta_0}^\prime(a_0))] ~.
\end{multline}
Letting $L$ denote the law of $\mathbb G_0$, and for $\mathbb G_1$ and $\mathbb G_2$ independent copies of $\mathbb G_0$, we have:
\begin{align}\label{lm:auxinvariance3}
\inf_{a \in B_{\epsilon}(a_0)}  E[f(\phi_{\theta_0}^\prime &(\mathbb G_1 + a)  - \phi_{\theta_0}^\prime(a))]P(\mathbb G_2 \in B_{\epsilon}(a_0)) \nonumber \\ & \leq \int_{B_\epsilon(a_0)} \int_{\mathbb D_L} f(\phi_{\theta_0}^\prime(z_1 + z_2) - \phi_{\theta_0}^\prime(z_2)) d L(z_1)dL(z_2) \nonumber \\ & \leq \sup_{a \in B_{\epsilon}(a_0)} E[f(\phi_{\theta_0}^\prime(\mathbb G_1 + a) - \phi_{\theta_0}^\prime(a))]P(\mathbb G_2 \in B_{\epsilon}(a_0)) ~.
\end{align}
In particular, if $a_0 \in \mathbb D_L$, then $P(\mathbb G_2 \in B_{\epsilon}(a_0)) > 0$ for all $\epsilon > 0$, and thus we conclude:
\begin{multline}\label{lm:auxinvariance4}
E[f(\phi_{\theta_0}^\prime(\mathbb G_0 + a_0) - \phi_{\theta_0}^\prime(a_0))] = \lim_{\epsilon \downarrow 0} E[f(\phi_{\theta_0}^\prime(\mathbb G_1 + \mathbb G_2) - \phi_{\theta_0}^\prime(\mathbb G_2))| \mathbb G_2 \in B_\epsilon(a_0)]  \\ = \lim_{\epsilon \downarrow 0} E[f(\phi_{\theta_0}^\prime(\mathbb G_1 + \mathbb G_2) - \phi_{\theta_0}^\prime(\mathbb G_2))| \mathbb G_2 \in B_\epsilon(0)] = E[f(\phi_{\theta_0}^\prime(\mathbb G_0))] ~,
\end{multline}
where the first equality follows from \eqref{lm:auxinvariance2} and \eqref{lm:auxinvariance3}, the second by $\phi_{\theta_0}^\prime(\mathbb G_1 + \mathbb G_2) - \phi_{\theta_0}^\prime(\mathbb G_2)$ being independent of $\mathbb G_2$ by hypothesis, and the final equality follows by results \eqref{lm:auxinvariance2}, \eqref{lm:auxinvariance3}, and $\phi_{\theta_0}^\prime(0) = 0$ due to $\phi_{\theta_0}^\prime$ being homogenous of degree one. \qed

\begin{lemma}\label{lm:auxcont}
Let Assumption \ref{ass:structure}(i) hold, $\psi : \mathbb D \rightarrow \mathbb E$ be continuous, and $K \subset \mathbb D$ be compact. It then follows that for every $\epsilon > 0$ there exist $\delta > 0, \eta > 0$ such that:
\begin{equation}\label{lm:auxcontdisp1}
\sup_{(a,b) \in K^\delta \times K^\delta : \|a - b\|_{\mathbb D} < \eta}  \|\psi(a) - \psi(b)\|_{\mathbb E} < \epsilon ~.
\end{equation}
\end{lemma}

\noindent {\sc Proof:} Fix $\epsilon > 0$ and note that since $\psi : \mathbb D \rightarrow \mathbb E$ is continuous, it follows that for every $a \in \mathbb D$ there exists a $\zeta_a$ such that $\|\psi(a) - \psi(b)\|_{\mathbb E} < \epsilon/2$ for all $b\in \mathbb D$ with $\|a - b\|_{\mathbb D} < \zeta_a$. Letting $B_{\zeta_a/4}(a) \equiv \{b \in \mathbb D : \|a - b\|_{\mathbb D} < \zeta_a/4\}$, then observe that $\{B_{\zeta_a/4}(a)\}_{a\in K}$ forms an open cover of $K$ and hence, by compactness of $K$, there exists a finite subcover $\{B_{\zeta_{a_j}/4}(a_j)\}_{j=1}^J$ for some $J < \infty$. To establish the Lemma, we then let
\begin{equation}\label{lm:auxcont1}
\eta \equiv \min_{1\leq j \leq J} \frac{\zeta_{a_j}}{4} \hspace{0.5 in} \delta \equiv \min_{1\leq j \leq J} \frac{\zeta_{a_j}}{4} ~.
\end{equation}
For any $a \in K^\delta$, there then exists a $\Pi a \in K$ such that $\|a - \Pi a\|_{\mathbb D} < \delta$, and since $\{B_{\zeta_{a_j}/4}(a_j)\}_{j=1}^J$ covers $K$, there also is a $\bar j$ such that $\Pi a \in B_{\zeta_{a_{\bar j}/4}}(a_{\bar j})$. Thus, we have
\begin{equation}\label{lm:auxcont2}
\|a - a_{\bar j}\|_{\mathbb D} \leq \|a - \Pi a\|_{\mathbb D} + \|\Pi a - a_{\bar j}\|_{\mathbb D} < \delta + \frac{\zeta_{a_{\bar j}}}{4} \leq \frac{\zeta_{a_{\bar j}}}{2} ~,
\end{equation}
due to the choice of $\delta$ in \eqref{lm:auxcont1}. Moreover, if $b \in \mathbb D$ satisfies $\|a-b\|_{\mathbb D} < \eta$, then:
\begin{equation}\label{lm:auxcont3}
\|b - a_{\bar j}\|_{\mathbb D} \leq \|a - b\|_{\mathbb D} + \|a - a_{\bar j}\|_{\mathbb D} < \eta + \frac{\zeta_{a_{\bar j}}}{2} \leq \zeta_{a_{\bar j}} ~,
\end{equation}
by the choice of $\eta$ in \eqref{lm:auxcont1}. We conclude from \eqref{lm:auxcont2}, \eqref{lm:auxcont3} that $a,b \in B_{\zeta_{a_{\bar j}}}(a_{\bar j})$, and
\begin{equation}\label{lm:auxcont4}
\|\psi(a) - \psi(b)\|_{\mathbb E} \leq \|\psi(a) - \psi(a_{\bar j})\|_{\mathbb E} + \|\psi(b) - \psi(a_{\bar j})\|_{\mathbb E} < \frac{\epsilon}{2} + \frac{\epsilon}{2} = \epsilon
\end{equation}
by our choice of $\zeta_{a_{\bar j}}$. Thus, the Lemma follows from result \eqref{lm:auxcont4}. \qed

\begin{lemma}\label{lm:auxmeas}
Let $(\Omega, \mathcal F, P)$ be a probability space, $c \in \mathbf R_+$ and $U:\Omega \rightarrow \mathbf R$ and $V:\Omega \rightarrow \mathbf R$ be arbitrary maps satisfying $U(\omega) \geq 0$ and $V(\omega) \geq 0$ for all $\omega \in \Omega$. If $E^*$ and $E_*$ denote outer and inner expectations respectively, then it follows that:

\vspace{-0.2 in}
\begin{packed_enum}
    \item[(i)] $E^*[U] - c \geq - E^*[|U - c|]$.
    \item[(ii)] $E_*[U] - c \leq E^*[|U - c|]$.
    \item[(iii)] $E^*[UV] - E^*[Uc] \geq - E^*[U|V-c|]$ whenever $\min\{E^*[UV], E^*[Uc]\} < \infty$.
    \item[(iv)] $E_*[UV] - E_*[Uc] \leq E^*[U|V-c|]$ whenever $\min\{E_*[UV], E_*[Uc]\} < \infty$.
    \item[(v)] $|E^*[UV] - E^*[Uc]| \leq E^*[U|V-c|]$ whenever $\min\{E_*[UV], E_*[Uc]\} < \infty$.
\end{packed_enum}
\end{lemma}

\noindent {\sc Proof:} The arguments are simple and tedious, but unfortunately necessary to address the possible nonlinearity of inner and outer expectations. Throughout, for a map $T:\Omega \rightarrow \mathbf R$, we let $T^*$ and $T_*$ denote the minimal measurable majorant and the maximal measurable minorant of $T$ respectively. We will also exploit the fact that:
\begin{equation}\label{lm:auxmeas1}
E_*[T] = - E^*[-T] ~,
\end{equation}
and that $E^*[T] = E[T^*]$ whenever $E[T^*]$ exists, which in the context of this Lemma is always satisfied since all variables are positive.

To establish the first claim of the Lemma, note that Lemma 1.2.2(i) in \cite{Vaart1996} implies $U^* - c = (U - c)^*$. Therefore, \eqref{lm:auxmeas1} and $E_* \leq E^*$ yield:
\begin{multline}\label{lm:auxmeas2}
E^*[U] - c = E[U^* - c] = E[(U - c)^*] = E^*[U- c] \\ \geq E^*[-|U-c|] = - E_*[|U- c|] \geq - E^*[|U- c|] ~.
\end{multline}
Similarly, for the second claim of the Lemma, exploit that $E_* \leq E^*$, and once again employ Lemma 1.2.2(i) in \cite{Vaart1996} to conclude that:
\begin{equation}\label{lm:auxmeas3}
E_*[U] - c \leq   E^*[U] - c = E[U^* - c] = E[(U - c)^*] \leq E^*[|U -c |] ~.
\end{equation}
For the third claim, note that Lemma 1.2.2(iii) in \cite{Vaart1996} implies $|(UV)^* - (Uc)^*| \leq |UV - Uc|^*$. Thus, since $|U(V-c)| = U|V-c|$ as a result of $U(\omega) \geq 0$ for all $\omega \in \Omega$, we obtain from relationship \eqref{lm:auxmeas1} and $E_*\leq E^*$ that:
\begin{multline}\label{lm:auxmeas4}
E^*[UV] - E^*[Uc] = E[(UV)^* - (Uc)^*] \geq E[-|(UV)^* - (Uc)^*|] \\ \geq E[-|UV - Uc|^*]
 = -E_*[U|V - c|] \geq -E^*[U|V - c|] ~.
\end{multline}
Similarly, for the fourth claim of the Lemma, employ \eqref{lm:auxmeas1}, that $|(-Uc)^* - (-UV)^*| \leq |(-Uc) - (-UV)|^*$ by Lemma 1.2.2(iii) in \cite{Vaart1996}, and that $|UV - Uc| = U|V-c|$ due to $U(\omega) \geq 0$ for all $\omega \in \Omega$ to obtain that:
\begin{multline}\label{lm:auxmeas5}
E_*[UV] - E_*[Uc] = E[(-Uc)^* - (-UV)^*]  \leq E[|(-Uc)^* - (-UV)^*|] \\ \leq E[|(-Uc) - (-UV)|^*] = E^*[U|V - c|] ~.
\end{multline}
Finally, for the fifth claim of the Lemma, note the same arguments as in \eqref{lm:auxmeas5} yield
\begin{multline}\label{lm:auxmeas6}
E^*[UV] - E^*[Uc] = E[(Uc)^* - (UV)^*]  \leq E[|(Uc)^* - (UV)^*|] \\ \leq E[|(Uc) - (UV)|^*] = E^*[U|V - c|] ~.
\end{multline}
Thus, part (v) of the Lemma follows from part (iii) and \eqref{lm:auxmeas6}. \qed

\begin{lemma}\label{lm:simpsuff}
Let Assumption \ref{ass:structure} hold, and suppose that for some $\kappa > 0$ and $C_0 < \infty$ we have $\|\hat \phi^\prime_n(h_1) - \hat \phi^\prime_n(h_2)\|_{\mathbb E} \leq C_0 \|h_1-h_2\|_{\mathbb D}^\kappa$ for all $h_1,h_2\in \mathbb D$ outer almost surely. Then, Assumption \ref{ass:derest} holds provided that for all $h \in \mathbb D_0$ we have:
\begin{equation}\label{lm:simpsuffdisp1}
\|\hat \phi_n^\prime(h) - \phi^\prime_{\theta_0}(h)\|_{\mathbb E} = o_p(1) ~.
\end{equation}
\end{lemma}

\noindent {\sc Proof:} Fix $\epsilon > 0$, let $K_0 \subseteq \mathbb D_0$ be compact, and for any $h\in \mathbb D$ let $\Pi : \mathbb D \rightarrow K_0$ satisfy $\|h - \Pi h\|_{\mathbb D} = \inf_{a \in K_0}\|h - a\|_{\mathbb D}$ -- here attainment is guaranteed by compactness. Since $\phi_{\theta_0}^\prime :\mathbb D\rightarrow \mathbb E$ is continuous, Lemma \ref{lm:auxcont} implies there exists a $\delta_1 > 0$ such that:
\begin{equation}\label{lm:simpsuff1}
\sup_{h \in K^{\delta_1}_0} \|\phi_{\theta_0}^\prime(h) - \phi_{\theta_0}^\prime(\Pi h)\|_{\mathbb E} < \epsilon ~.
\end{equation}
Next, set $\delta_2 < (\epsilon/C_0)^{1/\kappa}$ and note that by hypothesis we have outer almost surely that:
\begin{equation}\label{lm:simpsuff2}
\sup_{h \in K^{\delta_2}_0} \|\hat \phi_n^\prime(h) - \hat \phi_{n}^\prime(\Pi h)\|_{\mathbb E} \leq \sup_{h \in K^{\delta_2}_0} C_0\|h - \Pi h\|_{\mathbb E}^{\kappa} \leq C_0\delta_2^\kappa < \epsilon ~.
\end{equation}
Defining $\delta_3 \equiv \min\{\delta_1,\delta_2\}$, exploiting \eqref{lm:simpsuff1}, \eqref{lm:simpsuff2}, and $\Pi h\in K_0$ we then conclude:
\begin{align}\label{lm:simpsiff3}
\sup_{h \in K^{\delta_3}_0}& \|\hat \phi_n^\prime(h) - \phi^\prime_{\theta_0}(h)\|_{\mathbb E} \nonumber \\ & \leq \sup_{h \in K^{\delta_3}_0}\{\|\hat \phi_n^\prime(h) - \hat \phi_n^\prime(\Pi h)\|_{\mathbb E} +  \| \phi_{\theta_0}^\prime(h) - \phi_{\theta_0}^\prime(\Pi h)\|_{\mathbb E} + \|\hat \phi_n^\prime (\Pi h) - \phi_{\theta_0}^\prime(\Pi h)\|_{\mathbb E} \}\nonumber \\ & \leq \sup_{h \in K_0} \|\hat \phi_n^\prime(h) - \phi_{\theta_0}^\prime(h)\|_{\mathbb E} + 2\epsilon
\end{align}
outer almost surely. Thus, since $K_0^\delta \subseteq K_0^{\delta_3}$ for all $\delta \leq \delta_3$ we obtain from \eqref{lm:simpsiff3} that:
\begin{multline}\label{lm:simsiff4}
\lim_{\delta\downarrow 0} \limsup_{n\rightarrow \infty}  P(\sup_{h \in K_0^{\delta}} \|\hat \phi_n^\prime(h) - \phi_{\theta_0}^\prime(h)\|_{\mathbb E} > 5\epsilon ) \\
\leq \limsup_{n\rightarrow \infty} P(\sup_{h \in K_0} \|\hat \phi_n^\prime(h) - \phi_{\theta_0}^\prime(h)\|_{\mathbb E} > 3\epsilon )  ~.
\end{multline}
Next note that since $K_0$ is compact, $\phi_{\theta_0}^\prime$ is uniformly continuous on $K_0$, and thus we can find a finite collection $\{h_j\}_{j=1}^J$ with $J < \infty$ such that $h_j \in K_0$ for all $j$ and:
\begin{equation}\label{lm:simsiff5}
\sup_{h\in K_0} \min_{1\leq j \leq J} \max\{C_0\|h - h_j\|_{\mathbb D}^\kappa, \|\phi_{\theta_0}^\prime(h) - \phi_{\theta_0}^\prime(h_j)\|_{\mathbb E}\} < \epsilon ~.
\end{equation}
In particular, since $\|\hat \phi_{\theta_0}^\prime(h) - \hat \phi_{\theta_0}^\prime(h_j)\|_{\mathbb E} \leq C_0\|h - h_j\|_{\mathbb D}^\kappa$, we obtain from \eqref{lm:simsiff5} that:
\begin{equation}\label{lm:simsiff6}
\sup_{h\in K_0} \|\hat \phi_{\theta_0}^\prime(h) - \phi_{\theta_0}^\prime(h)\|_{\mathbb E} \leq \max_{1\leq j \leq J} \|\hat \phi_{\theta_0}^\prime(h_j) - \phi_{\theta_0}^\prime(h_j)\|_{\mathbb E} + 2\epsilon ~.
\end{equation}
Thus, we can conclude from \eqref{lm:simsiff6} and $\hat \phi_{\theta_0}^\prime$ satisfying \eqref{lm:simpsuffdisp1} for any $h \in \mathbb D_0$ that:
\begin{multline}\label{lm:simsiff7}
\limsup_{n\rightarrow \infty} P(\sup_{h \in K_0} \|\hat \phi_n^\prime(h) - \phi_{\theta_0}^\prime(h)\|_{\mathbb E} > 3\epsilon ) \\
\leq \limsup_{n\rightarrow \infty} P(\max_{1\leq j \leq J} \|\hat \phi_n^\prime(h_j) - \phi_{\theta_0}^\prime(h_j)\|_{\mathbb E} > \epsilon ) = 0 ~.
\end{multline}
Since $\epsilon$ and $K_0$ were arbitrary, the Lemma follows from \eqref{lm:simsiff4} and \eqref{lm:simsiff5}. \qed

\begin{lemma}\label{aux:supp}
Let Assumptions \ref{ass:structure}(i)-(ii), \ref{ass:paramest}(ii) hold, and $\mathbb G_0$ be a Gaussian measure. If the support of $\mathbb G_0$ is a vector subspace of $\mathbb D$, then it is also a separable Banach space under $\|\cdot\|_{\mathbb D}$ and it includes the mean of $\mathbb G_0$.
\end{lemma}

\noindent {\sc Proof:} By Assumption \ref{ass:structure} and Theorem 7.1.7 in \cite{bogachev2:2007} it follows that $\mathbb G_0$ is regular. Hence, since in addition $\mathbb G_0$ is tight by Assumption \ref{ass:paramest}(ii), we can further conclude that $\mathbb G_0$ is Radon. Letting $\mathbb D_L$ and $\mu_0$ respectively denote the support and the mean of $\mathbb G_0$, we then obtain by Theorem 3.6.1 in \cite{bogachevgauss} that
\begin{equation}\label{aux:supp1}
\mathbb D_L = \mu_0 + \mathbb D_A ~,
\end{equation}
where $\mathbb D_A$ is a closed separable subspace of $\mathbb D$. However, since $\mathbb D_L$ is also a vector subspace of $\mathbb D$ by hypothesis, it follows that $\mu_0 \in \mathbb D_A$ and hence $\mathbb D_L = \mathbb D_A$ and $\mu_0 \in \mathbb D_L$. Moreover, since $\mathbb D_A$ is a closed separable subspace of $\mathbb D$, we further conclude $\mathbb D_A$ is a separable Banach space under $\|\cdot\|_{\mathbb D}$ and the Lemma follows from $\mathbb D_L = \mathbb D_A$. \qed


\begin{center}
{\Large {\sc \hypertarget{examplesupp}{Appendix B}} - Results for Examples \ref{ex:absmean}-\ref{ex:lrorder}}
\end{center}

\renewcommand{\theequation}{B.\arabic{equation}}
\renewcommand{\thelemma}{B.\arabic{lemma}}
\renewcommand{\thecorollary}{B.\arabic{corollary}}
\renewcommand{\thetheorem}{B.\arabic{theorem}}
\setcounter{lemma}{0}
\setcounter{theorem}{0}
\setcounter{corollary}{0}
\setcounter{equation}{0}
\setcounter{remark}{0}

\begin{lemma}\label{lm:auxmax}
Let $\mathbf A$ be totally bounded under a norm $\|\cdot\|_{\mathbf A}$, and $\bar {\mathbf A}$ denote its closure under $\|\cdot\|_{\mathbf A}$. Further let $\phi : \ell^\infty(\mathbf A) \rightarrow \mathbf R$ be given by $\phi(\theta) = \sup_{a\in \mathbf A} \theta(a)$, and define $\Psi_{\bar{\mathbf A}}(\theta) \equiv \arg\max_{a\in \bar {\mathbf A}} \theta(a)$ for any $\theta \in \mathcal C(\bar{\mathbf A})$. Then, $\phi$ is Hadamard directionally differentiable tangentially to $\mathcal C(\bar{ \mathbf A})$ at any $\theta \in \mathcal C(\bar{\mathbf A})$, and $\phi_\theta^\prime :\mathcal C(\bar{\mathbf A}) \rightarrow \mathbf R$ satisfies: $$\phi_\theta^\prime(h) = \sup_{a\in\Psi_{\bar {\mathbf A}}(\theta)} h(a) \quad \quad h\in \mathcal C(\bar {\mathbf A}) ~.$$
\end{lemma}

\noindent {\sc Proof:} First note Corollary 3.29 in \cite{aliprantis:border:2006} implies $\bar {\mathbf A}$ is compact under $\|\cdot\|_{\mathbf A}$. Next, let $\{t_n\}$ and $\{h_n\}$ be sequence with $t_n \in \mathbf R$, $h_n \in \ell^\infty(\mathbf A)$ for all $n$ and $\|h_n-h\|_\infty = o(1)$ for some $h \in \mathcal C(\bar{\mathbf A})$. Then note that for any $\theta \in \mathcal C(\bar{\mathbf A})$ we have:
\begin{equation}\label{lm:auxmax1}
|\sup_{a \in \mathbf A} \{\theta(a) + t_n h_n(a)\} - \sup_{a \in \mathbf A} \{\theta(a) + t_n h(a)\}| \leq t_n \|h_n - h\|_\infty = o(t_n) ~.
\end{equation}
Further note that since $\bar{\mathbf A}$ is compact, $\Psi_{\bar {\mathbf A}}(\theta)$ is well defined for any $\theta \in \mathcal C(\bar{\mathbf A})$. Defining $\Gamma_\theta : \mathcal C(\bar{\mathbf A}) \rightarrow \mathcal C(\bar{\mathbf A})$ to be given by $\Gamma_\theta(g) = \theta + g$, then note that $\Gamma_\theta$ is trivially continuous. Therefore, Theorem 17.31 in \cite{aliprantis:border:2006} and the relation
\begin{equation}\label{lm:auxmax2}
\Psi_{\bar {\mathbf A}}(\theta + g) = \arg\max_{a \in \bar{\mathbf A}} \Gamma_\theta(g)(a)
\end{equation}
imply that $\Psi_{\bar {\mathbf A}}(\theta + g)$ is upper hemicontinuous in $g$. In particular, for $\Psi_{\bar {\mathbf A}}(\theta)^\epsilon \equiv \{a \in \bar{\mathbf A }: \inf_{a_0 \in \Psi_{\bar {\mathbf A}}(\theta)} \|a - a_0\|_{\mathbf A} \leq \epsilon\}$, it follows from $\|t_n h\|_\infty = o(1)$ that $\Psi_{\bar {\mathbf A}}(\theta +t_n h) \subseteq \Psi_{\bar {\mathbf A}}(\theta)^{\delta_n}$ for some $\delta_n \downarrow 0$. Thus, since $\Psi_{\bar {\mathbf A}}(\theta) \subseteq \Psi_{\bar {\mathbf A}}(\theta)^{\delta_n}$ we can conclude that
\begin{align}\label{lm:auxmax3}
|\sup_{a\in \bar{\mathbf A}} \{\theta(a) + t_n h(a)\} - & \sup_{a\in \Psi_{\bar {\mathbf A}}(\theta)} \{\theta(a) + t_n h(a)\}| \nonumber \\ & = \sup_{a\in \Psi_{\bar {\mathbf A}}(\theta)^{\delta_n}} \{\theta(a) + t_n h(a)\} - \sup_{a\in \Psi_{\bar {\mathbf A}}(\theta)} \{\theta(a) + t_n h(a)\} \nonumber \\ & \leq \sup_{a_0,a_1 \in \bar{\mathbf A }: \|a_0 - a_1\|_{\mathbf A} \leq \delta_n} t_n|h(a_0) - h(a_1)| \nonumber \\ &= o(t_n) ~,
\end{align}
where the final result follows from $h$ being uniformly continuous by compactness of $\bar{\mathbf A}$. Therefore, exploiting \eqref{lm:auxmax1}, \eqref{lm:auxmax3} and $\theta$ being constant on $\Psi_{\bar {\mathbf A}}(\theta)$ yields
\begin{multline}\label{lm:auxmax4}
 |\sup_{a\in \mathbf A} \{\theta(a) + t_nh_n(a)\} - \sup_{a \in \mathbf A} \theta(a) - t_n \sup_{a \in \Psi_{\bar {\mathbf A}}(\theta)} h(a)| \\ \leq |\sup_{a\in \Psi_{\bar {\mathbf A}}(\theta)} \{\theta(a) + t_nh(a)\} - \sup_{a \in \Psi_{\bar {\mathbf A}}(\theta)} \theta(a) - t_n \sup_{a \in \Psi_{\bar {\mathbf A}}(\theta)} h(a)| + o(t_n) = o(t_n) ~,
\end{multline}
which verifies the claim of the Lemma. \qed

\begin{lemma}\label{lm:auxstochdiff}
Let $w:\mathbf R \rightarrow \mathbf R_+$ satisfy $\int_{\mathbf R} w(u)du < \infty$ and $\phi:\ell^\infty(\mathbf R)\times \ell^\infty(\mathbf R) \rightarrow \mathbf R$ be given by $\phi(\theta) = \int_{\mathbf R}\max\{\theta^{(1)}(u) - \theta^{(2)}(u),0\}w(u)du$ for any $\theta = (\theta^{(1)},\theta^{(2)}) \in \ell^\infty(\mathbf R)\times \ell^\infty(\mathbf R)$. Then, $\phi$ is Hadamard directionally differentiable at any $\theta \in \ell^\infty(\mathbf R)\times \ell^\infty(\mathbf R)$ with $\phi_\theta^\prime:\ell^\infty(\mathbf R)\times \ell^\infty(\mathbf R) \rightarrow \mathbf R$ satisfying for any $h = (h^{(1)},h^{(2)})\in \ell^\infty(\mathbf R)\times \ell^\infty(\mathbf R)$
$$\phi^\prime_\theta(h) = \int_{B_0(\theta)}\max\{h^{(1)}(u) - h^{(2)}(u),0\}w(u)du + \int_{B_+(\theta)}(h^{(1)}(u)-h^{(2)}(u))w(u)du ~, $$
where $B_+(\theta) \equiv \{u \in \mathbf R: \theta^{(1)}(u) > \theta^{(2)}(u)\}$ and $B_0(\theta) \equiv \{u \in \mathbf R: \theta^{(1)}(u) = \theta^{(2)}(u)\}$.
\end{lemma}

\noindent {\sc Proof:} Let $\{h_n\} = \{(h_{n}^{(1)},h_{n}^{(2)})\}$ be a sequence in $\ell^\infty(\mathbf R) \times \ell^\infty(\mathbf R)$ satisfying $\|h_{n}^{(1)} - h^{(1)}\|_\infty  \vee \|h_{n}^{(2)} - h^{(2)}\|_\infty = o(1)$ for some $h = (h^{(1)},h^{(2)}) \in\ell^\infty(\mathbf R)\times \ell^\infty(\mathbf R)$. Further let:
\begin{equation}\label{lm:auxstochdiff1}
B_-(\theta) \equiv \{u \in \mathbf R : \theta^{(1)}(u) < \theta^{(2)}(u)\} ~.
\end{equation}
Next, observe that since $\theta^{(1)}(u) - \theta^{(2)}(u) < 0$ for all $u \in B_-(\theta)$, and $\|h_{n}^{(1)} - h_{n}^{(2)}\|_\infty = O(1)$ due to $\|h^{(1)}- h^{(2)}\|_\infty < \infty$, the dominated convergence theorem yields that:
\begin{multline}\label{lm:auxstochdiff2}
\int_{B_-(\theta)}\max\{(\theta^{(1)}(u) - \theta^{(2)}(u)) + t_n(h_{n}^{(1)}(u) - h_{n}^{(2)}(u)),0\}w(u)du \\ \lesssim t_n  \int_{B_-(\theta)} 1\{t_n(h_{n}^{(1)}(u) - h_{n}^{(2)}(u)) \geq - (\theta^{(1)}(u) - \theta^{(2)}(u))\}w(u)du = o(t_n) ~.
\end{multline}
Thus, \eqref{lm:auxstochdiff2}, $B_-(\theta)^c = B_+(\theta) \cup B_0(\theta)$ and the dominated convergence theorem imply
\begin{multline*}
\frac{1}{t_n}\{\phi(\theta + t_nh_n) - \phi(\theta)\}  \\  = \int_{B_-(\theta)^c}\max\{h_{n}^{(1)}(u)-h_{n}^{(2)}(u), - \frac{\theta(u)^{(1)} - \theta^{(2)}(u)}{t_n}\}w(u) du + o(1) =\phi_\theta^\prime(h) + o(1)
\end{multline*}
which establishes the claim of the Lemma. \qed

\begin{lemma}\label{lm:auxsupest}
Let Assumption \ref{ass:structure} hold, and $\mathbf A$ be compact under $\|\cdot\|_{\mathbf A}$. Further suppose $\phi : \ell^\infty(\mathbf A) \rightarrow \mathbf R$ is Hadamard directionally differentiable tangentially to $\mathcal C(\mathbf A)$ at $\theta_0 \in \mathcal C(\mathbf A)$, and that for some $A_0 \subseteq \mathbf A$, its derivative $\phi_{\theta_0}^\prime:\mathcal C(\mathbf A) \rightarrow \mathbf R$ is given by:
\begin{equation}\label{lm:auxsupestdisp1}
\phi_{\theta_0}^\prime (h) = \sup_{a \in A_0} h(a) ~.
\end{equation}
If $\hat A_0 \subseteq \mathbf A$ outer almost surely, and $d_H(\hat A_0,A_0,\|\cdot\|_{\mathbf A}) = o_p(1)$, then it follows that $\hat \phi_n^\prime : \ell^\infty(\mathbf A) \rightarrow \mathbf R$ given by $\hat \phi_n^\prime(h) = \sup_{a\in \hat A_0} h(a)$ for any $h\in \ell^\infty(\mathbf A)$ satisfies \eqref{ass:derestdisp1}.
\end{lemma}

\noindent {\sc Proof:} First note that $\hat \phi_n^\prime$ is outer almost surely Lipschitz since $|\hat \phi_n^\prime(h_1) - \hat \phi_n^\prime(h_2)| \leq \|h_1 - h_2\|_{\infty}$ for all $h_1,h_2\in \ell^\infty(\mathbf A)$ due to $\hat A_0 \subseteq \mathbf A$ outer almost surely. Therefore, by Lemma \ref{lm:simpsuff} it suffices to verify that for any $h \in \mathcal C(\mathbf A)$, $\hat \phi_n^\prime$ satisfies
\begin{equation}\label{lm:auxsupest1}
 |\hat \phi_n^\prime(h) - \phi^\prime_{\theta_0}(h)| = o_p(1) ~.
\end{equation}
Towards this end, fix an arbitrary $\epsilon_0 > 0$ and note $h$ is uniformly continuous on $\mathbf A$ due to $\mathbf A$ being compact. Hence, we conclude there exists an $\eta > 0$ such that
\begin{equation}\label{lm:auxsupest2}
\sup_{\|a_1 - a_2\|_{\mathbf A} < \eta} |h(a_1) - h(a_2)| < \epsilon_0 ~.
\end{equation}
Moreover, given the definitions of $\hat \phi_n^\prime$ and $\phi_{\theta_0}^\prime$ it also follows that for any $h \in \ell^\infty(\mathbf A)$:
\begin{equation}\label{lm:auxsupest3}
|\hat \phi_n^\prime(h) - \phi^\prime_{\theta_0}(h)| \leq \sup_{\|a_1-a_2\|_{\mathbf A}\leq d_H(\hat A_0,A_0,\|\cdot\|_{\mathbf A})} |h(a_1) - h(a_2)| ~.
\end{equation}
Thus, by results \eqref{lm:auxsupest2} and \eqref{lm:auxsupest3}, and the Hausdorff consistency of $\hat A_0$, we obtain:
\begin{equation}\label{lm:auxsupest4}
\limsup_{n\rightarrow \infty} P(|\hat \phi_n^\prime (h) - \phi_{\theta_0}^\prime(h)| > \epsilon_0) \leq \limsup_{n\rightarrow \infty} P(d_H(\hat A_0,A_0,\|\cdot\|_{\mathbf A}) > \eta) = 0 ~.
\end{equation}
It follows that \eqref{lm:auxsupest1} indeed holds, and the claim of the Lemma follows. \qed


\begin{center}
{\Large {\sc \hypertarget{examplesupp}{Appendix C}} - Results for Section \ref{sec:convex}}
\end{center}

\renewcommand{\theequation}{C.\arabic{equation}}
\renewcommand{\thelemma}{C.\arabic{lemma}}
\renewcommand{\thecorollary}{C.\arabic{corollary}}
\renewcommand{\thetheorem}{C.\arabic{theorem}}
\setcounter{lemma}{0}
\setcounter{theorem}{0}
\setcounter{corollary}{0}
\setcounter{equation}{0}
\setcounter{remark}{0}

\noindent {\sc Proof of Proposition \ref{prop:convexdist}:} We proceed by verifying Assumptions \ref{ass:structure} and \ref{ass:paramest}, and then employing Theorem \ref{th:delta} to obtain \eqref{prop:convexdistdisp}. To this end, define the maps $\phi_1 : \mathbb H \rightarrow \mathbb H$ to be given by $\phi_1(\theta) = \theta - \Pi_{\Lambda}\theta$, and $\phi_2 : \mathbb H \rightarrow \mathbf R$ by $\phi_2(\theta) \equiv \|\theta\|_{\mathbb H}$. Letting $\phi \equiv \phi_2 \circ \phi_1$ and noting $\phi_1(\theta_0) = 0$ due to $\theta_0 \in \Lambda$, we then obtain the equality:
\begin{equation}\label{prop:convexdist1}
r_n\|\hat \theta_n - \Pi_{\Lambda} \hat \theta_n\|_{\mathbb H} = r_n\{\phi(\hat \theta_n) - \phi(\theta_0)\} ~.
\end{equation}
By Lemma 4.6 in \cite{zarantonello}, $\phi_1$ is then Hadamard directionally differentiable at $\theta_0$ with derivative $\phi_{1,\theta_0}^\prime : \mathbb H \rightarrow \mathbb H$ given by $\phi_{1,\theta_0}^\prime(h) = h - \Pi_{T_{\theta_0}} h$; see also \cite[p. 135]{shapiro:1994}. Moreover, since $\phi_{2}$ is Hadamard directionally differentiable at $0\in \mathbb H$ with derivative $\phi_{2,0}^\prime(h) = \|h\|_{\mathbb H}$, Proposition 3.6 in \cite{Shapiro1990} implies $\phi$ is Hadamard directionally differentiable at $\theta_0$ with $\phi_{\theta_0}^\prime = \phi_{2,0}^\prime \circ \phi_{1,\theta_0}^\prime$. In particular, we have
\begin{equation}\label{prop:convexdist2}
\phi_{\theta_0}^\prime(h) = \|h - \Pi_{T_{\theta_0}}h\|_{\mathbb H} ~,
\end{equation}
for any $h \in \mathbb H$. Thus, \eqref{prop:convexdist2} verifies Assumptions \ref{ass:structure}(i)-(ii) and, because in this case $\mathbb D = \mathbb D_0 = \mathbb H$, we conclude Assumption \ref{ass:structure}(iii) holds as well. Since Assumption \ref{ass:paramest} was directly imposed, the Proposition then follows form Theorem \ref{th:delta}. \qed

\noindent {\sc Proof of Proposition \ref{prop:convexinf}:} In order to establish the first claim of the Proposition, we first observe that for any $h_1,h_2 \in \mathbb H$ we must have that:
\begin{multline}\label{prop:convexinf1}
\hat \phi_n^\prime(h_1) - \hat \phi_n^\prime(h_2)  \leq  \sup_{\theta \in \Lambda : \|\theta - \Pi_{\Lambda}  \hat \theta_n\|_{\mathbb H} \leq \epsilon_n} \{ \|h_1 - \Pi_{T_\theta} h_1\|_{\mathbb H} - \|h_2 - \Pi_{T_\theta} h_2\|_{\mathbb H}\} \\ \leq \sup_{\theta \in \Lambda : \|\theta - \Pi_{\Lambda}  \hat \theta_n\|_{\mathbb H} \leq \epsilon_n} \{ \|h_1 - \Pi_{T_\theta} h_2\|_{\mathbb H} - \|h_2 - \Pi_{T_\theta} h_2\|_{\mathbb H}\} \leq \|h_1 - h_2\|_{\mathbb H}  ~,
\end{multline}
where the first inequality follows from the definition of $\hat \phi_n^\prime(h)$, the second inequality is implied by $\|h_1 - \Pi_{T_\theta}h_1\|_{\mathbb H} \leq \|h_1 - \Pi_{T_\theta}h_2\|_{\mathbb H}$ for all $\theta  \in \Lambda$, and the third inequality holds by the triangle inequality.  Result \eqref{prop:convexinf1} further implies $\hat \phi_n^\prime(h_2) - \hat \phi_n^\prime(h_1) \leq \|h_1 - h_2\|_{\mathbb H}$, and hence we can conclude $\hat \phi_n^\prime : \mathbb H \rightarrow \mathbf R$ is Lipschitz -- i.e. for any $h_1,h_2 \in \mathbb H$:
\begin{equation}\label{propconvexinf2}
|\hat \phi_n^\prime(h_1) - \hat \phi_n^\prime(h_2)| \leq \|h_1 - h_2\|_{\mathbb H} ~.
\end{equation}
Thus, by Lemma \ref{lm:simpsuff}, in verifying $\hat \phi_n^\prime$ satisfies Assumption \ref{ass:derest} it suffices to show that:
\begin{equation}\label{propconvexinf3}
|\hat \phi_n^\prime(h) - \phi^\prime_{\theta_0}(h)| = o_p(1)
\end{equation}
for all $h \in \mathbb H$. To this end,  note that convexity of $\Lambda$ and Proposition 46.5(2) in \cite{zeidlerIII} imply $\|\Pi_{\Lambda} \theta_0 - \Pi_{\Lambda} \theta\|_{\mathbb H} \leq \|\theta_0 - \theta\|_{\mathbb H}$ for any $\theta \in \mathbb H$. Thus, since $r_n\{\hat \theta_n - \theta_0\}$ is asymptotically tight by Assumption \ref{ass:paramest} and $r_n \epsilon_n \uparrow \infty$ by hypothesis, we conclude that:
\begin{equation}\label{prop:convexinf4}
\liminf_{n\rightarrow \infty} P( \|\Pi_{\Lambda}\theta_0 - \Pi_{\Lambda} \hat \theta_n\|_{\mathbb H} \leq \epsilon_n) \geq \liminf_{n\rightarrow \infty} P( r_n\|\theta_0 - \hat \theta_n\|_{\mathbb H} \leq r_n\epsilon_n) = 1 ~.
\end{equation}
Moreover, the same arguments as in \eqref{prop:convexinf4} and the triangle inequality further imply that:
\begin{multline}\label{prop:convexinf5}
\liminf_{n\rightarrow \infty} P(\|\theta - \Pi_{\Lambda} \theta_0\|_{\mathbb H} \leq 2\epsilon_n \text{ for all } \theta \in \Lambda \text{ s.t. } \|\theta - \Pi_{\Lambda}\hat \theta_n\|_{\mathbb H} \leq \epsilon_n) \\
\geq \liminf_{n\rightarrow \infty} P(\|\Pi_{\Lambda} \theta_0 - \Pi_{\Lambda} \hat \theta_n\|_{\mathbb H} \leq \epsilon_n) = 1 ~.
\end{multline}
Hence, from the definition of $\hat \phi_n^\prime$ and results \eqref{prop:convexinf4} and \eqref{prop:convexinf5} we obtain for any $h \in \mathbb H$:
\begin{equation}\label{prop:convexinf6}
\liminf_{n\rightarrow \infty} P( \|h - \Pi_{T_{\theta_0}} h\|_{\mathbb H} \leq \hat \phi_n^\prime(h) \leq \sup_{\theta \in \Lambda : \|\theta - \Pi_\Lambda \theta_0\|_{\mathbb H}\leq 2\epsilon_n} \|h - \Pi_{T_\theta}h\|_{\mathbb H}) = 1 ~.
\end{equation}
Next, select a sequence $\{\theta_n\}$ with $\theta_n \in \Lambda$ and $\|\theta_n - \Pi_{\Lambda} \theta_0\|_{\mathbb H} \leq 2\epsilon_n$ for all $n$, such that:
\begin{equation}\label{prop:convexinf7}
\limsup_{n\rightarrow \infty} \{\sup_{\theta \in \Lambda : \|\theta - \Pi_\Lambda \theta_0\|_{\mathbb H}\leq 2\epsilon_n} \|h - \Pi_{T_\theta} h\|_{\mathbb H} \} = \lim_{n\rightarrow \infty} \|h - \Pi_{T_{\theta_n}} h\|_{\mathbb H} ~.
\end{equation}
By Theorem 4.2.2 in \cite{aubin:frankowska}, the cone valued map $\theta \mapsto T_\theta$ is lower semicontinuous on $\Lambda$ and hence since $\|\theta_n - \Pi_{\Lambda}\theta_0\|_{\mathbb H} = o(1)$, it follows that there exists a sequence $\{\tilde h_n\}$ such that $\tilde h_n \in T_{\theta_n}$ for all $n$ and $\|\Pi_{T_{\theta_0}} h - \tilde h_n\|_{\mathbb H} = o(1)$. Thus,
\begin{multline}\label{prop:convexinf8}
\limsup_{n\rightarrow \infty} \{\sup_{\theta \in \Lambda : \|\theta - \Pi_\Lambda \theta_0\|_{\mathbb H}\leq 2\epsilon_n} \|h - \Pi_{T_\theta} h\|_{\mathbb H} \} \\
= \lim_{n\rightarrow \infty} \|h - \Pi_{T_{\theta_n}} h\|_{\mathbb H} \leq \lim_{n\rightarrow \infty} \|h - \tilde h_n\|_{\mathbb H} = \|h - \Pi_{T_{\theta_0}}h\|_{\mathbb H} ~,
\end{multline}
where the first equality follows from \eqref{prop:convexinf7}, the inequality by $\tilde h_n \in T_{\theta_n}$, and the second equality by $\|\tilde h_n - \Pi_{T_{\theta_0}} h\|_{\mathbb H} = o(1)$. Hence, combining \eqref{prop:convexinf6} and \eqref{prop:convexinf8} we conclude that \eqref{propconvexinf3} holds, and by Lemma \ref{lm:simpsuff} and \eqref{propconvexinf2} that $\hat \phi_n^\prime$ satisfies Assumption \ref{ass:derest}.

For the second claim, first observe that $\Lambda$ being convex implies $T_{\theta_0}$ is a closed convex cone. Hence, by Proposition 46.5(4) in \cite{zeidlerIII}, it follows that $\|\Pi_{T_{\theta_0}} h\|^2_{\mathbb H} = \langle h, \Pi_{T_{\theta_0}}h\rangle_{\mathbb H}$ for any $h \in \mathbb H$. In particular, for any $h_1,h_2 \in \mathbb H$ we must have:
\begin{equation}\label{prop:convexinf9}
\|h_1 + h_2 - \Pi_{T_{\theta_0}}(h_1 + h_2)\|^2_{\mathbb H} = \langle h_1 + h_2, h_1 + h_2 - \Pi_{T_{\theta_0}}(h_1 + h_2)\rangle_{\mathbb H}  ~.
\end{equation}
However, Proposition 46.5(4) in \cite{zeidlerIII} further implies that $\langle c,h_1+h_2 - \Pi_{T_{\theta_0}} (h_1+h_2)\rangle \leq 0$ for any $h_1,h_2\in \mathbb H$ and $c\in T_{\theta_0}$. Therefore, since $\Pi_{T_{\theta_0}}h_1,\Pi_{T_{\theta_0}} h_2 \in T_{\theta_0}$, we can conclude from result \eqref{prop:convexinf9} and the Cauchy Schwarz inequality that
\begin{multline}\label{prop:convexinf10}
\|h_1 + h_2 - \Pi_{T_{\theta_0}}(h_1 + h_2)\|^2_{\mathbb H}  \leq \langle h_1 - \Pi_{T_{\theta_0}}h_1 + h_2 - \Pi_{T_{\theta_0}}h_2, h_1 + h_2 - \Pi_{T_{\theta_0}}(h_1 + h_2)\rangle_{\mathbb H} \\ \leq \|h_1 + h_2 - \Pi_{T_{\theta_0}}(h_1 + h_2)\|_{\mathbb H} \times \|(h_1 - \Pi_{T_{\theta_0}} h_1) + (h_2 - \Pi_{T_{\theta_0}} h_2)\|_{\mathbb H} ~.
\end{multline}
Thus, the Proposition follows from \eqref{prop:convexinf10} and the triangle inequality. \qed

\newpage

\phantomsection
\bibliography{mybibliography}

\begin{thebibliography}{82}
\expandafter\ifx\csname natexlab\endcsname\relax\def\natexlab#1{#1}\fi
\expandafter\ifx\csname url\endcsname\relax
  \def\url#1{\texttt{#1}}\fi
\expandafter\ifx\csname urlprefix\endcsname\relax\def\urlprefix{URL }\fi
\providecommand{\eprint}[2][]{\url{#2}}

\bibitem[{Abadie et~al.(2002)Abadie, Angrist and
  Imbens}]{abadie:angrist:imbens}
\textsc{Abadie, A.}, \textsc{Angrist, J.} and \textsc{Imbens, G.} (2002).
\newblock Instrumental variables estimates of the effects of subsidized
  trainining on the quantiles of trainee earnings.
\newblock \textit{Econometrica}, \textbf{70} 91--117.

\bibitem[{Aliprantis and Border(2006)}]{aliprantis:border:2006}
\textsc{Aliprantis, C.~D.} and \textsc{Border, K.~C.} (2006).
\newblock \textit{Infinite Dimensional Analysis -- A Hitchhiker's Guide}.
\newblock Springer-Verlag, Berlin.

\bibitem[{Andrews(2000)}]{Andrews2000Bootstrap}
\textsc{Andrews, D. W.~K.} (2000).
\newblock {I}nconsistency of the {B}ootstrap {W}hen a {P}arameter {I}s {O}n the
  {B}oundary of the {P}arameter {S}pace.
\newblock \textit{Econometrica}, \textbf{68} 399--405.

\bibitem[{Andrews and Shi(2013)}]{andrews:shi:2013}
\textsc{Andrews, D. W.~K.} and \textsc{Shi, X.} (2013).
\newblock Inference based on conditional moment inequatlies.
\newblock \textit{Econometrica}, \textbf{81} 609--666.

\bibitem[{Andrews and Soares(2010)}]{andrews:soares:2010}
\textsc{Andrews, D. W.~K.} and \textsc{Soares, G.} (2010).
\newblock Inference for parameters defined by moment inequalities using
  generalized moment selection.
\newblock \textit{Econometrica}, \textbf{78} 119--157.

\bibitem[{Angrist et~al.(2006)Angrist, Chernozhukov and
  Fernandez-Val}]{angrist:Chernozhukov:fernanez}
\textsc{Angrist, J.}, \textsc{Chernozhukov, V.} and \textsc{Fernandez-Val, I.}
  (2006).
\newblock Quantile regression under misspecification, with an application to
  the u.s. wage structure.
\newblock \textit{Econometrica}, \textbf{74} 539--563.

\bibitem[{Arellano et~al.(2012)Arellano, Hansen and
  Sentana}]{arellano:hansen:sentana:2012}
\textsc{Arellano, M.}, \textsc{Hansen, L.~P.} and \textsc{Sentana, E.} (2012).
\newblock Underidentification?
\newblock \textit{Journal of Econometrics}, \textbf{170} 256--280.

\bibitem[{Armstrong and Chan(2012)}]{armstrong:chan:2012}
\textsc{Armstrong, T.~B.} and \textsc{Chan, H.~P.} (2012).
\newblock Multiscale adaptive infernece on conditional moment inequalities.
\newblock Working paper. Yale University.

\bibitem[{Aubin and Frankowska(1990)}]{aubin:frankowska}
\textsc{Aubin, J.-P.} and \textsc{Frankowska, H.} (1990).
\newblock \textit{Set-Valued Analysis}.
\newblock Birkhauser, Boston.

\bibitem[{Beare and Moon(2015)}]{BeareandMoon2012}
\textsc{Beare, B.} and \textsc{Moon, J.} (2015).
\newblock {N}onparametric {T}ests of {D}ensity {R}atio {O}rdering.
\newblock \textit{Econometric Theory}, \textbf{31} 471--492.

\bibitem[{Beare and Shi(2015)}]{BrendanPrelim}
\textsc{Beare, B.} and \textsc{Shi, X.} (2015).
\newblock An imporved bootstrap test of density ratio ordering.
\newblock Working paper, University of California - San Diego.

\bibitem[{Beran(1997)}]{beran:1997}
\textsc{Beran, R.} (1997).
\newblock Diagnosing bootstrap success.
\newblock \textit{Annals of the Institute of Statistical Mathematics},
  \textbf{49} 1--24.

\bibitem[{Beresteanu and Molinari(2008)}]{beresteanu:molinari:2008}
\textsc{Beresteanu, A.} and \textsc{Molinari, F.} (2008).
\newblock Asymptotic properties for a class of partially identified models.
\newblock \textit{Econometrica}, \textbf{76} 763--814.

\bibitem[{Bickel et~al.(1997)Bickel, G\"{o}tze and
  Zwet}]{Bickel_Gotze_Zwet1997}
\textsc{Bickel, P.}, \textsc{G\"{o}tze, F.} and \textsc{Zwet, W.} (1997).
\newblock {R}esampling {F}ewer {T}han $n$ {O}bservations: {G}ains, {L}osses,
  and {R}emedies for {L}osses.
\newblock \textit{Statistica Sinica}, \textbf{7} 1--31.

\bibitem[{Bickel and Freedman(1981)}]{bickel1981some}
\textsc{Bickel, P.~J.} and \textsc{Freedman, D.~A.} (1981).
\newblock Some asymptotic theory for the bootstrap.
\newblock \textit{The Annals of Statistics} 1196--1217.

\bibitem[{Bickel et~al.(1998)Bickel, Ritov and Rydén}]{bickel1998}
\textsc{Bickel, P.~J.}, \textsc{Ritov, Y.} and \textsc{Rydén, T.} (1998).
\newblock Asymptotic normality of the maximum-likelihood estimator for general
  hidden markov models.
\newblock \textit{Ann. Statist.}, \textbf{26} 1614--1635.
\newblock \urlprefix\url{http://dx.doi.org/10.1214/aos/1024691255}.

\bibitem[{Bogachev(1998)}]{bogachevgauss}
\textsc{Bogachev, V.~I.} (1998).
\newblock \textit{Gaussian Measures}.
\newblock American Mathematical Society, Providence.

\bibitem[{Bogachev(2007)}]{bogachev2:2007}
\textsc{Bogachev, V.~I.} (2007).
\newblock \textit{Measure Theory: Volume II}.
\newblock Springer-Verlag, Berlin.

\bibitem[{Bontemps et~al.(2012)Bontemps, Magnac and
  Maurin}]{Bontemps_Magnac_Maurin2012}
\textsc{Bontemps, C.}, \textsc{Magnac, T.} and \textsc{Maurin, E.} (2012).
\newblock {S}et {I}dentified {L}inear {M}odels.
\newblock \textit{Econometrica}, \textbf{80} 1129--1155.

\bibitem[{Buchinsky(1994)}]{buchinsky:1994}
\textsc{Buchinsky, M.} (1994).
\newblock Changes in the u.s. wage structure 1963-1987: Application of quantile
  regression.
\newblock \textit{Econometrica}, \textbf{62} 405--458.

\bibitem[{Bugni(2010)}]{bugni:2010}
\textsc{Bugni, F.~A.} (2010).
\newblock Bootstrap inference in partially identified models defined by moment
  inequalities: Coverage of the identified set.
\newblock \textit{Econometrica}, \textbf{78} 735--753.

\bibitem[{Canay(2010)}]{canay2010}
\textsc{Canay, I.~A.} (2010).
\newblock El inference for partially identified models defined by moment
  inqeualities: Coverage of the identified set.
\newblock \textit{The Journal of Econometrics}, \textbf{156} 408--425.

\bibitem[{Carolan and Tebbs(2005)}]{carolan:tebbs:2005}
\textsc{Carolan, C.~A.} and \textsc{Tebbs, J.~M.} (2005).
\newblock Nonparametric tests for and against likelihood ratio ordering in the
  two-sample problem.
\newblock \textit{Biometrika}, \textbf{92} 159--171.

\bibitem[{Chandrasekarh et~al.(2013)Chandrasekarh, Chernozhukov, Molinari and
  Schrimpf}]{chandrasekhar:chernozhukov:molinari:schrimpf}
\textsc{Chandrasekarh, A.}, \textsc{Chernozhukov, V.}, \textsc{Molinari, F.}
  and \textsc{Schrimpf, P.} (2013).
\newblock Working Paper, Massachusetts Institute of Technology.

\bibitem[{Chen and Fang(2015)}]{chen:fang:2015}
\textsc{Chen, Q.} and \textsc{Fang, Z.} (2015).
\newblock Inference on functionals under first order degeneracy.
\newblock Working paper, University of California - San Diego.

\bibitem[{Chen et~al.(2011)Chen, Tamer and Torgovitsky}]{chen:tamer:torgo:2011}
\textsc{Chen, X.}, \textsc{Tamer, E.} and \textsc{Torgovitsky, A.} (2011).
\newblock Sensitivity analysis in partially identified semiparametric models.
\newblock Cowles Foundation Discussion Paper No. 1836.

\bibitem[{Chernozhukov and Hansen(2005)}]{chernozhukov:hansen:2005}
\textsc{Chernozhukov, V.} and \textsc{Hansen, C.} (2005).
\newblock An iv model of quantile treatment effects.
\newblock \textit{Econometrica}, \textbf{73} 245--261.

\bibitem[{Chernozhukov et~al.(2007)Chernozhukov, Hong and
  Tamer}]{chernozhukov:hong:tamer:2007}
\textsc{Chernozhukov, V.}, \textsc{Hong, H.} and \textsc{Tamer, E.} (2007).
\newblock Estimation and confidence regions for parameter sets in econometric
  models.
\newblock \textit{Econometrica}, \textbf{75} 1243--1284.

\bibitem[{Chernozhukov et~al.(2013)Chernozhukov, Lee and
  Rosen}]{chernozhukov:lee:rosen:2013}
\textsc{Chernozhukov, V.}, \textsc{Lee, S.~S.} and \textsc{Rosen, A.~M.}
  (2013).
\newblock Intersection bounds: Estimation and inference.
\newblock \textit{Econometrica}, \textbf{81} 667--737.

\bibitem[{Chetverikov(2012)}]{chetverikov:2012}
\textsc{Chetverikov, D.} (2012).
\newblock Adaptive test of conditional moment inequalities.
\newblock Working paper. University of California - Los Angeles.

\bibitem[{Ciliberto and Tamer(2009)}]{ciliberto:tamer:2009}
\textsc{Ciliberto, F.} and \textsc{Tamer, E.} (2009).
\newblock Market structure and multiple equilibria in airline markets.
\newblock \textit{Econometrica}, \textbf{77} 1791--1828.

\bibitem[{Davydov et~al.(1998)Davydov, Lifshits and
  Smorodina}]{davydov:lifshits:smorodina:1998}
\textsc{Davydov, Y.~A.}, \textsc{Lifshits, M.~A.} and \textsc{Smorodina, N.~V.}
  (1998).
\newblock \textit{Local Properties of Distribuions of Stochastic Functionals}.
\newblock American Mathematical Society, Providence.

\bibitem[{Dorfman(1938)}]{dorfman1938}
\textsc{Dorfman, R.} (1938).
\newblock A note on the $\delta$-method for finding variance formulae.
\newblock \textit{The Biometric Bulletin}, \textbf{1} 129--137.

\bibitem[{Dufour and Taamouti(2005)}]{dufour2005projection}
\textsc{Dufour, J.-M.} and \textsc{Taamouti, M.} (2005).
\newblock Projection-based statistical inference in linear structural models
  with possibly weak instruments.
\newblock \textit{Econometrica}, \textbf{73} 1351--1365.

\bibitem[{Dugundji(1951)}]{dugundji:1951}
\textsc{Dugundji, J.} (1951).
\newblock An extension of tietze's theorem.
\newblock \textit{Pacific Journal of Mathematics}, \textbf{1} 353--367.

\bibitem[{D\"{u}mbgen(1993)}]{Dumbgen1993}
\textsc{D\"{u}mbgen, L.} (1993).
\newblock {O}n {N}ondifferentiable {F}unctions and the {B}ootstrap.
\newblock \textit{Probability Theory and Related Fields}, \textbf{95} 125--140.

\bibitem[{Efron(1979)}]{efron1979}
\textsc{Efron, B.} (1979).
\newblock Bootstrap methods: Another look at the jacknife.
\newblock \textit{Annals of Statistics}, \textbf{7} 1--26.

\bibitem[{Escanciano and Zhu(2013)}]{escanciano:zhu:2013}
\textsc{Escanciano, J.~C.} and \textsc{Zhu, L.} (2013).
\newblock Inference in semiparametric partially identified models.
\newblock Working paper, Indiana University.

\bibitem[{Fang(2015)}]{fang:2014}
\textsc{Fang, Z.} (2015).
\newblock Optimal plug-in estimators of directionally differentiable
  functionals.
\newblock Working paper, University of California - San Diego.

\bibitem[{Garel and Hallin(1995)}]{garel1995local}
\textsc{Garel, B.} and \textsc{Hallin, M.} (1995).
\newblock Local asymptotic normality of multivariate arma processes with a
  linear trend.
\newblock \textit{Annals of the Institute of Statistical Mathematics},
  \textbf{47} 551--579.

\bibitem[{Gill et~al.(1989)Gill, Wellner and Pr{\ae}stgaard}]{gill1989non}
\textsc{Gill, R.~D.}, \textsc{Wellner, J.~A.} and \textsc{Pr{\ae}stgaard, J.}
  (1989).
\newblock Non-and semi-parametric maximum likelihood estimators and the von
  mises method (part 1)[with discussion and reply].
\newblock \textit{Scandinavian Journal of Statistics} 97--128.

\bibitem[{Hall(1992)}]{HallBoot}
\textsc{Hall, P.} (1992).
\newblock \textit{The Bootstrap and Edgeworth Expansion}.
\newblock Springer-Verlag, Berlin.

\bibitem[{Hansen(2015)}]{hansen2015regression}
\textsc{Hansen, B.~E.} (2015).
\newblock Regression kink with an unknown threshold.
\newblock Working paper, University of Wisconsin - Madison.

\bibitem[{Hansen(2005)}]{hansen:2005}
\textsc{Hansen, P.~R.} (2005).
\newblock A test for superior predictive ability.
\newblock \textit{Journal of Business and Economic Statistics}, \textbf{23}
  365--380.

\bibitem[{Hirano and Porter(2012)}]{hirano:porter:2009}
\textsc{Hirano, K.} and \textsc{Porter, J.~R.} (2012).
\newblock Impossibility results for nondifferentiable functionals.
\newblock \textit{Econometrica}, \textbf{80} 1769--1790.

\bibitem[{Hong and Li(2015)}]{hong:li:2014}
\textsc{Hong, H.} and \textsc{Li, J.} (2015).
\newblock The numerical directional delta method.
\newblock Working paper, Stanford University.

\bibitem[{Horowitz(2001)}]{HorowitzBoot}
\textsc{Horowitz, J.~L.} (2001).
\newblock The bootstrap.
\newblock In \textit{Handbook of Econometrics V} (J.~J. Heckman and E.~Leamer,
  eds.). Elsevier, 3159--3228.

\bibitem[{Imbens and Manski(2004)}]{imbens:manski:2004}
\textsc{Imbens, G.~W.} and \textsc{Manski, C.~F.} (2004).
\newblock Confidence intervals for partially identified parameters.
\newblock \textit{Econometrica}, \textbf{72} 1845--1857.

\bibitem[{Jackwerth(2000)}]{jackwerth:2000}
\textsc{Jackwerth, J.~C.} (2000).
\newblock Recovering risk aversion from option prices and realized returns.
\newblock \textit{Review of Financial Studies}, \textbf{13} 433--451.

\bibitem[{Jha and Wolak(2015)}]{jha2015testing}
\textsc{Jha, A.} and \textsc{Wolak, F.~A.} (2015).
\newblock Testing for market efficiency with transaction costs: An application
  to financial trading in wholesale electricity markets.

\bibitem[{Kaido(2013)}]{Kaido2013dual}
\textsc{Kaido, H.} (2013).
\newblock {A} {D}ual {A}pproach to {I}nference for {P}artially {I}dentified
  {E}conometric {M}odels.
\newblock Working paper, Boston University.

\bibitem[{Kaido and Santos(2014)}]{Kaido_Santos2013}
\textsc{Kaido, H.} and \textsc{Santos, A.} (2014).
\newblock {A}symptotically {E}fficient {E}stimation of {M}odels {D}efined by
  {C}onvex {M}oment {I}nequalities.
\newblock \textit{Econometrica}, \textbf{82} 387--413.

\bibitem[{Kitamura and Stoye(2013)}]{kitamura:stoye:2013}
\textsc{Kitamura, Y.} and \textsc{Stoye, J.} (2013).
\newblock Nonparametric analysis of random utility models: Testing.
\newblock Cemmap Working Paper (CWP36/13). Center for Microdata Methods and
  Practice.

\bibitem[{Kline and Santos(2013)}]{kline:santos:2013}
\textsc{Kline, P.} and \textsc{Santos, A.} (2013).
\newblock Sensitivity to missing data assumptions: Theory and an evaluation of
  the u.s. wage structure.
\newblock \textit{Quantitative Economics}, \textbf{4} 231--267.

\bibitem[{Kosorok(2008)}]{Kosorok2008}
\textsc{Kosorok, M.} (2008).
\newblock \textit{{I}ntroduction to {E}mpirical {P}rocesses and
  {S}emiparametric {I}nference}.
\newblock Springer.

\bibitem[{Ledoux and Talagrand(1991)}]{ledoux:talagrand}
\textsc{Ledoux, M.} and \textsc{Talagrand, M.} (1991).
\newblock \textit{Probability in Banach Spaces}.
\newblock Springer-Verlag, Berlin.

\bibitem[{Lee and Bhattacharya(2015)}]{lee:bhattacharya}
\textsc{Lee, Y.-Y.} and \textsc{Bhattacharya, D.} (2015).
\newblock Welfare analysis for discrete choice with interval-data on income.
\newblock Tech. rep., Working paper, Oxford University.

\bibitem[{Linton et~al.(2010)Linton, Song and Whang}]{Linton2010}
\textsc{Linton, O.}, \textsc{Song, E., K.} and \textsc{Whang, Y.-J.} (2010).
\newblock {A}n {I}mproved {B}ootstrap {T}est of {S}tochastic {D}ominance.
\newblock \textit{Journal of Econometrics}, \textbf{154} 186 -- 202.

\bibitem[{Manski(2003)}]{manski:2003}
\textsc{Manski, C.~F.} (2003).
\newblock \textit{Partial Identification of Probability Distributions}.
\newblock Springer-Verlag, New York.

\bibitem[{Muralidharan and Sundararaman(2011)}]{karthik2011}
\textsc{Muralidharan, K.} and \textsc{Sundararaman, V.} (2011).
\newblock Teacher performance pay: Experimental evidence from india.
\newblock \textit{Journal of Political Economy}, \textbf{9} 39--77.

\bibitem[{Pakes et~al.(2006)Pakes, Porter, Ho and
  Ishii}]{Pakes_Porter_etal2006aWP}
\textsc{Pakes, A.}, \textsc{Porter, J.}, \textsc{Ho, K.} and \textsc{Ishii, J.}
  (2006).
\newblock Moment inequalities and their application.
\newblock Working Paper, Harvard University.

\bibitem[{Politis et~al.(1999)Politis, Romano and Wolf}]{PWR1999}
\textsc{Politis, D.~N.}, \textsc{Romano, J.} and \textsc{Wolf, M.} (1999).
\newblock \textit{{S}ubsampling}.
\newblock Springer, New York.

\bibitem[{Reeds(1976)}]{reeds1976definition}
\textsc{Reeds, J.~A.} (1976).
\newblock \textit{On the definition of von Mises functionals}.
\newblock Thesis, Harvard University.

\bibitem[{Romano(1988)}]{romano1988bootstrap}
\textsc{Romano, J.~P.} (1988).
\newblock A bootstrap revival of some nonparametric distance tests.
\newblock \textit{Journal of the American Statistical Association}, \textbf{83}
  698--708.

\bibitem[{Romano and Shaikh(2008)}]{romano:shaikh:2008}
\textsc{Romano, J.~P.} and \textsc{Shaikh, A.~M.} (2008).
\newblock Inference for identifiable parameters in partially identified
  econometric models.
\newblock \textit{Journal of Statistical Planning and Inference -- Special
  Issue in Honor of Ted Anderson}, \textbf{138} 2786--2807.

\bibitem[{Romano and Shaikh(2010)}]{Romano_Shaikh2010}
\textsc{Romano, J.~P.} and \textsc{Shaikh, A.~M.} (2010).
\newblock {I}nference for the {I}dentified {S}et in {P}artially {I}dentified
  {E}conometric {M}odels.
\newblock \textit{Econometrica}, \textbf{78} 169--211.

\bibitem[{Seo(2014)}]{seo2014tests}
\textsc{Seo, J.} (2014).
\newblock Tests of stochastic monotonicity with improved size and power
  properties.
\newblock Tech. rep., Working paper, University of California - San Diego.

\bibitem[{Shao(1994)}]{shao1994bootstrap}
\textsc{Shao, J.} (1994).
\newblock Bootstrap sample size in nonregular cases.
\newblock \textit{Proceedings of the American Mathematical Society},
  \textbf{122} 1251--1262.

\bibitem[{Shapiro(1990)}]{Shapiro1990}
\textsc{Shapiro, A.} (1990).
\newblock {O}n {C}oncepts of {D}irectional {D}ifferentiability.
\newblock \textit{Journal of Optimization Theory and Applications}, \textbf{66}
  477--487.

\bibitem[{Shapiro(1991)}]{Shapiro1991}
\textsc{Shapiro, A.} (1991).
\newblock {A}symptotic {A}nalysis of {S}tochastic {P}rograms.
\newblock \textit{Annals of Operations Research}, \textbf{30} 169--186.

\bibitem[{Shapiro(1994)}]{shapiro:1994}
\textsc{Shapiro, A.} (1994).
\newblock Existence and differentiability of metric projections in hilbert
  spaces.
\newblock \textit{Siam Journal of Optimization}, \textbf{4} 130--141.

\bibitem[{Song(2014)}]{Song2013Minimax}
\textsc{Song, K.} (2014).
\newblock {L}ocal {A}symptotic {M}inimax {E}stimation of {N}onregular
  {P}arameters with {T}ranslation-{S}cale {E}quivariant {M}aps.
\newblock \textit{Journal of Multivariate Analysis}, \textbf{125} 136--158.

\bibitem[{Strasser(1985)}]{strasser1985mathematical}
\textsc{Strasser, H.} (1985).
\newblock \textit{Mathematical theory of statistics: statistical experiments
  and asymptotic decision theory}, vol.~7.
\newblock Walter de Gruyter.

\bibitem[{van~der Vaart(1991)}]{Vaart1991differentibility}
\textsc{van~der Vaart, A.} (1991).
\newblock {O}n {D}ifferentiable {F}unctionals.
\newblock \textit{The Annals of Statistics}, \textbf{19} pp. 178--204.

\bibitem[{van~der Vaart(1998)}]{Vaart1998}
\textsc{van~der Vaart, A.} (1998).
\newblock \textit{{A}symptotic {S}tatistics}.
\newblock Cambridge University Press.

\bibitem[{van~der Vaart and Wellner(1996)}]{Vaart1996}
\textsc{van~der Vaart, A.} and \textsc{Wellner, J.} (1996).
\newblock \textit{{W}eak {C}onvergence and {E}mpirical {P}rocesses}.
\newblock Springer Verlag.

\bibitem[{ver Hoef(2012)}]{hoef2012}
\textsc{ver Hoef, J.~M.} (2012).
\newblock Who invented the delta method?
\newblock \textit{The American Statistician}, \textbf{66} 124--127.

\bibitem[{White(2000)}]{white:2000}
\textsc{White, H.} (2000).
\newblock A reality check for data snooping.
\newblock \textit{Econometrica}, \textbf{68} 1097--1126.

\bibitem[{Wolak(1988)}]{wolak:1988}
\textsc{Wolak, F.~A.} (1988).
\newblock Duality in testing multivariate hypotheses.
\newblock \textit{Biometrika}, \textbf{75} 611--615.

\bibitem[{Woutersen and Ham(2013)}]{ham:woutersen:2013}
\textsc{Woutersen, T.} and \textsc{Ham, J.~C.} (2013).
\newblock Calculating confidence intervals for continuous and discontinuous
  functions of parameters.
\newblock Cemmap Working Paper (CWP23/13). Center for Microdata Methods and
  Practice.

\bibitem[{Zarantonello(1971)}]{zarantonello}
\textsc{Zarantonello, E.~H.} (1971).
\newblock Projections on convex sets and hilbert spaces and spectral theory.
\newblock In \textit{Contributions to Nonlinear Functional Analysis} (E.~H.
  Zaranotello, ed.). Academic Press.

\bibitem[{Zeidler(1984)}]{zeidlerIII}
\textsc{Zeidler, E.} (1984).
\newblock \textit{Nonlinear Functional Analysis and its Applications III}.
\newblock Springer-Verlag, New York.

\end{thebibliography}
\addcontentsline{toc}{section}{References}

\end{document}